# CICHOŃ'S DIAGRAM AND LOCALISATION CARDINALS

MARTIN GOLDSTERN AND LUKAS DANIEL KLAUSNER

Abstract. We reimplement the creature forcing construction used by Fischer et al. [FGKS17] to separate Cichoń's diagram into five cardinals as a countable support product. Using the fact that it is of countable support, we augment our construction by adding uncountably many additional cardinal characteristics, sometimes referred to as localisation cardinals.

## 1. Introduction

Set theory began as a mathematical subject when Georg Cantor discovered the notion of infinite cardinality and proved that the cardinality of the set of real numbers (the *continuum* $2^{\aleph_0}$) is different from the cardinality of the set of natural numbers $\aleph_0$. The question of "how different?" immediately became a focal point of the new subject and has kept its central place for more than a century. Even before Paul Cohen proved that Cantor's well-known continuum hypothesis cannot be refuted, i. e. that there can consistently be infinite sets of reals of intermediate cardinality, several cardinal numbers of potentially "intermediate" size (so-called cardinal characteristics, such as the unbounding number $\mathfrak{b}$ and the dominating number $\mathfrak{d}$, and of course $\aleph_1$) were known, and the inability of mathematicians to prove equalities between them already hinted at the vast range of unprovability results that emerged as Cohen's forcing method was developed and refined.

For a general overview of cardinal characteristics, see [Bla10] and [Vau90] as well as [BJ95]. Some of the most popular cardinal characteristics are collected in Cichoń's diagram. The paper [FGKS17] is one in a series of progressively more difficult results showing that more and more of the cardinals from Cichoń's diagram can in fact simultaneously be different, in suitably constructed models of set theory. In that particular paper it was shown that those cardinals in Cichoń's diagram which are neither $\mathrm{cov}(\mathcal{N})$ nor provably below $\mathfrak{d}$ (specifically: $\mathrm{non}(\mathcal{M})$, $\mathrm{non}(\mathcal{N})$, $\mathrm{cof}(\mathcal{N})$, and $2^{\aleph_0}$) can have quite arbitrary values (subject to the known inequalities which the diagram expresses).

The older paper [GS93] presented a consistency result about infinitely many pairwise different cardinal characteristics of the continuum with particularly simple definitions, answering a question of Blass [Bla93, p. 78]. Specifically, they constructed a set-theoretic universe where the so-called *localisation cardinals* $\mathfrak{c}_{f,g}$ take

2010 *Mathematics Subject Classification.* Primary 03E17; Secondary 03E35, 03E40.
*Key words and phrases.* cardinal characteristics of the continuum, localisation cardinals, Cichoń's diagram, creature forcing.
Both authors were supported by the Austrian Science Fund (FWF) project P29575 "Forcing Methods: Creatures, Products and Iterations". We are grateful to Diego Alejandro Mejía and the anonymous referee for suggesting numerous improvements to both the content and the presentation of this paper.

uncountably many pairwise different values. (For functions $f, g \in \omega^\omega$, the cardinal $\mathfrak{c}_{f,g}$ is the answer to the question "How many slaloms of width $g$ do we need to cover all functions bounded by $f$?", or formally:

$$\mathfrak{c}_{f,g} := \min \left\{ |\mathcal{S}| \;\middle|\; \mathcal{S} \subseteq \prod_{k<\omega} [f(k)]^{\leq g(k)}, \forall\, x \in \prod_{k<\omega} f(k) \; \exists\, S \in \mathcal{S} \colon x \in^* S \right\}$$

For more detail on $\mathfrak{c}_{f,g}$, see section 10; for a more general treatment of localisation and anti-localisation cardinals, see [KM18, Definition 1.4].)

The common method used in both of these papers is *creature forcing*, the standard reference work for which is by Rosłanowski and Shelah [RS99], but both papers are actually self-contained. While the method of [GS93] was a rather straightforward countable support product of natural tree-like forcing posets, the elements of the forcing poset in [FGKS17] were sequences of so-called *compound creatures*, and the forcing poset was not obviously decomposable as a product of simpler forcing posets. The apparent complexity of that construction may have deterred some readers from taking a closer look at this method.

In the current paper, we will revisit the construction of [FGKS17], but in a more modular way. Using (mostly) a countable support product of lim sup creature forcing posets, together with a lim inf creature forcing poset, we construct a ZFC universe in which the cardinal characteristics $\aleph_1$, $\mathrm{non}(\mathcal{M})$, $\mathrm{non}(\mathcal{N})$, $\mathrm{cof}(\mathcal{N})$ and $2^{\aleph_0}$ are all distinct, and moreover distinct from uncountably many localisation cardinals.

We give a brief outline of the construction. The original forcing construction from [FGKS17] can be decomposed and modified to become a product consisting of four factors:

- a countable support power $\mathbb{Q}_{\mathrm{nn}}^{\kappa_{\mathrm{nn}}}$ of a comparatively simple lim sup creature forcing poset $\mathbb{Q}_{\mathrm{nn}}$, designed to increase $\mathrm{non}(\mathcal{N})$ to $\kappa_{\mathrm{nn}}$;
- a countable support power $\mathbb{Q}_{\mathrm{cn}}^{\kappa_{\mathrm{cn}}}$ of a quite similar creature forcing poset $\mathbb{Q}_{\mathrm{cn}}$, designed to increase $\mathrm{cof}(\mathcal{N})$ to $\kappa_{\mathrm{cn}}$;
- a lim sup creature forcing poset $\mathbb{Q}_{\mathrm{ct}}$ that is not further decomposable for technical reasons, responsible for increasing the continuum to $\kappa_{\mathrm{ct}}$;
- and a lim inf creature forcing poset $\mathbb{Q}_{\mathrm{nm}}$ to increase $\mathrm{non}(\mathcal{M})$ (and hence also $\mathrm{cof}(\mathcal{N}) = \max(\mathrm{non}(\mathcal{M}, \mathfrak{d})) = \mathrm{non}(\mathcal{M})$ to to $\kappa_{\mathrm{nm}}$.

The latter two are still simpler than the parts of the original creature forcing construction corresponding to them; we believe they cannot be replaced by countable support products of creature forcing posets. This new representation allows describing the methods and proofs used in a more modular way, which can then more easily be combined with other lim sup creature forcing posets. As a motivating example, we show how to add a variant of the lim sup creature forcing posets used to separate the localisation cardinals $\mathfrak{c}_{f,g}$ from [GS93] to this construction.

The main result is the following:

**Theorem 1.1.** *Let*

$$\begin{aligned}
\textbf{types} &:= \{\mathrm{nm}, \mathrm{nn}, \mathrm{cn}, \mathrm{ct}\} \cup \bigcup_{\xi<\omega_1} \{\xi\} \\
\textbf{types}_{\lim\sup} &:= \{\mathrm{nn}, \mathrm{cn}, \mathrm{ct}\} \cup \bigcup_{\xi<\omega_1} \{\xi\} \\
\textbf{types}_{\mathrm{modular}} &:= \{\mathrm{nn}, \mathrm{cn}\} \cup \bigcup_{\xi<\omega_1} \{\xi\} \\
\textbf{types}_{\mathrm{slalom}} &:= \bigcup_{\xi<\omega_1} \{\xi\}
\end{aligned}$$

*Assume* CH *in the ground model. Assume we are given cardinals* $\kappa_{\mathrm{nm}} \leq \kappa_{\mathrm{nn}} \leq \kappa_{\mathrm{cn}} \leq \kappa_{\mathrm{ct}}$ *as well as a sequence of cardinals* $\langle \kappa_\xi \mid \xi < \omega_1 \rangle$ *with* $\kappa_\xi \leq \kappa_{\mathrm{cn}}$ *for all* $\xi < \omega_1$ *such that for each* $\mathrm{t} \in \textbf{types}$*,* $\kappa_{\mathrm{t}}^{\aleph_0} = \kappa_{\mathrm{t}}$*. Further assume we are given a congenial sequence*[1] *of function pairs* $\langle f_\xi, g_\xi \mid \xi < \omega_1 \rangle$*.*

*Then there are natural* $\lim\sup$ *creature forcing posets* $\mathbb{Q}_{\mathrm{t}}$ *for each* $\mathrm{t} \in \textbf{types}_{\mathrm{modular}}$*, a* $\lim\sup$ *creature forcing poset* $\mathbb{Q}_{\mathrm{ct},\kappa_{\mathrm{ct}}}$ *and a* $\lim\inf$ *creature forcing poset* $\mathbb{Q}_{\mathrm{nm},\kappa_{\mathrm{nm}}}$ *such that*

$$\mathbb{Q} := \left( \prod_{\mathrm{t} \in \textbf{types}_{\mathrm{modular}}} \mathbb{Q}_{\mathrm{t}}^{\kappa_{\mathrm{t}}} \right) \times \mathbb{Q}_{\mathrm{ct},\kappa_{\mathrm{ct}}} \times \mathbb{Q}_{\mathrm{nm},\kappa_{\mathrm{nm}}}$$

*(where all products and powers have countable support) forces:*

*(M1)* $\mathrm{cov}(\mathcal{N}) = \mathfrak{d} = \aleph_1$*,*
*(M2)* $\mathrm{non}(\mathcal{M}) = \kappa_{\mathrm{nm}}$*,*
*(M3)* $\mathrm{non}(\mathcal{N}) = \kappa_{\mathrm{nn}}$*,*
*(M4)* $\mathrm{cof}(\mathcal{N}) = \kappa_{\mathrm{cn}}$*,*
*(M5)* $\mathfrak{c}_{f_\xi,g_\xi} = \kappa_\xi$ *for all* $\xi < \omega_1$*, and*
*(M6)* $2^{\aleph_0} = \kappa_{\mathrm{ct}}$*.*

*Moreover,* $\mathbb{Q}$ *preserves all cardinals and cofinalities.*

See Figure 1 for a graphical representation of our results.

We give a brief outline of the paper.

- In section 3, we define all the constituent parts of the forcing construction, and in section 4, we show how to put them together and prove a few fundamental properties.
- We then introduce and prove the main properties of the forcing construction which will be used throughout the paper – bigness in section 5 and continuous and rapid reading in section 6 and section 7. The latter section also contains proofs of properness and $\omega^\omega$-bounding, as well as the "easy" parts of the main theorem ((M1) and (M6)).
- The following sections contain the proofs to the remaining parts of the main theorem:
    - section 8 and section 11 prove (M4),

[1] We will define this term in Definition 3.2; informally, it means that for all $\xi \neq \zeta$ we have either $g_\xi < f_\xi \ll g_\zeta < f_\zeta$ or $g_\zeta < f_\zeta \ll g_\xi < f_\xi$.

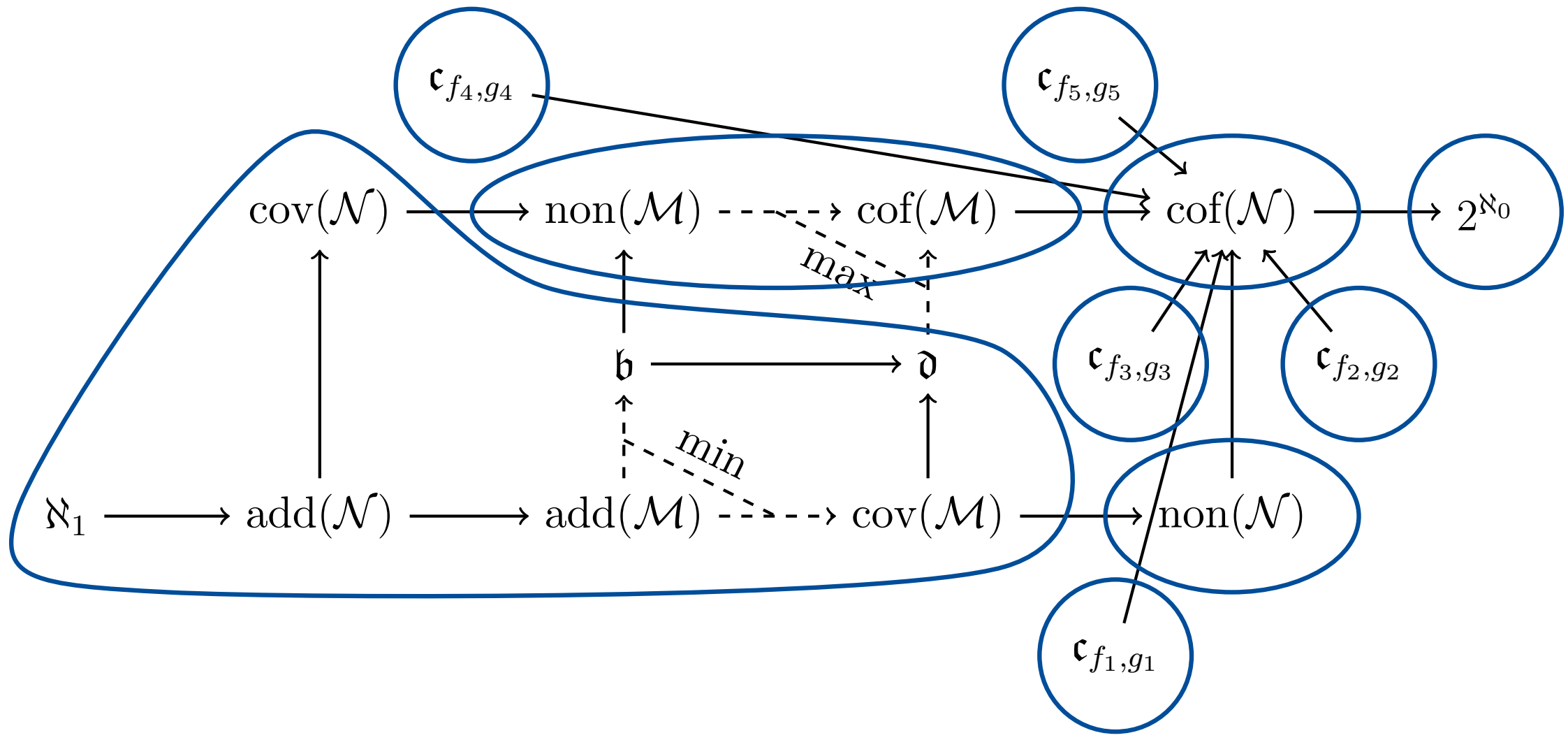

FIGURE 1. Cichoń's diagram with some exemplary $\mathfrak{c}_{f,g}$ added to it; cardinals which are forced to be equal are grouped together, and each such group can be forced to be different from the others subject to the usual constraints.

  - section 9 proves (M2),
  - section 10 proves (M5), and
  - section 11 and section 12 prove (M3).
- Finally, in section 13, we give a brief account of the limitations of the method (and some of our failed attempts to add factors to the construction) and open questions.

## 2. Motivational Preface

We now define the basic framework of the forcing poset. We will not be defining each and every cog of the machinery right from the start; we will instead fill in the blanks one by one, to reduce the complexity and allow for more easily digestible reading.

At the most elementary level, our forcing poset is a product of four parts, each of which employs creature forcing constructions. In such a creature forcing construction, conditions are $\omega$-sequences of so-called creatures, where each creature holds some finite amount of information on the generic real. For technical reasons, we will separate these forcing posets into different sets of *levels* – the (compound) creatures in the lim inf forcing poset $\mathbb{Q}_{\mathrm{nm},\,\kappa_{\mathrm{nm}}}$ will be enumerated by integers of the form $4k$, the creatures in the modular lim sup forcing posets $\mathbb{Q}_{\mathrm{nn}}$ and $\mathbb{Q}_{\mathrm{cn}}$ will be enumerated by integers of the form $4k+1$, the creatures in the modular lim sup forcing posets $\mathbb{Q}_{\xi}$ will be enumerated by integers of the form $4k+2$ and the creatures in the Sacks-like lim sup forcing poset $\mathbb{Q}_{\mathrm{ct},\,\kappa_{\mathrm{ct}}}$ will be enumerated by integers of the form $4k+3$.

We fix an "*index set*" $A$, which is a disjoint union $A_{\mathrm{nm}} \cup A_{\mathrm{nn}} \cup A_{\mathrm{cn}} \cup A_{\mathrm{ct}} \cup \bigcup_{\xi<\omega_1} A_{\xi}$, where each $A_{\mathrm{t}}$ has the appropriate cardinality $\kappa_{\mathrm{t}}$. To each condition $p$ we will associate a countable subset of $A$, called the *support* of $p$.

The modular lim sup forcing posets are not too complicated, having just a creature $C_\ell$ at each level $\ell$ for each index in the support. Each such $C_\ell$ is a subset of some finite set of so-called possibilities $\mathbf{POSS}_{\mathrm{t},\ell}$.

The lim sup forcing poset $\mathbb{Q}_{\mathrm{ct},\kappa_{\mathrm{ct}}}$ cannot be separated into a countable support product of factors. (To be precise, *we* cannot separate it into a countable support product of factors or replace it by such a forcing poset.) We are quite certain that this is due to fundamental structural reasons, namely that in order to prove Lemma 8.1, we have to group the levels (and hence the associated creatures) in this forcing poset together in a certain way, and these partitions need to be compatible, i. e. there must be a single level partition shared by all indices in the support of $\mathbb{Q}_{\mathrm{ct},\kappa_{\mathrm{ct}}}$.

Each element of the lim inf forcing poset $\mathbb{Q}_{\mathrm{nm},\kappa_{\mathrm{nm}}}$ consists of a sequence $(C_\ell)_\ell$, where each $C_\ell$ is a two-dimensional grid of creatures. Each such grid has a finite support $S_\ell \subseteq A_{\mathrm{nm}}$. For each $\ell$, there is a finite set $J_\ell$ (i. e. some natural number), and the grid consists of a $|S_\ell|$-tuple $\langle C_{\ell,\alpha} \mid \alpha \in S_\ell \rangle$ of stacked creatures $C_{\ell,\alpha}$; each stacked creature $C_{\ell,\alpha}$, in turn, is a finite sequence of creatures $C_{(\ell,0),\alpha}, \dots, C_{(\ell,J_\ell-1),\alpha}$. We view this sitation as dividing the level $\ell$ into "*sublevels*" $(\ell,0),\dots,(\ell,J_\ell-1)$. Additionally, each lim inf level $\ell$ also has a so-called "*halving parameter*" $d(\ell)$, a natural number.[2]

For easier reading, we will be using the term "height" to mean "level" for the lim sup forcing posets or "sublevel" for the lim inf forcing poset. A height $L \in$ **heights** is thus either a level $\ell = 4k+1$, $\ell = 4k+2$ or $\ell = 4k+3$ or a sublevel $(\ell,i)$ with $\ell = 4k$ and $i \in J_\ell$. (For $\mathbb{Q}_{\mathrm{ct},\kappa_{\mathrm{ct}}}$, we will consider all creatures within the same class of the level partition as a unit, which complicates induction on the heights a little bit.)

The descriptions of these forcing posets as "lim inf forcing" or "lim sup forcing" refers to the kind of requirements we demand of the sequences of creatures. Each of these creature forcing posets has a norm (a sequence of functions from $2^{\mathbf{POSS}_{\mathrm{t},L}}$ to the non-negative reals) associated with it. As one would expect from the nomenclature, we will demand that for any given condition $p \in \mathbb{Q}$, for each lim sup forcing poset $\mathbb{Q}_{\mathrm{t}}$ we have $\limsup_{\ell\to\infty} \|p(\alpha,\ell)\|_{\mathrm{t},\ell} = \infty$ for each $\alpha \in \operatorname{supp}(p)$ (again, for $\mathbb{Q}_{\mathrm{ct},\kappa_{\mathrm{ct}}}$, this will look a tiny bit different) as well as that $\liminf_{\ell\to\infty} \|p(\ell)\|_{\mathrm{nm}} = \infty$. (Note that in this statement, we are deliberately referring to levels and not to heights, and the limits are to be understood as limits in terms of the tg-appropriate levels.)

The forcing posets involved will depend on certain parameters, which we will define iteratively by induction on the heights. For each height $L$, we will also inductively define natural numbers $n^P_{<L} < n^R_{<L} < n^B_L < n^S_L$.

We want to briefly explain the purpose of these sequences:

- $n^P_{<L}$ will be an upper bound on the *number of possibilities* below $L$ (corresponding to maxposs from [FGKS17] and $n^-$ from [GS93]). By this we mean that $n^P_{<L}$ will bound the number of different possible maximal strengthenings[3] of a condition $p$ below $L$ and hence e. g. the number of iterations we have to go through whenever we want to consider all possible such strengthenings.

[2] See Figure 4 for a graphical representation of this structure.

- $n^B_L$ will be a lower bound on the *bigness* of a creature at height $L$ (corresponding to $b$ from [FGKS17] and also $n^-$ from [GS93]), which we will be defining a bit later. For now, think of this as follows: Whenever we partition a creature $C_L$ into at most $n^B_L$ many sets (e. g. according to which value they force some name to have), there is always one set such that strengthening $C_L$ to a subcreature corresponding to that set will only very slightly decrease the norm.
- $n^S_L$ will be an upper bound on the *size* of $\mathbf{POSS}_{\mathrm{t},L}$ for all $\mathrm{t} \in \mathbf{types}$ (corresponding to $M$ from [FGKS17] and $n^+$ from [GS93]).
- $n^R_{<L}$ will be used to control how quickly a condition $p$ decides finite initial segments of reals (corresponding to $H$ from [FGKS17]), i. e. its *rapidity.* This decision of initial segments will be referred to as "reading" in the sequel.

## 3. Defining the Forcing Factors

Let us now begin to define the framework of the forcing construction.

**Definition 3.1.** Assume we are given cardinals $\kappa_{\mathrm{nm}} \leq \kappa_{\mathrm{nn}} \leq \kappa_{\mathrm{cn}} \leq \kappa_{\mathrm{ct}}$ and a sequence of cardinals $\langle \kappa_\xi \mid \xi < \omega_1 \rangle$ with $\kappa_\xi \leq \kappa_{\mathrm{cn}}$ for all $\xi < \omega_1$ such that for each $\mathrm{t} \in \mathbf{types}$, $\kappa_{\mathrm{t}}^{\aleph_0} = \kappa_{\mathrm{t}}$.[4]

(i) Choose disjoint index sets $A_{\mathrm{t}}$ of size $\kappa_{\mathrm{t}}$ for each $\mathrm{t} \in \mathbf{types}$. We will use the shorthand notations
- $A_{\mathrm{slalom}} := \bigcup_{\xi<\omega_1} A_\xi$,
- $A_{*\mathrm{n}} := A_{\mathrm{cn}} \cup A_{\mathrm{nn}}$, and
- $A := A_{\mathrm{nm}} \cup A_{*\mathrm{n}} \cup A_{\mathrm{slalom}} \cup A_{\mathrm{ct}}$,

as well as the notations
- $\mathbf{typegroups} := \{\mathrm{nm}, *\mathrm{n}, \mathrm{slalom}, \mathrm{ct}\}$, and
- $\mathbf{typegroups}_{\limsup} := \mathbf{typegroups} \setminus \{\mathrm{nm}\}$.

(ii) For each $\ell = 4k$, we will fix some $J_\ell$ with $0 < J_\ell < \omega$. We will refer to the set of *heights*

$$\mathbf{heights} := \bigcup_{k<\omega} (\{(4k,i) \mid i \in J_{4k}\} \cup \{4k+1, 4k+2, 4k+3\})$$

as well as its subsets
- $\mathbf{heights}_{\mathrm{nm}} := \bigcup_{k<\omega} \{(4k,i) \mid i \in J_{4k}\}$,
- $\mathbf{heights}_{*\mathrm{n}} := \{4k+1 \mid k < \omega\}$,
- $\mathbf{heights}_{\mathrm{slalom}} := \{4k+2 \mid k < \omega\}$, and
- $\mathbf{heights}_{\mathrm{ct}} := \{4k+3 \mid k < \omega\}$.

[3] This will be explained in more detail in the following. By "strengthening of $p$ below $L$", we mean conditions $q \leq p$ such that $p$ and $q$ are identical at all heights $K \geq L$, and by "maximal strengthening" we mean that there is no stronger $q$ with this property.

[4] Recall that $\mathbf{types} = \{\mathrm{nm}, \mathrm{nn}, \mathrm{cn}, \mathrm{ct}\} \cup \bigcup_{\xi<\omega_1} \{\xi\}$.

(iii) The heights will be ordered in the obvious way, that is:

$$\ldots < 4k-1 < (4k,0) < (4k,1) < \ldots < (4k, J_{4k}-1)$$
$$< 4k+1 < 4k+2 < 4k+3 < (4k+4,0) < \ldots$$

We will also use $L^+$ and $L^-$ to refer to the successor and predecessor of a height $L$ in this order.

(iv) The creatures of our forcing poset $\mathbb{Q}$ will "live" on (some subset of)

$$A_{\mathrm{nm}} \times \mathbf{heights}_{\mathrm{nm}} \cup A_{*\mathrm{n}} \times \mathbf{heights}_{*\mathrm{n}} \cup A_{\mathrm{slalom}} \times \mathbf{heights}_{\mathrm{slalom}} \cup A_{\mathrm{ct}} \times \mathbf{heights}_{\mathrm{ct}},$$

that is, each $p \in \mathbb{Q}$ will have creatures for each $\alpha$ in a countably infinite $\mathrm{supp}(p) \subseteq A$ and for each height $L$ of the corresponding type (though for each height, only finitely many will be non-trivial). For each

$$(\alpha, L) \in \bigcup_{\mathrm{tg} \in \mathbf{typegroups}} (A_{\mathrm{tg}} \cap \mathrm{supp}(p)) \times \mathbf{heights}_{\mathrm{tg}},$$

there will be a finite set $\mathbf{POSS}_{\alpha,L}$, and the creatures $C_{\alpha,L}$ will be some non-empty subsets of these. (See Figure 2 for a schematic representation of the structure of $\mathbb{Q}$.)

(v) Given some index $\alpha \in A$ respectively some height $L \in \mathbf{heights}$, we will use $\mathrm{tg}(\alpha)$ respectively $\mathrm{tg}(L)$ to denote the appropriate group of types, i. e. the tg such that $\alpha \in A_{\mathrm{tg}}$ respectively the tg such that $L \in \mathbf{heights}_{\mathrm{tg}}$.

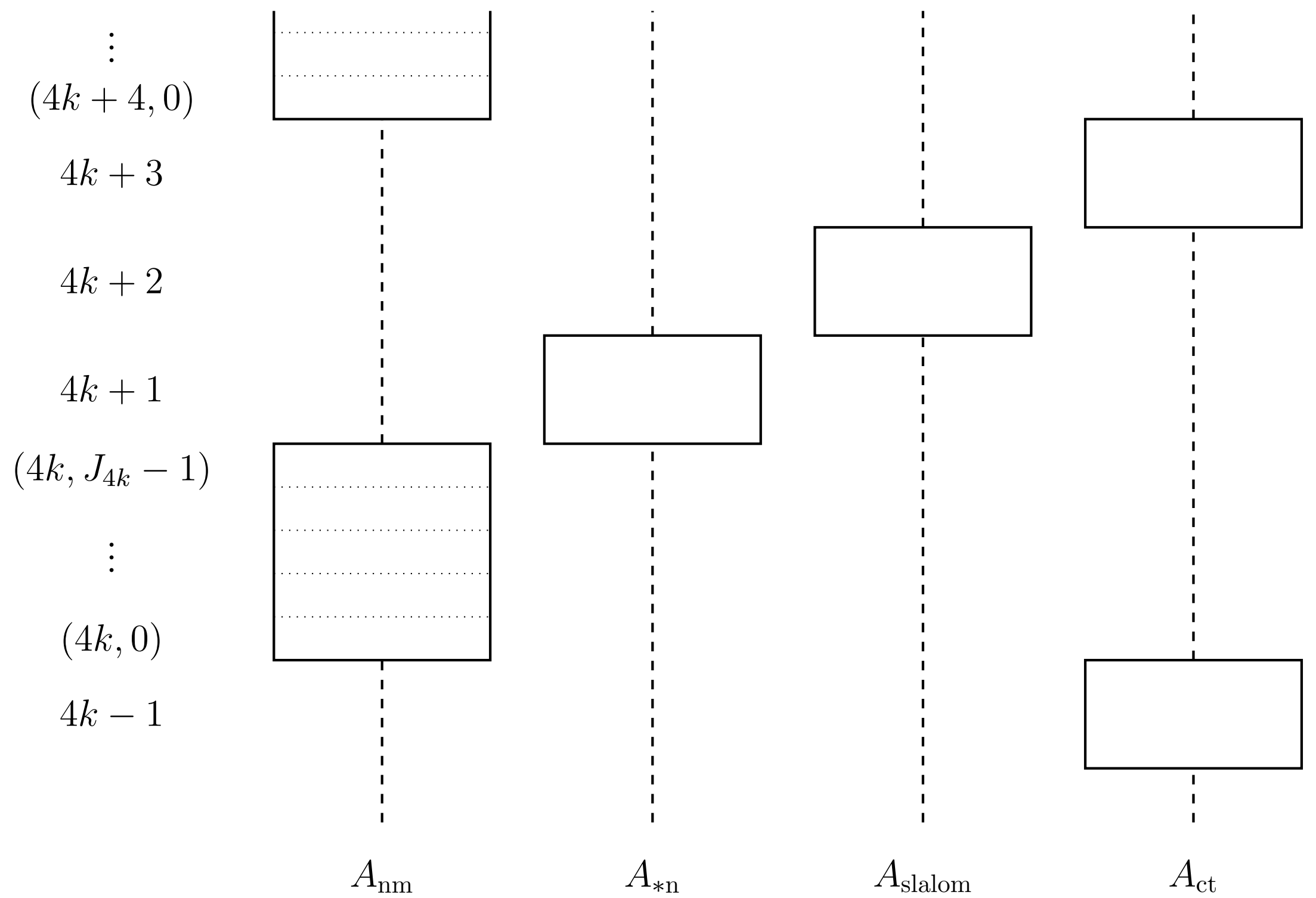

FIGURE 2. A diagram of the basic structure of the forcing poset $\mathbb{Q}$.

We will now first define the forcing posets themselves. However, the inductive definitions of the forcing posets and those of the auxiliary sequences $n^P_{<L}, n^B_L, n^S_L, n^R_{<L}$ mentioned above are actually intertwined. We will be using the auxiliary functions as parameters here and very diligently make sure in section 4 that when inductively defining them, we will not be using anything not previously defined up to that step of the induction process (mostly, this means taking care not to commit off-by-one errors). For now, think of these four sequences as growing very, very quickly and fulfilling $n^P_{<L} \ll n^R_{<L} \ll n^B_L \ll n^S_L \ll n^P_{<L^+}$.[5]

Keep the following in mind: To define creature forcing posets, we mainly have to define the sets of *possibilities* $\mathbf{POSS}_{t,L}$ and the associated norms. The reasons for the specific choices of the norms will only become clear later in section 5, when we define the concept of bigness.

We will start with $\mathbb{Q}_{\text{slalom}}$.

**Definition 3.2.** Given the sequences $n^P_{<L}, n^B_L, n^S_L$, we call a sequence of function pairs $(f_\xi, g_\xi)$ in $\omega^\omega$ *congenial* if:

(i) for each $\xi$ and for all $k < \omega$, $n^B_{4k+2} \leq g_\xi(k) < f_\xi(k) \leq n^S_{4k+2}$,

(ii) for each $\xi$, $\lim_{k\to\infty} \frac{\log f_\xi(k)}{n^B_{4k+2} \cdot \log g_\xi(k)} = \infty$, and

(iii) for all $\xi, \zeta$ with $\xi \neq \zeta$, either $\lim_{k\to\infty} \frac{f_\zeta(k)^2}{g_\xi(k)} = 0$ or $\lim_{k\to\infty} \frac{f_\xi(k)^2}{g_\zeta(k)} = 0$.[6]

When referring to a single pair of functions in a congenial sequence, we will call this a *congenial pair of functions.*

The choice of $n^B_L \ll n^S_L$ will ensure that there are sufficiently many different such function pairs.

**Definition 3.3.** Given the sequences $n^P_{<L}, n^B_L, n^S_L, n^R_{<L}$ and a congenial pair of functions $(f_\xi, g_\xi)$, the forcing factor $\mathbb{Q}_\xi$ is defined as the set of all conditions $p$ fulfilling the following:

(i) $p$ consists of a sequence of creatures $\langle p(L) \mid L \in \mathbf{heights}_{\text{slalom}} \rangle$. Each such $L$ is of the form $4k+2$.

(ii) The possibilities at level $L$ are given by $\mathbf{POSS}_{\xi,L} := f_\xi((L-2)/4) = f_\xi(k)$. This means that for each such $L$, $p(L) = p(4k+2) \subseteq f_\xi(k) = f_\xi((L-2)/4)$ (and $p(L) \neq \varnothing$).

(iii) The norm $\|\cdot\|_{\xi,L}$ on (of course nonempty) subsets of $\mathbf{POSS}_{\xi,L}$ is given by $\|M\|_{\xi,L} := \frac{\log|M|}{n^B_L \cdot \log g_\xi((L-2)/4)} = \frac{\log|M|}{n^B_{4k+2} \cdot \log g_\xi(k)}$.[7]

(iv) There is an increasing sequence of $L_i \in \mathbf{heights}_{\text{slalom}}$ such that $\|p(L_i)\|_{\xi,L_i} \geq i$. Equivalently, $\limsup_{L\to\infty} \|p(L)\|_{\xi,L} = \infty$. This means that for these $L_i$, $|p(L_i)|$ is much larger than $g_\xi((L_i-2)/4)$ (in more legible notation: for

---

[5] As a matter of fact, in the inductive definitions of these sequences, we will only be demanding that they be far larger than some other term, and we define them in some appropriate way to ensure this property; making them even larger would not pose any problems.

[6] Property (iii) here corresponds to the assumption in [GS93, Theorem 3.1], but is more specific. While [GS93] only demands (if we ignore the distinction between $g$ and $h$) that for each $k < \omega$, either $f_\zeta(k)$ is much smaller than $g_\xi(k)$ or $g_\zeta(k)$ is much bigger than $f_\xi(k)$ (but the relations could be different for different $k$'s), we actually demand that the functions are eventually ordered the same way. We could just as well work with the more general property, but we believe our restriction makes the proofs somewhat easier to digest.

[7] As stated previously, the specific definition of the norm here is not really important, and other definitions might work equally well; we mainly require the norms to have a property called "bigness" (defined in section 5), which will be proved in Theorem 5.6.

these $k_i$, i. e. such that $L_i =: 4k_i + 2$, we have that $|p(4k_i + 2)|$ is much larger than $g_\xi(k_i)$).

A condition $q$ is stronger than a condition $p$ if $q(L) \subseteq p(L)$ holds for each $L \in \textbf{heights}_{\text{slalom}}$.

Note that Definition 3.2 (ii) ensures that $\mathbb{Q}_\xi$ is non-empty.

A generic filter for $\mathbb{Q}_\xi$ will add a new real $y_\xi \in \prod_n f_\xi(n)$. The purpose of the specific definition of $\mathbb{Q}_\xi$ is to ensure that $y_\xi$ will infinitely often avoid any slalom of size $g_\xi$ from the ground model (i. e. for any $\langle S(n) \mid n \in \omega \rangle \in V$ with $|S(n)| \leq g_\xi(n)$, there will be infinitely many $n$ such that $y_\xi(n) \notin S(n)$). We will show that $y_\xi$ will moreover avoid all $g$-slaloms added by factors other than $\mathbb{Q}_\xi$.

Next, we define $\mathbb{Q}_{\text{nn}}$.

**Definition 3.4.** Given the sequences $n^P_{<L}, n^B_L, n^S_L, n^R_{<L}$, the forcing factor $\mathbb{Q}_{\text{nn}}$ is defined as the set of all conditions $p$ fulfilling the following:

(i) For each $L \in \textbf{heights}_{*\text{n}}$, we fix a finite interval $I_L \subseteq \omega$ (for notational simplicity, disjoint from all $I_K$ for $K < L$) such that with the definitions given below, $\| \textbf{POSS}_{\text{nn},L} \|_{\text{nn},L} > n^B_L$.

(ii) $p$ consists of a sequence of creatures $\langle p(L) \mid L \in \textbf{heights}_{*\text{n}} \rangle$.

(iii) The possibilities at level $L$ are given by

$$\textbf{POSS}_{\text{nn},L} := \left\{ X \subseteq 2^{I_L} \;\middle|\; |X| = \left(1 - \frac{1}{2^{n^B_L}}\right) \cdot |2^{I_L}| \right\},$$

that is, all subsets $X$ of $2^{I_L}$ of relative size $1 - 2^{-n^B_L}$.

(iv) The norm $\| \cdot \|_{\text{nn},L}$ on subsets of $\textbf{POSS}_{\text{nn},L}$ is given by

$$\|M\|_{\text{nn},L} := \frac{\log \|M\|^{\text{intersect}}_L}{n^B_L \log n^B_L}$$

with $\|M\|^{\text{intersect}}_L := \min\{|Y| \mid Y \subseteq 2^{I_L}, \forall X \in M \colon X \cap Y \neq \varnothing\}$.

(v) There is an increasing sequence of $L_i \in \textbf{heights}_{*\text{n}}$ such that $\|p(L_i)\|_{\text{nn},L_i} \geq i$. Equivalently, $\limsup_{L\to\infty} \|p(L)\|_{\text{nn},L} = \infty$.

Assuming $|I_L| \geq n^B_L$, the minimum in the definition above is equal to $2^{|I_L| - n^B_L} + 1$ for $M = \textbf{POSS}_{\text{nn},L}$. Therefore, fulfilling $\| \textbf{POSS}_{\text{nn},L} \|_{\text{nn},L} > n^B_L$ (and the lim sup condition on the norms for conditions) is achievable by choosing $I_L$ sufficiently large, and hence $\mathbb{Q}_{\text{nn}}$ is non-empty.

A condition $q$ is stronger than a condition $p$ if $q(L) \subseteq p(L)$ holds for each $L \in \textbf{heights}_{*\text{n}}$.

A $\mathbb{Q}_{\text{nn}}$-generic filter will define a sequence $\langle X_L \mid L \in \textbf{heights}_{*\text{n}} \rangle$ of quite large sets $X_L \subseteq 2^{I_L}$, which we can identify with corresponding sets $X'_L \subseteq 2^\omega$ of measure $1 - 2^{-n^B_L}$. The set $\bigcap_K \bigcup_{L>K} (2^\omega \setminus X'_L)$ will then be a null set. We will ensure in section 11 that this null set covers all reals from the ground model (and more).

Next, we define $\mathbb{Q}_{\text{cn}}$. The norm we give here is technically different from (and hopefully simpler than) the one given in [FGKS17], but fulfils the same purpose.

**Definition 3.5.** Given the sequences $n^P_{<L}, n^B_L, n^S_L, n^R_{<L}$, the forcing factor $\mathbb{Q}_{\text{cn}}$ is defined as the set of all conditions $p$ fulfilling the following:

(i) For each $L \in \mathbf{heights}_{*\mathrm{n}}$, we fix a finite interval $I_L \subseteq \omega$ (for notational simplicity, disjoint from all $I_K$ for $K < L$) such that with the definitions given below, $\| \mathbf{POSS}_{\mathrm{cn},L} \|_{\mathrm{cn},L} > n^B_L$.

(ii) $p$ consists of a sequence of creatures $\langle p(L) \mid L \in \mathbf{heights}_{*\mathrm{n}} \rangle$.

(iii) The possibilities at level $L$ are again given by

$$\mathbf{POSS}_{\mathrm{cn},L} := \left\{ X \subseteq 2^{I_L} \;\middle|\; |X| = \left(1 - \frac{1}{2^{n^B_L}}\right) \cdot |2^{I_L}| \right\},$$

that is, all subsets $X$ of $2^{I_L}$ of relative size $1 - 2^{-n^B_L}$. This is the same kind of possibility set as for $\mathbb{Q}_{\mathrm{nn}}$, but the norm is different:

(iv) The norm $\| \cdot \|_{\mathrm{cn},L}$ on subsets of $\mathbf{POSS}_{\mathrm{cn},L}$ is given by

$$\|M\|_{\mathrm{cn},L} := \frac{\log |M| - \log \binom{2^{|I_L|-1}}{2^{n^B_L}-1}}{2^{\min I_L} \cdot (n^B_L)^2 \cdot \log 3n^B_L}.$$

(v) There is an increasing sequence of $L_i \in \mathbf{heights}_{*\mathrm{n}}$ such that $\|p(L_i)\|_{\mathrm{cn},L_i} \geq i$. Equivalently, $\limsup_{L\to\infty} \|p(L)\|_{\mathrm{cn},L} = \infty$.

A condition $q$ is stronger than a condition $p$ if $q(L) \subseteq p(L)$ holds for each $L \in \mathbf{heights}_{*\mathrm{n}}$.

Note that if the $I_L$ are chosen as above, then $\mathbb{Q}_{\mathrm{cn}}$ is non-empty; see the observation below on why such a choice of $I_L$ is possible.

Similar to $\mathbb{Q}_{\mathrm{nn}}$ above, the generic object for $\mathbb{Q}_{\mathrm{cn}}$ will give us a generic null set; in this case, we will ensure that this null set is not covered by any null set from $V$.

This is the only forcing poset which we have substantially modified as compared to [FGKS17], so let us briefly explain what we have done and why that is fine. (We will omit the rounding to integers in the following calculations.)

**Observation 3.6.** The construction in [FGKS17] combines two different norms which provide properties required for the proofs, $\mathrm{nor}^{\cap}_b$ and $\mathrm{nor}^{\div}_{I,b}$. One can easily see that $\mathrm{nor}^{\cap}_b(x) = \left\lfloor \frac{\log x}{\log 3b} \right\rfloor$; this is not explicitly stated in [FGKS17], but is is straightforward from the definitions (setting $M(\delta, \ell) := 3\,{}^{\ell}/_{\delta}$).

On the other hand, $\mathrm{nor}^{\div}_{I,b}(x) = {}^{x}/\binom{2^{|I|-1}}{2^{b}-1}$. Consider $\mathrm{nor}^{\div}_{I,b}(\mathbf{POSS}_{\mathrm{cn},L})$ for $b := n^B_L$; for appropriately large (with respect to $n^B_L$) choices of $I$, this can become arbitrarily large. But then the same holds for $\log \mathrm{nor}^{\div}_{I,b}(x)$ and also for $\frac{\log \mathrm{nor}^{\div}_{I,b}(x)}{\log 3b}$, and we have only decreased the norm by modifying it this way. This norm now is almost the same as $\mathrm{nor}^{\cap}_b$ except for the subtrahend, but we have already established that this norm still goes to infinity. Hence, if we replace the norm in [FGKS17, Definition 10.1.1 (3)] (which defines $\mathbb{Q}_{\mathrm{cn}}$) by this instead, all relevant properties are preserved and we have used a slightly nicer, closed form instead.

Next, we define $\mathbb{Q}_{\mathrm{ct},\kappa_{\mathrm{ct}}}$. As mentioned before, this forcing poset is a $\limsup$ forcing poset, but not decomposable into factors.

**Lemma 3.7.** *There is a monotone function $f\colon (\omega \setminus \{0\})^3 \to \omega$ with the following properties: For all $j, n, c \in (\omega \setminus \{0\})^3$, whenever a product $B_1 \times \cdots \times B_j$ (where each $B_i$ has cardinality $\geq f(j, n, c)$) is coloured with $c$ colours, there are subsets*

*$A_1 \subseteq B_1$, ..., $A_j \subseteq B_j$ such that the set $A_1 \times \cdots \times A_j$ is homogeneous for this colouring and $|A_i| \geq n$ for $i = 1, \ldots, j$.*

*Proof.* Given $j, n, c \in (\omega \setminus \{0\})^3$, define $b_1 := c \cdot n$, inductively define $b_{i+1} := c^{b_1 \cdots b_i} \cdot n$ for $i = 1, \ldots, j-1$, and let $f(j, n, c) := b_j$. This guarantees that any map from a set of size $b_{i+1}$ into a set of size $c^{b_1 \cdots b_i}$ will be constant on a set of size $\geq n$.

Now let $C\colon B_1 \times \cdots \times B_j \to c$ be a colouring; without loss of generality, assume all $B_i$ have the same size $b_j = f(j, n, c)$. For $i = 1, \ldots, j$, let $B'_i \subseteq B_i$ be a set of cardinality $b_i$.

Write $C_j$ for the map $C{\restriction}_{B'_1 \times \cdots \times B'_j}$, which we now view as a map from $B'_j$ into $c^{B'_1 \times \cdots \times B'_{j-1}}$. This map is constant on a set $A_j \subseteq B'_j$ of size $\geq n$, say with value $C_{j-1}$. $C_{j-1}$ is a map from $B'_1 \times \cdots \times B'_{j-1}$ into $c$, but we view it as a map from $B'_{j-1}$ into $c^{B'_1 \times \cdots \times B'_{j-2}}$; as above, we find $A_{j-1} \subseteq B'_{j-1}$ of size $\geq n$ such that $C_{j-1}$ is constant on $A_{j-1}$ with value $C_{j-2}$.

We continue by induction until we reach $C_0$, which is a map from the empty product into $c$, with some value $c^*$. Then $C$ is constant on $A_1 \times \cdots \times A_j$ with value $c^*$. □

**Definition 3.8.** We define a norm $\mathrm{nor}^{B,m}_{\mathrm{Sacks}}(X)$ as

$$\mathrm{nor}^{B,m}_{\mathrm{Sacks}}(X) := \max\left(\{i \mid F^B_m(i) \leq |X|\} \cup \{0\}\right),$$

where the function $F^B_m(i)$ is defined as an "iterate" of the function $f$ from the previous lemma as follows: $F^B_m(0) := 1$ and $F^B_m(n+1) := f(m, F^B_m(n), B)$ for all $n \geq 0$.

**Observation 3.9.** For fixed $B$ and $m$, we can make $\mathrm{nor}^{B,m}_{\mathrm{Sacks}}(X)$ as large as we want if we allow arbitrarily large $X$.

This norm satisfies a version[8] of *bigness*:

**Lemma 3.10.** *Assume that $X_1, \ldots, X_j$ are sets satisfying $\mathrm{nor}^{B,j}_{\mathrm{Sacks}}(X_i) \geq k+1$ for $i = 1, \ldots, j$ and that $C\colon X_1 \times \cdots \times X_j \to \{1, \ldots, B\}$ is a colouring. Then there are subsets $Y_i \subseteq X_i$ for $i = 1, \ldots, j$ of norm $\geq k$ such that $C{\restriction}_{Y_1 \times \cdots \times Y_j}$ is constant.*

*Proof.* Our assumption implies $|X_i| \geq F^B_j(k+1)$, so $|X_i| \geq f(j, F^B_j(k), B)$ for all $i$. By the characteristic property of $f$ from Lemma 3.7, we can find subsets $Y_i \subseteq X_i$ such that $C$ is constant on $Y_1 \times \cdots \times Y_j$, and all $Y_i$ have size $\geq F^B_j(k)$, hence norm $\geq k$. □

**Definition 3.11.** Given a cardinal $\kappa_{\mathrm{ct}}$ with $\kappa_{\mathrm{ct}}^{\omega} = \kappa_{\mathrm{ct}}$, an index set $A_{\mathrm{ct}}$ of size $\kappa_{\mathrm{ct}}$ and the sequences $n^P_{<L}, n^B_L, n^S_L, n^R_{<L}$, the forcing poset $\mathbb{Q}_{\mathrm{ct},\kappa_{\mathrm{ct}}}$ is defined as the set of all conditions $p$ with countable $\mathrm{supp}(p) \subseteq A_{\mathrm{ct}}$ fulfilling the following:

(i) There is a partition of $\mathbf{heights}_{\mathrm{ct}}$ into a sequence of consecutive intervals, which we will call a *frame*. To avoid confusion with the intervals $I_L$, we will refer to the intervals (i. e. partition classes) of the frame as *segments*.

[8] This is a simplified variant of [FGKS17, 2.3.6(6)]

(ii) We formalise the frame as a function segm: $\textbf{heights}_{\mathrm{ct}} \to \textbf{heights}_{\mathrm{ct}}^{<\omega}$ mapping each height to the finite tuple of heights constituting the segment it belongs to. Using floor$(L)$ to refer to $\min(\mathrm{segm}(L))$, we then have

$$\mathrm{segm}(L) = [\mathrm{floor}(L), \mathrm{floor}(L^*) - 4] \cap \textbf{heights}_{\mathrm{ct}},$$

where $L^*$ is the minimal $L' \in \textbf{heights}_{\mathrm{ct}}$ above $L$ such that $\mathrm{segm}(L') \neq \mathrm{segm}(L)$. (See Figure 3 for the structure of a frame.)

(iii) For each $L \in \textbf{heights}_{\mathrm{ct}}$, we fix a finite interval $I_L \subseteq \omega$ (for notational simplicity, disjoint from all $I_K$ for $K < L$) such that with the definitions given below, $\| \textbf{POSS}_{\mathrm{ct},L} \|_{\mathrm{ct},L} > n^B_L$. This is possible by Observation 3.9; we just have to choose a sufficiently long interval $I_L$. (This ensures that even for the trivial frame consisting of only singleton segments, there are valid conditions.)

(iv) For each $\alpha \in \mathrm{supp}(p)$, $p(\alpha)$ consists of a sequence of creatures $\langle p(\alpha, L) \mid L \in \textbf{heights}_{\mathrm{ct}} \rangle$.

(v) Given a segment $\bar{M} := \langle M_1, \dots, M_m \rangle$, we will use the abbreviated notation $p(\alpha, \bar{M})$ to denote $\langle p(\alpha, M_1), \dots, p(\alpha, M_m) \rangle$. We will call $p(\alpha, \bar{M})$ a *creature segment*.

(vi) The possibilities at level $L$ are given by $\textbf{POSS}_{\mathrm{ct},L} := 2^{I_L}$. This means that for each such $L$ and $\alpha \in \mathrm{supp}(p)$, $p(\alpha, L) \subseteq 2^{I_L}$ (and $p(L) \neq \varnothing$).

(vii) We will treat each creature segment as a unit and define the norm of a condition in $\mathbb{Q}_{\mathrm{ct},\kappa_{\mathrm{ct}}}$ on creature segments. Let $\bar{X} := \langle X_1, \dots, X_m \rangle$ be a creature segment of $p(\alpha)$ (for some $\alpha \in \mathrm{supp}(p)$) associated with the segment $\bar{K} := \langle K_1, \dots, K_m \rangle$. This means that for some $i < \omega$, we have $K_j = 4(i + j) + 3$ and $X_j \subseteq \textbf{POSS}_{\mathrm{ct},K_j}$ for $j \in \{1, \dots, m\}$.

(viii) The norm $\| \cdot \|_{\mathrm{ct},L}$ on a creature segment $\bar{X}$ is given by

$$\|\bar{X}\|_{\mathrm{ct},K_1} := \max_{j \in \{1,\dots,m\}} \mathrm{nor}_{\mathrm{Sacks}}^{n^B_{K_1},k}(X_j),$$

where $k$ is such that $K_1 = 4k + 3$.

(ix) For each $\alpha \in \mathrm{supp}(p)$, there is an increasing sequence of $L_i \in \textbf{heights}_{\mathrm{ct}}$ (each of which is a segment's initial height) such that $\|p(\alpha, \mathrm{segm}(L_i))\|_{\mathrm{ct},L_i} \geq i$. Equivalently, $\limsup_{L \to \infty} \|p(\alpha, \mathrm{segm}(L))\|_{\mathrm{ct},\mathrm{floor}(L)} = \infty$.

A condition $q$ is stronger than a condition $p$ if $\mathrm{supp}(q) \supseteq \mathrm{supp}(p)$, $q(\alpha, L) \subseteq p(\alpha, L)$ holds for each $\alpha \in \mathrm{supp}(p)$ and each $L \in \textbf{heights}_{\mathrm{ct}}$, and the frame of $q$ is coarser than the frame of $p$.

Note that the choice of the $I_L$ above ensures that $\mathbb{Q}_{\mathrm{ct},\kappa_{\mathrm{ct}}}$ is non-empty. Also note that we will sometimes for brevity write $\|X\|_{\mathrm{ct},K_1}$ to mean $\mathrm{nor}_{\mathrm{Sacks}}^{n^B_{K_1},k}(X)$, to avoid having to single out the ct case when it is not strictly necessary.

We remark that if $q \leq p$ only differs from $p$ in that its frame is coarser, then the norms on the creature segments in $q$ are greater or equal to the norms on the corresponding creature segments in $p$.

**Observation 3.12.** One could decompose the forcing poset $\mathbb{Q}_{\mathrm{ct},\kappa_{\mathrm{ct}}}$ as the composition of a forcing poset $\mathbb{F}$ defining a frame partition $\mathcal{F}$ and a parametrised version $\mathbb{Q}^{\mathcal{F}}_{\mathrm{ct},\kappa_{\mathrm{ct}}}$ with a fixed frame, somewhat analogous to the well-known decomposition of Mathias forcing $\mathbb{R}$ into a forcing poset $\mathbb{U}$ adding an ultrafilter $\mathcal{U}$ and the parametrised Mathias forcing $\mathbb{R}_{\mathcal{U}}$ (cf. [Hal17, Lemma 26.10]). However, this decomposition of $\mathbb{Q}_{\mathrm{ct},\,\kappa_{\mathrm{ct}}}$ neither simplifies nor generalises our constructions, so we will not use it.

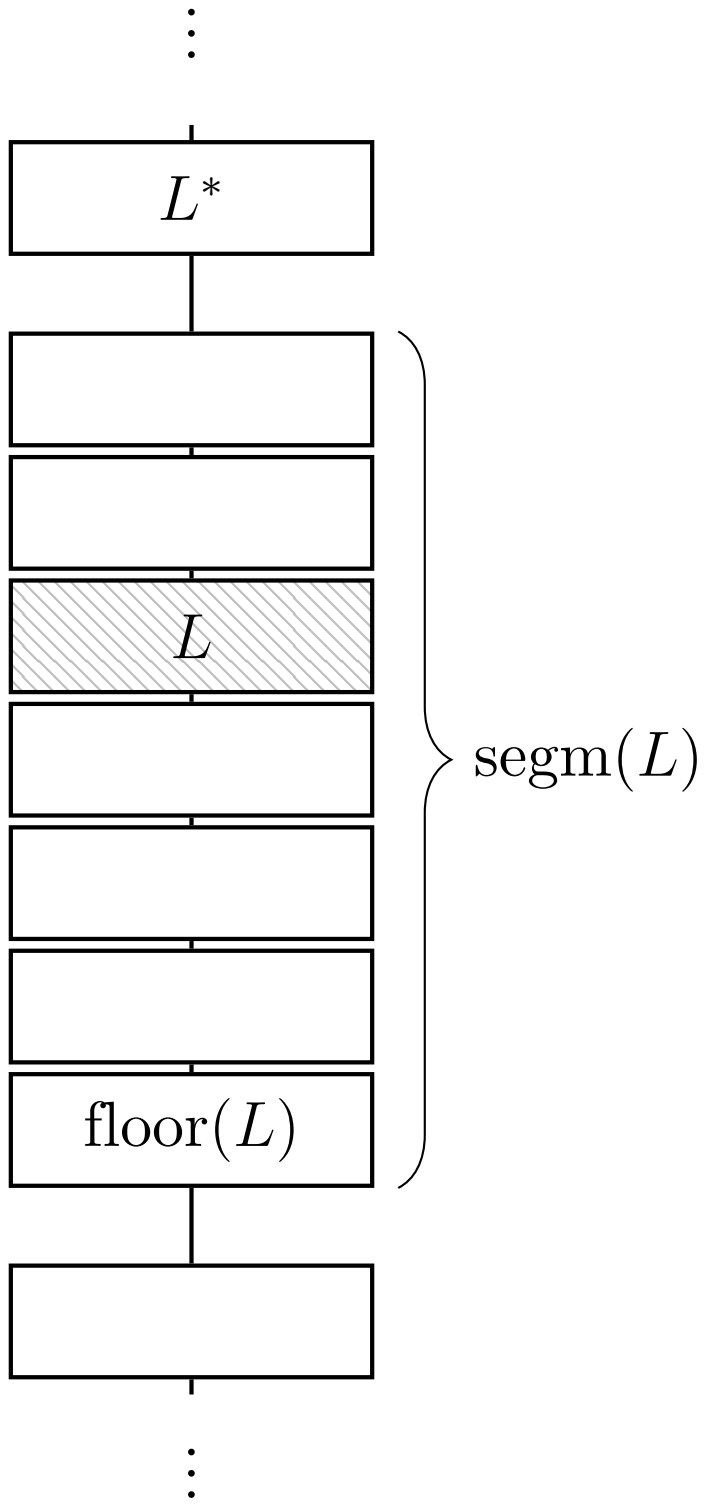

FIGURE 3. A visualisation of the segment segm($L$) of the height $L$ in a frame, as well as the last and first heights of the preceding and succeeding segment, respectively.

Finally, we define $\mathbb{Q}_{\mathrm{nm},\,\kappa_{\mathrm{nm}}}$ (the only lim inf forcing poset), which we will define *en bloc* instead of as a countable support product.

**Definition 3.13.** Given a cardinal $\kappa_{\mathrm{nm}}$ with $\kappa_{\mathrm{nm}}^{\omega} = \kappa_{\mathrm{nm}}$, an index set $A_{\mathrm{nm}}$ of size $\kappa_{\mathrm{nm}}$ and the sequences $n^P_{<L}, n^B_L, n^S_L, n^R_{<L}$, the forcing poset $\mathbb{Q}_{\mathrm{nm},\,\kappa_{\mathrm{nm}}}$ is defined as the set of all conditions $p$ fulfilling the following:

(i) We first define the finite sets of sublevels to be

$$J_{4k} := 3^{(4k+1)\cdot 2^{4k\cdot n^P_{<(4k,0)}}}.$$

(The reason for this definition will become clear in (xi) below.)

(ii) $p$ consists of a sequence $\langle p(4k) \mid k < \omega\rangle$ of *compound creatures*, each of which has a finite support $S_{4k} \subseteq \mathrm{supp}(p) \subseteq A_{\mathrm{nm}}$, together with a sequence of reals $d(4k)$ which are called *halving parameters*. The supports $S_{4k}$ are non-decreasing and $\bigcup_{k<\omega} S_{4k} = \mathrm{supp}(p)$.

(iii) For each $\alpha \in S_{4k}$, $C_{\alpha,4k}$ is a *stacked creature*, i. e. a finite sequence consisting of $|J_{4k}|$ many creatures $C_{\alpha,(4k,i)}$, $i \in J_{4k}$. So the compound creature $p(4k)$ is indexed by $S_{4k} \times \{(4k, i) \mid i \in J_{4k}\}$, and each $p(\alpha, 4k) = C_{\alpha,4k}$ is a stacked creature. (See Figure 4 for an example of a compound creature.)

(iv) In the following, $L = (4k, i)$ will refer to some sublevel height of $p(4k)$.

(v) Define the *cell norm* $\|\cdot\|^{\text{cell}}_L$ by $\|M\|^{\text{cell}}_L := \frac{\log|M|}{n^B_L \log n^B_L}$.

(vi) For each $L = (4k, i)$, fix a finite interval $I_L \subseteq \omega$ (for notational simplicity, disjoint from all $I_K$ for $K < L$) such that $\|2^{I_L}\|^{\text{cell}}_L > n^B_L$.

(vii) The possibilities at height $L$ are given by $\mathbf{POSS}_{\text{nm},L} := 2^{I_L}$; this means that for each such $L$ and all $\alpha \in S_{4k}$, $p(\alpha, L) \subseteq 2^{I_L}$.

(viii) Call the minimal $4k < \omega$ such that there is an $\alpha \in \text{supp}(p)$ and a $K = (4k, i) \in \mathbf{heights}_{\text{nm}}$ with $|p(\alpha, K)| > 1$ the *trunk length* of $p$, denoted by $\text{trklgth}(p)$. We call the part of $p$ below $\text{trklgth}(p)$ the *trunk* and denote it by $\text{trunk}(p)$; the trunk of $p$ consists of singletons $p(\alpha, L)$ in $\mathbf{POSS}_{\text{nm},L}$ for each $\alpha \in \text{supp}(p)$ and each $L = (4j, i)$ with $j < k$ and $i \in J_{4j}$. However, there may be (and, on a dense set of conditions, there will be) many additional singleton creatures above $\text{trklgth}(p)$, namely all creatures $p(\alpha, (4\ell, j))$ with $\alpha \in \text{supp}(p) \setminus S_{4\ell}$ and $(4\ell, j) > \text{trklgth}(p)$, which we do not consider part of the trunk. By definition, we let $S_{4j} = \varnothing$ for $j < k$, and to avoid unnecessary complications, we demand that the halving parameters of $p$ have to be 0 below $\text{trklgth}(p)$.

(ix) Each sublevel fulfils a condition called *modesty*,[9] which means that for each $L = (4k, i)$, there is at most one index $\alpha \in S_{4k}$ such that $p(\alpha, L)$ is non-trivial, i. e. $|p(\alpha, L)| > 1$.

(x) Define the *stack norm* $\|\cdot\|^{\text{stack}}_{4k}$ on stacked creatures by the following: $\|p(\alpha, 4k)\|^{\text{stack}}_{4k}$ is the maximal $r$ such that there is an $X \subseteq J_{4k}$ with $\mu_{4k}(X) := \frac{\log_3|X|}{4k+1} \geq r$ and such that $\|p(\alpha, (4k, x))\|^{\text{cell}}_{(4k,x)} \geq r$ for all $x \in X$.

Note that $\mu_{4k}(J_{4k}) = 2^{4k \cdot n^P_{<(4k,0)}}$. We will later choose $n^B_L$ such that $n^B_L > 2^{4k \cdot n^P_{<(4k,0)}}$ for $L > (4k, 0)$, so the stack norm of a maximal stacked creature having the full $2^{I_L}$ at each height then also is $2^{4k \cdot n^P_{<(4k,0)}}$, .

(xi) Define $\|\cdot\|_{\text{nm},4k}$ on compound creatures by

$$\|p(4k)\|_{\text{nm},4k} := \frac{\log_2\left(\min\{\|p(\alpha, 4k)\|^{\text{stack}}_{4k} \mid \alpha \in S_{4k}\} - d(4k)\right)}{n^P_{<(4k,0)}}$$

or 0 if the above is ill-defined (for instance, because the minimal stacked creature norm is smaller than the halving parameter $d(4k)$).

Note that for the trunk, applying this norm to any subset of $\text{supp}(p)$ and any level $4j < \text{trklgth}(p)$ also just yields 0. Also note that the norm of the maximal compound creature consisting of the maximal stacked creatures thus is

$$\frac{\log_2(2^{4k \cdot n^P_{<(4k,0)}} - d(4k))}{n^P_{<(4k,0)}},$$

which for $d(4k) = 0$ is exactly $4k$. (However, such a maximal compound creature is not modest, of course.)

[9] In [FGKS17], the set of modest conditions was introduced as a dense subset of the conditions instead; while we will do this similarly in the following section in Lemma 4.4, for sake of easier presentation, we prefer to define the lim inf conditions as modest right from the start. Note that if we drop modesty from the definition, applying Lemma 2.2.2 from [FGKS17] to an arbitrary condition $p$ easily yields a stronger modest condition $q$.

(xii) There is an increasing sequence of $k_i < \omega$ such that $\|p(4\ell)\|_{\mathrm{nm},4\ell} \geq i$ for all $\ell \geq k_i$. Equivalently, $\liminf_{k\to\infty} \|p(4k)\|_{\mathrm{nm},4k} = \infty$.

(xiii) The relative widths of the compound creatures converge to 0, i. e. $\lim_{k\to\infty} \frac{|S_{4k}|}{4k+1} = 0$.

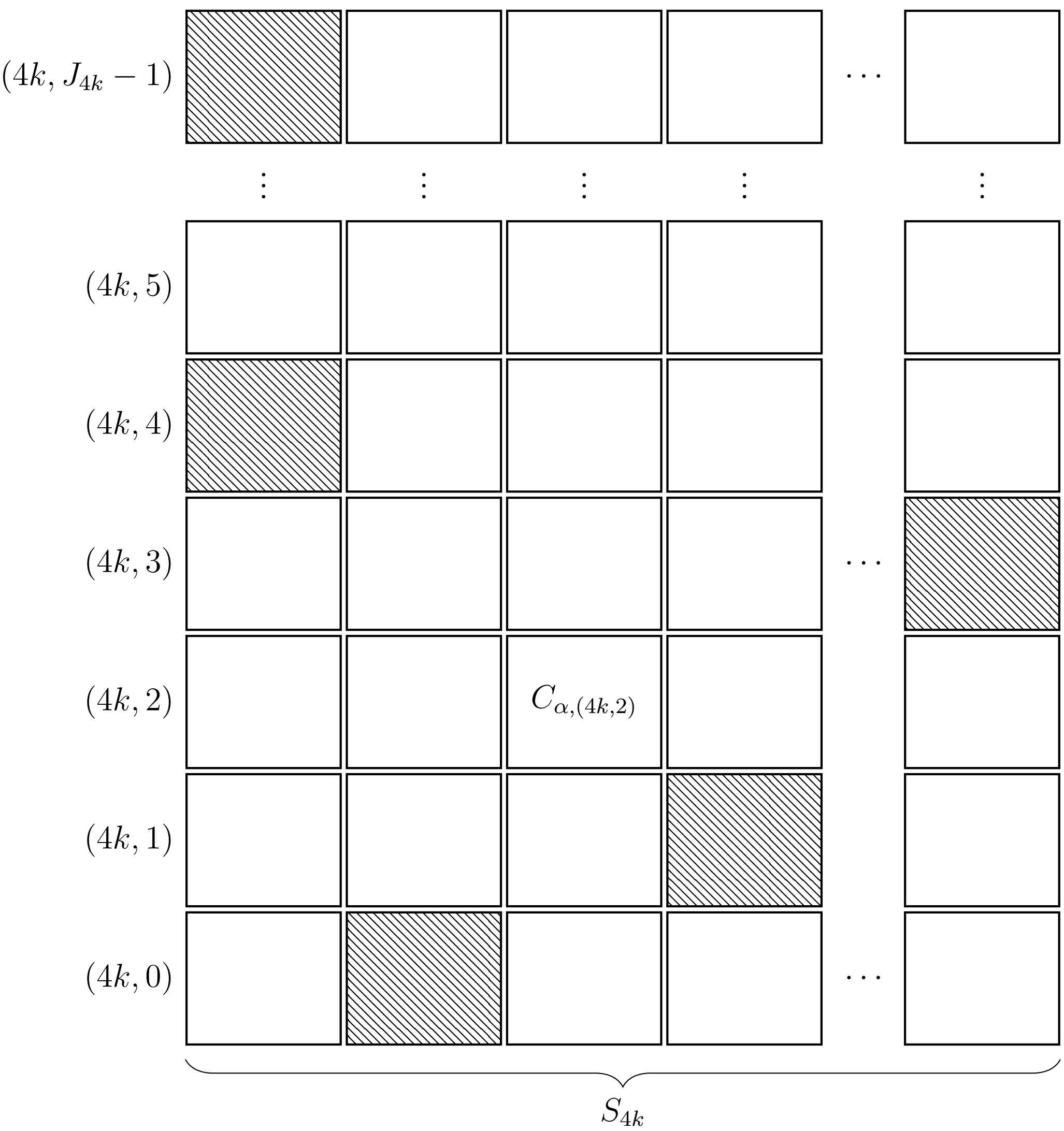

FIGURE 4. An example of a compound creature $C := p(4k)$ of a condition $p \in \mathbb{Q}_{\mathrm{nm},\kappa_{\mathrm{nm}}}$. A possible pattern of cells containing non-trivial creatures is hatched.

A condition $q$ is stronger than a condition $p$ if

- $\mathrm{trklgth}(q) \geq \mathrm{trklgth}(p)$ (the trunk may grow),
- $S_{4k}(q) \supseteq S_{4k}(p)$ for each $4k \geq \mathrm{trklgth}(q)$ (above the trunk, the supports do not shrink),
- for each $k < \omega$, for each $\alpha \in S_{4k}(p)$ and for each $i \in J_{4k}$, $q(\alpha, (4k, i)) \subseteq p(\alpha, (4k, i))$, and
- $d(q)(4k) \geq d(p)(4k)$ (the halving parameters do not decrease).

Note that for reasonably small halving parameters (namely, such that for some $k_0 < \omega$ and some $\varepsilon > 0$

$$d(4k) < 2^{4k \cdot n^P_{<(4k,0)}} \cdot (1 - \varepsilon)$$

holds for all $k > k_0$), the choice of the $I_L$ above ensures that there are conditions in $\mathbb{Q}_{\mathrm{nm},\kappa_{\mathrm{nm}}}$ with such a sequence of small but positive halving parameters. (This is, of course, trivially true for a halving parameter sequence of all 0s.)

We want to briefly remark on the terminology: Our compound creatures are the smallest possible kind of compound creatures in [FGKS17], since there compound creatures could span multiple levels. Our cells and stacks are the subatoms and atoms of [FGKS17].

## 4. Putting the Parts Together

We remark that we still have not shown that the definitions we make are possible, as we require the sequences $n^P_{<L}, n^B_L, n^S_L, n^R_{<L}$ to make the definitions. Before we rectify that omission, we define the full forcing poset.

We fix an "*index set*" $A$, which is a disjoint union $A_{\mathrm{nm}} \cup A_{\mathrm{nn}} \cup A_{\mathrm{cn}} \cup A_{\mathrm{ct}} \cup \bigcup_{\xi<\omega_1} A_\xi$, where each $A_{\mathrm{t}}$ has the appropriate cardinality $\kappa_{\mathrm{t}}$.

Our forcing poset will be a product of five factors. Three of these factors will themselves be countable support products of simpler forcings, using the index sets $A_{\mathrm{nn}}$, $A_{\mathrm{cn}}$ and $\bigcup_{\xi<\omega_1} A_\xi$, respectively; in the representation below, we will collect them into one larger product.

The other two factors will not literally be products, but will (when ignoring the norms) superficially resemble products, using the index sets $A_{\mathrm{nm}}$ and $A_{\mathrm{ct}}$, respectively.

**Definition 4.1.** Recall $\mathbf{typegroups} = \{\mathrm{nm}, {*}\mathrm{n}, \mathrm{slalom}, \mathrm{ct}\}$, $\mathbf{typegroups}_{\limsup} = \mathbf{typegroups} \setminus \{\mathrm{nm}\}$ from Definition 3.1 (i). Let

$$\begin{aligned} \mathbf{types} &:= \{\mathrm{nm}, \mathrm{nn}, \mathrm{cn}, \mathrm{ct}\} \cup \bigcup_{\xi<\omega_1} \{\xi\} \\ \mathbf{types}_{\limsup} &:= \{\qquad \mathrm{nn}, \mathrm{cn}, \mathrm{ct}\} \cup \bigcup_{\xi<\omega_1} \{\xi\} \\ \mathbf{types}_{\mathrm{modular}} &:= \{\qquad \mathrm{nn}, \mathrm{cn} \qquad\} \cup \bigcup_{\xi<\omega_1} \{\xi\} \\ \mathbf{types}_{\mathrm{slalom}} &:= \bigcup_{\xi<\omega_1} \{\xi\} \end{aligned}$$

Assume we are given cardinals $\kappa_{\mathrm{nm}} \leq \kappa_{\mathrm{nn}} \leq \kappa_{\mathrm{cn}} \leq \kappa_{\mathrm{ct}}$ as well as a sequence of cardinals $\langle \kappa_\xi \mid \xi < \omega_1 \rangle$ with $\kappa_\xi \leq \kappa_{\mathrm{cn}}$ for all $\xi < \omega_1$ such that for each $\mathrm{t} \in \mathbf{types}$, $\kappa_{\mathrm{t}}^{\aleph_0} = \kappa_{\mathrm{t}}$. Then our forcing poset is defined as follows:

$$\mathbb{Q} := \prod_{\mathrm{t} \in \mathbf{types}_{\mathrm{modular}}} \mathbb{Q}_{\mathrm{t}}^{\kappa_{\mathrm{t}}} \times \mathbb{Q}_{\mathrm{ct},\kappa_{\mathrm{ct}}} \times \mathbb{Q}_{\mathrm{nm},\kappa_{\mathrm{nm}}},$$

where all products and powers have countable support.

Since $\mathbb{Q}$ is a product, a condition $q$ is stronger than a condition $p$ if each factor of $q$ is stronger than the corresponding factor of $p$. See Fact 4.10 for a detailed description of all the properties subsumed by the statement "$q \leq p$".

For any tg $\in$ **typegroups** and any t $\in$ **types**, we will use $p(\mathrm{tg})$ and $p(\mathrm{t})$ to refer to the restriction of $p$ to $A_{\mathrm{tg}}$ and $\mathrm{A}_{\mathrm{t}}$, respectively.

**Definition 4.2.** Given $p \in \mathbb{Q}$, tg $\in$ $\textbf{typegroups}_{\limsup}$ and $L \in \textbf{heights}_{\mathrm{tg}}$, define $\mathrm{supp}(p, \mathrm{tg}, L)$ to be the set of all $\alpha \in A_{\mathrm{tg}}$ such that for some $K \leq L$ in $\textbf{heights}_{\mathrm{tg}}$, $|p(\alpha, K)| > 1$. This means that the tg-specific support of a condition at some height $L$ is the set of all indices of that group of types such that $p$ has already had a non-trivial creature at that index up to $L$.

For tg $=$ ct, we usually will refer to the support of a segment $\bar{K} = \langle K_1, \dots, K_m \rangle$ (since we treat each creature segment as a whole) and mean $\mathrm{supp}(p, \mathrm{ct}, \bar{K}) = \mathrm{supp}(p, \mathrm{ct}, K_m)$.

For tg $=$ nm, we define $\mathrm{supp}(p, \mathrm{nm}, 4k) := S_{4k}(p)$ and $\mathrm{supp}(p, \mathrm{nm}, (4k, i)) := \mathrm{supp}(p, \mathrm{nm}, 4k)$ for all $i \in J_{4k}$.

We define $\mathrm{supp}(p, L)$ to be the union of all appropriate $\mathrm{supp}(p, \mathrm{tg}, K)$ with $K \leq L$, and $\mathrm{supp}(p)$ to be the union of all $\mathrm{supp}(p, L)$ with $L \in$ **heights**.

Finally, we will use the notation $p(\mathrm{tg})$ to refer to the tg part of a condition $p$ as well as the notation $p(\mathrm{nm}, 4k)$, $p(\mathrm{slalom}, 4k+2)$ and $p(\mathrm{ct}, \bar{L})$ to refer to the respective heights and segments of the corresponding part, respectively.

We immediately remark that we will instead work with a dense subset of $\mathbb{Q}$:

**Definition 4.3.** We call a condition $p \in \mathbb{Q}$ *modest* if

(i) for each tg $\in$ $\textbf{typegroups}_{\limsup}$, $\mathrm{supp}(p, \mathrm{tg}, \ell) = \varnothing$ for all $\ell < \mathrm{trklgth}(p)$,[10]
(ii) for each $L \in$ **heights**, there is at most one index $\alpha \in A$ such that $p(L, \alpha)$ is non-trivial, i. e. $|p(L, \alpha)| > 1$,
(iii) the segments of $p(\mathrm{ct})$ are such that for each segment $\bar{L} = \langle L_1, \dots, L_m \rangle$ with $L_1 = 4k+3$, for all $\alpha \in \mathrm{supp}(p, \mathrm{ct}, \bar{L})$ we have $\|p(\alpha, \bar{L})\|_{\mathrm{ct}, L_1} \geq k$ as well as $|\mathrm{supp}(p, \mathrm{ct}, \bar{L})| = |\mathrm{supp}(p, \mathrm{ct}, L_m)| < k$, and
(iv) for each segment $\bar{L} = \langle L_1, \dots, L_m \rangle$ of the frame of $p(\mathrm{ct})$ (with $L_1 = 4k+3$) and $\alpha \in A_{\mathrm{ct}}$ such that $p(\alpha, \bar{L})$ is non-trivial, there is exactly one $L^* \in \bar{L}$ such that $p(\alpha, L^*)$ is non-trivial, and furthermore $\|p(\alpha, \bar{L})\|_{\mathrm{ct}, L_1} = k$. We furthermore demand that the size of each such $p(\alpha, L^*)$ is already minimal; in particular, this means that there are exactly $c = F_k^{n^B_{L_1}}(k)$ many possibilities in $p(\alpha, L^*)$.

**Lemma 4.4.** *The set of modest conditions is dense in $\mathbb{Q}$; moreover, for any $p \in \mathbb{Q}$ there is even a modest $q \leq p$ with the same support.*

*Proof.* Given an arbitrary $p \in \mathbb{Q}$, we have to find a modest $q \leq p$. We first pick arbitrary singletons in each non-trivial creature below $\mathrm{trklgth}(p)$ to fulfil (i). Then, we define $q$ piecewise for each tg $\in$ **typegroups**:

[10] In the preceding section (in Definition 3.13), $\mathrm{trklgth}(p)$ only took the nm part of $p$ into account, but now we want to make sure that $p$ has no non-trivial creatures below $\mathrm{trklgth}(p)$ at all.

- For $\mathrm{tg} = \mathrm{nm}$, we have already defined the compound creatures such that they fulfil (ii).
- For $\mathrm{tg} \in \mathbf{typegroups}_{\mathrm{modular}}$, finding $q(\mathrm{tg})$ is just a matter of diagonalisation and bookkeeping (picking arbitrary singletons within creatures as required to fulfil (ii)).
- To achieve (iii), we coarsen the frame to encompass sufficiently large $p(\alpha, K)$ into the creature segments and/or strengthen to arbitrary singletons whenever necessary (plus bookkeeping, again).
- Property (iv) is fulfilled by choosing, for each $\alpha \in \mathrm{supp}(p, \bar{L})$, a single $L^* \in \bar{L}$ such that $p(\alpha, L^*)$ has norm at least $k$, then shrinking $p(\alpha, L^*)$ such that its cardinality is minimal with still the same norm, and shrinking all other $p(\alpha, L')$ to arbitrary singletons; by definition, all of this leaves the ct norms of such segments at least $k$ and the resulting $q(\mathrm{ct})$ is still a valid (part of a) condition.

It is clear that $\mathrm{supp}(q) = \mathrm{supp}(p)$. □

Note that for any modest $p \in \mathbb{Q}$, property (ii) immediately implies that $\mathrm{supp}(p, L)$ is finite for any $L \in \mathbf{heights}$.

We will extend the meaning of the word "trunk" to refer to the entire single possibility of a modest condition $p$ below the trunk length of $p$.

We will only ever work with modest conditions; whenever we speak of conditions, the qualifier "modest" is implied. Though the results of some constructions may not be modest conditions themselves, we can find stronger conditions with the same support by the preceding lemma; and moreover, if a condition is already partially modest (i. e. modest up to a certain height), we can keep that part when making it modest.

We remark that modesty properties (iii) and (iv) roughly correspond to the concept of "Sacks pruning" in [FGKS17, subsection 3.4] and [FGKS17, Lemma 2.3.6].

Modesty properties (ii)–(iv) are of vital importance to the entire construction. Without them, it would not be possible to define the sequences $n^P_{<L}, n^B_L, n^S_L, n^R_{<L}$ in a sensible manner, which we are now finally able to do. Before we do so, we have to introduce the "maximal strengthenings of a condition $p$ below a height $L$" mentioned in the introductory remarks.

**Definition 4.5.** Given a condition $p \in \mathbb{Q}$, we call a height $L$ *relevant* if either

- $L \in \mathbf{heights}_{\mathrm{ct}}$ is the minimum $L_1$ of a segment $\bar{L} = \langle L_1, \dots, L_m \rangle$ of the frame of $p(\mathrm{ct})$; or
- $L \in \mathbf{heights} \setminus \mathbf{heights}_{\mathrm{ct}}$ and there is an $\alpha \in \mathrm{supp}(p)$ such that $p(\alpha, L)$ is non-trivial, i. e. such that $|p(\alpha, L)| > 1$.

We will use this terminology to simplify the structure of proofs in which we iterate over the heights and modify a condition at each height; naturally, we will only need to do this at the relevant heights. (The reason why it is the minimum height of a segment that we consider "relevant" is that in an inductive construction over all heights, the minimum height of a compound creature in the ct part is the inductive step in which it has to be dealt with.)

**Definition 4.6.** We define the *possibilities of a condition* $p \in \mathbb{Q}$ up to some height $L$ as follows:

- Recall $\mathrm{tg}(\alpha)$ and $\mathrm{tg}(K)$ as defined in Definition 3.1 (v). For each sensible choice of $\alpha \in \mathrm{supp}(p)$ and $K \in \mathbf{heights}$ (i. e. such that $\mathrm{tg}(\alpha) = \mathrm{tg}(K)$), let $\mathrm{poss}(p, \alpha, K) := p(\alpha, K)$.
- For each $\alpha \in \mathrm{supp}(p) \smallsetminus A_{\mathrm{ct}}$ and each $L \in \mathbf{heights}$, let
$$\mathrm{poss}(p, \alpha, {<}L) := \prod_{\substack{K<L \\ \mathrm{tg}(K)=\mathrm{tg}(\alpha)}} \mathrm{poss}(p, \alpha, K).$$
- For each $\alpha \in \mathrm{supp}(p) \cap A_{\mathrm{ct}}$ and each segment $\bar{L} = \langle L_1, \dots, L_m \rangle$ of the frame of $p(\mathrm{ct})$, let
$$\mathrm{poss}(p, \alpha, {<}L_1) := \prod_{K<L_1} \mathrm{poss}(p, \alpha, K),$$
$$\mathrm{poss}(p, \alpha, {<}L_i) := \prod_{K \leq L_m} \mathrm{poss}(p, \alpha, K)$$
for all $i \in \{2, \dots, m\}$. This means that when talking about possibilities of $p(\mathrm{ct})$ below some $L_i$, we have to take the whole segment of the frame into account unless we are at the lower boundary $L_1$ of such a segment.
- (For easier notation, consider $\mathrm{poss}(p, \alpha, {<}L)$ for $L \in \mathbf{heights} \smallsetminus \mathbf{heights}_{\mathrm{ct}}$ and $\alpha \in \mathrm{supp}(p) \cap A_{\mathrm{ct}}$ to mean $\mathrm{poss}(p, \alpha, {<}L^*)$ with $L^* := \min\{K \in \mathbf{heights}_{\mathrm{ct}} \mid L < K\}$.)
- For each $L \in \mathbf{heights}$, let
$$\mathrm{poss}(p, {<}L) := \prod_{\alpha \in \mathrm{supp}(p) \smallsetminus A_{\mathrm{ct}}} \mathrm{poss}(p, \alpha, {<}L) \times \prod_{\alpha \in \mathrm{supp}(p) \cap A_{\mathrm{ct}}} \mathrm{poss}(p, \alpha, {<}L).$$

Note that while both factors in the product above are technically infinite products themselves, thanks to modesty actually only finitely (even boundedly) many of the factors will be non-trivial. The fact that for each $p$ and $L$ iterating over all $\eta \in \mathrm{poss}(p, {<}L)$ only takes boundedly many steps (with the bound depending only on $L$) will be very important in many of the following proofs. Also note that for $L \leq \mathrm{trklgth}(p)$, $|\mathrm{poss}(p, {<}L)| = 1$.[11]

**Definition 4.7.** Given $p \in \mathbb{Q}$, $L \in \mathbf{heights}$ and $\eta \in \mathrm{poss}(p, {<}L)$, we define $p \wedge \eta =: q$ as the condition resulting from replacing all creatures below $L$ as well as those above $L$ in the current segment of the frame of $p(\mathrm{ct})$ with the singletons from $\eta$. Formally, $q$ is defined by

- $q(\alpha, K) := \{\eta(\alpha, K)\}$ for all $K < L$ and $q(\alpha, M) := p(\alpha, K)$ for all $K \geq L$ and all $\alpha \in \mathrm{supp}(p) \smallsetminus A_{\mathrm{ct}}$, and
- $q(\alpha, K) := \{\eta(\alpha, K)\}$ for all $K < L^*$ and $q(\alpha, K) := p(\alpha, K)$ for all $K \geq L^*$ and all $\alpha \in \mathrm{supp}(p) \cap A_{\mathrm{ct}}$, where $L^* := \min\{M \in \mathbf{heights}_{\mathrm{ct}} \mid M \geq L$ and $M$ is the minimum of a segment of $p(\mathrm{ct})\}$.

---

[11] A further note: While as in [FGKS17], the "shapes" of possibilities are not really "nice", this is less of a conceptual problem in this paper, as due to the compartmentalisation of the creatures to different heights depending on the factor they belong to, possibilities are by necessity tiered; the fact that the ct possibilities may be further "down" in the heights structure is less of a conceptual stretch now.

In some proofs, we will use the notation $p^{<L}$ or $q^{\geq L}$ to denote partial initial or terminal (pseudo-)conditions in the obvious sense of $p^{<L} := \langle p(K) \mid K < L\rangle$ and $q^{\geq L} := \langle q(M) \mid M \geq L\rangle$. We will denote the join of such partial conditions by $p^{<L} {}^\frown q^{\geq L}$; we will at those times take special care to make sure what we are writing down actually ends up being a proper condition.

We will now finally show that the definition of the sequences $n^P_{<L}, n^B_L, n^S_L, n^R_{<L}$ is possible in a consistent way. What we are actually doing is the following: We define the base sets in each level/height $L$ of the forcing posets $\mathbb{Q}_{\mathrm{t}}$, $\mathrm{t} \in \mathbf{types}_{\limsup}$, respectively in each sublevel/height $L$ of $\mathbb{Q}_{\mathrm{nm},\,\kappa_{\mathrm{nm}}}$ iteratively by induction on the levels and also define the four sequences for that $L$ in that step, assuming we already know the four sequences for $K < L$. The order of definitions is as follows:

(1) $n^P_{<L}$,
(2) $n^R_{<L}$,
(3) $n^B_L$,
(4) the section of the forcing poset for the height $L$, and finally
(5) $n^S_L$.

**Definition 4.8.** Recall that $L^-$ and $L^+$ denote the predecessor and successor of a height $L$, respectively. We define the sequences $n^P_{<L}, n^R_{<L}, n^B_L, n^S_L$ as follows:

$\boldsymbol{n^P_{<L}}$: We recall that $n^P_{<L}$ is meant to be an upper bound on the number of possibilities below the height $L$, hence $n^P_{<L^-} \cdot n^S_{L^-} < n^P_{<L}$ must hold. (Note that as an immediate consequence, we also get $\prod_{K<L} n^S_K < n^P_{<L}$.) For the initial step, simply let $n^P_{<(0,0)} := 1$. (The interpretation of this number still makes sense, as there is exactly one trivial – empty – possibility "below the first height".) For any height $L$ not in $\mathbf{heights}_{\mathrm{ct}}$, let $n^P_{<L^+}$ be the minimal integer fulfilling the inequality

$$n^P_{<L} \cdot n^S_L < n^P_{<L^+}.$$

For a height $L = 4k+3 \in \mathbf{heights}_{\mathrm{ct}}$, let $n^P_{<L^+}$ be the minimal integer fulfilling the inequality

$$n^P_{<L} \cdot (n^S_L)^{k-1} < n^P_{<L^+}.$$

This will indeed be an upper bound on the number of possibilities below below height $L$, at least if we consider only modest conditions. This follows from modesty properties 4.3 (iii) and (iv). Note that there are at most $k-1$ many of those creatures in the segment starting at $L = 4k+3$.[12]

$\boldsymbol{n^R_{<L}}$: The definition of this sequence coding the rapidity of the reading depends on the forcing factor. (Technically, both $n^P_{<L}$ and $n^R_{<L}$ only require information about the previous height $L^-$, but it makes more sense to define them both at the beginning of the following height's definitions.) For technical reasons, we require $n^P_{<L} < n^R_{<L}$; apart from that, the definition's motivations should be clear once the concepts of rapid reading (Definition 6.1) respectively punctual reading (Definition 10.12) have been introduced.

As a general requirement, for any $L^- \in \mathbf{heights}$ we demand the following: Let $\ell := 4k$ (if $L^- = (4k, i)$) or $\ell := L^-$ (otherwise). We then require $n^R_{<L}$ to be at least

[12] In many cases (whenever the segments of the frame are not trivially short), we are actually way too generous here, but that does not matter.

large enough that it fulfils the inequality

$$n^P_{<L} < n^R_{<L} < \frac{2^{n^R_{<L}}}{\ell}.$$

(While this is not strictly necessary, it makes the proof of Lemma 7.8 slightly nicer.) In most cases, this is easily fulfilled already, anyways, but in the case that for any lower heights the subsequent definitions are smaller than would be required by the above, we just pick $n^R_{<L}$ larger instead.

Depending on the specific typegroup,

- For $L^- \in \mathbf{heights}_{\mathrm{nm}}$: Let $n^R_{<L} := n^P_{<L} + 2^{\max I_{L^-}+1}$.
- For $L^- \in \mathbf{heights}_{*\mathrm{n}}$: Let $n^R_{<L} := n^P_{<L} + 2 \nearrow (\max I_{L^-} + 2)$, where $2 \nearrow x := 2^{2^x}$.
- For $L^- \in \mathbf{heights}_{\mathrm{slalom}}$: In Definition 10.13, we will define a function $z \in \omega^\omega$ in which the value of $z(k)$ only depends on the value of $n^S_{4k+2}$;[13] we let $n^R_{<L} := n^P_{<L} + z(k)$.
- For $L^- \in \mathbf{heights}_{\mathrm{ct}}$: We do not have any additional demands for this part of the sequence.

$\boldsymbol{n^B_L}$: This is a straightforward definition; let

$$n^B_L := (n^B_{L^-})^{n^P_{<L} \cdot n^R_{<L}}$$

(with $n^B_{(0,0)^-} := 2$).[14]

Note that having defined these three numbers, all of the definitions of the various forcing factors can be made, though see the next paragraph regarding the slalom forcing posets.

$\boldsymbol{n^S_L}$: The definition of this sequence also depends on the specific forcing factor. We recall that this sequence is meant to be an upper bound on the size of the base sets at this height, i. e. the number of possibilities at this height.[15]

- For $L \in \mathbf{heights}_{\mathrm{nm}}$: The base sets for this factor is $2^{I_L}$ for some $I_L$, so let $n^S_L := 2^{|I_L|}$.
- For $L \in \mathbf{heights}_{*\mathrm{n}}$: For both cn and nn, the base set for these factors is the set of all subsets of $2^{I_L}$ of relative size $1 - 2^{-n^B_L}$; there are of course equally many of relative size $2^{-n^B_L}$, so let

$$n^S_L := \binom{2^{|I_L|}}{2^{|I_L| - n^B_L}}.$$

- For $L \in \mathbf{heights}_{\mathrm{slalom}}$: This is a bit different from the other cases. While for the other factors, the bound on the size is an a posteriori observation, for the slalom forcing factor, we actually define the bound $n^S_L$ on the size a priori and then (in Lemma 10.10) define the congenial sequence of function pairs $\langle f_\xi, g_\xi \mid \xi < \omega_1 \rangle$ such that they fit between $n^B_L$ and $n^S_L$. For $L = 4k+2$, we hence pick

$$n^S_L := (n^B_L)^{e_k^{10 \cdot 2^k}}$$

[13] The definition there is as follows: We let $\langle y_k \mid k \in \omega \rangle$ be a sequence of numbers such that $n^S_{4k+2} \leq 2^{y_k}$ and then let $z(k) := \sum_{\ell \leq k} y_\ell$.

[14] Note that this definition implies $(n^B_{L^-})^{n^P_{<L}} < n^B_L$ as well as $n^B_{L^-} \cdot 2^{n^S_{L^-}+1} < n^B_L$ and $2^{n^P_{<L} \cdot n^R_{<L}} \leq (n^B_{L^-})^{n^P_{<L} \cdot n^R_{<L}}$.

[15] Also note that the definitions of the intervals $I_L$ are such that $n^B_L < n^S_L$ holds for all $L$. However, in case the reader prefers not to verify this fact, they can just assume that the $I_L$ are chosen even larger such that this inequality holds.

(letting $e_k := n^B_{4k+2} + 1$, as defined in Lemma 10.10).[16]

- For $L \in \mathbf{heights}_{\mathrm{ct}}$: The base set for this factor is also $2^{I_L}$ for some $I_L$ (though with very different requirements on the size of $I_L$), so let $n^S_L := 2^{|I_L|}$.

We immediately see that $n^P_{<L}$ and $n^S_L$ work as intended:

**Lemma 4.9.** *For all $p \in \mathbb{Q}$ and $L \in$ **heights**, $|\operatorname{poss}(p, {<}L)| \leq n^P_{<L}$.*

*Proof.* For $L = (0,0)$, $\operatorname{poss}(p, {<}(0,0))$ is trivial and $n^P_{<(0,0)} = 1$. The rest follows from modesty by induction:

- for $L \in \mathbf{heights} \setminus \mathbf{heights}_{\mathrm{ct}}$, it follows from $n^P_{<L} \cdot n^S_L < n^P_{<L^+}$; and
- for $L = 4k+3 \in \mathbf{heights}_{\mathrm{ct}}$, it follows from $n^P_{<L} \cdot (n^S_L)^{k-1} < n^P_{<L^+}$.

The case distinction is necessary (here and in the definition of $n^P_{<L}$) because for $\mathrm{tg} \neq \mathrm{ct}$, modesty ensures that there is at most one non-trivial creature at any given height, while for ct, there can be up to $k-1$ many non-trivial objects of size at most $n^S_L$ (recall Definition 4.3 (iii)). □

(The function of $n^B_L$ and $n^R_{<L}$ will be shown in detail in section 5 and section 6, respectively.)

Having finally defined all parameters required for the forcing poset, we will now first remark on a few simple properties.

**Fact 4.10.** Since $\mathbb{Q}$ is a product, a condition $q$ is stronger than a condition $p$ if $q(\mathrm{tg})$ is stronger than $p(\mathrm{tg})$ for each $\mathrm{tg} \in$ **typegroups**; moreover, for each $\mathrm{t} \in \mathbf{types}_{\mathrm{modular}}$ (i. e. all but ct and nm), this statement can be broken down further to "$q(\alpha)$ is stronger than $p(\alpha)$ for each $\alpha \in A_{\mathrm{t}} \cap \operatorname{supp}(p)$".

To briefly summarise, "$q \leq p$" hence means that

- $\operatorname{trklgth}(q) \geq \operatorname{trklgth}(p)$ (the trunk may grow),
- $\operatorname{supp}(q) \supseteq \operatorname{supp}(p)$ (the support may grow),
- $\operatorname{supp}(q, \mathrm{nm}, 4k) \supseteq \operatorname{supp}(p, \mathrm{nm}, 4k)$ for each $k < \omega$ (above the trunk, the supports do not shrink for the lim inf factor),[17]
- the frame[18] of $q(\mathrm{ct})$ is coarser than the frame of $p(\mathrm{ct})$,
- for each $\alpha \in \operatorname{supp}(p)$ and each $L \in$ **heights** with $\mathrm{tg}(\alpha) = \mathrm{tg}(L)$, $q(\alpha, L) \subseteq p(\alpha, L)$ (strengthening the creatures on the old support), and

[16] Actually, any increasing sequence $e_k$ strictly greater than $n^B_L$ would work here and in Lemma 10.10.

[17] In the forcing construction of [FGKS17], this was true in a more general sense, but we have restricted the concept of the trunk to the lim inf factor and defined the support at a height slightly differently. These changes mean that in fact, the support at a certain height may shrink in the lim sup factors in a stronger condition, because the non-trivial creatures witnessing that a certain index $\alpha$ was already in the support by height $L$ may have been eliminated when extending the trunk, so $\alpha$ will then only enter the support at a later height. This conceptual change does not cause any problems, however.

[18] Recall Definition 3.11 for the definition of frames.

- for each $k < \omega$, $d(q)(4k) \geq d(p)(4k)$ (the halving parameters do not decrease).

**Lemma 4.11.** *For any given countable set of indices $B \subseteq A$, there is a condition $p$ such that* $\operatorname{supp}(p) = B$*. In particular, given any $\alpha \in A$, there is a condition $p$ such that* $\operatorname{supp}(p) = \{\alpha\}$*.*

*Proof.* We prove the simple case first: Given any $\alpha \in A$, define $p$ by letting $p(\alpha, L)$ be equal to the full base set for each $L \in \mathbf{heights}_{\mathrm{tg}(\alpha)}$. (If $\alpha \in A_{\mathrm{nm}}$, let the halving parameter sequence be equal to the constant 0 sequence. If $\alpha \in A_{\mathrm{ct}}$, let the frame be the trivial partition of $\mathbf{heights}_{\mathrm{ct}}$ into singleton segments.)

Given an arbitrary countable $B \subseteq A$ (without loss of generality such that $B$ has infinite intersection with $A_{\mathrm{t}}$ for each $\mathrm{t} \in \mathbf{types}$) instead, we first enumerate $B_{\mathrm{tg}} = A_{\mathrm{tg}} \cap B$ for each $\mathrm{tg} \in \mathbf{typegroups}_{\limsup}$ as $B_{\mathrm{tg}} =: \{\alpha_x, \alpha_{4+x}, \alpha_{2\cdot 4+x}, \alpha_{3\cdot 4+x}, \dots\}$ (with the $x$ depending on the tg, in such a way that the $4k + x$ correspond to the appropriate levels for this tg). Also enumerate $B_{\mathrm{nm}} = A_{\mathrm{nm}} \cap B$ as $B_{\mathrm{nm}} =: \{\beta_1, \beta_2, \dots\}$.

We then define the condition $p$ as follows:

- For $\mathrm{tg} \in \mathbf{typegroups}_{\limsup}$, let $x \in \{1, 2, 3\}$ be the appropriate value. We first let $p^*(\alpha_i, 4k + x)$ be equal to the full base sets for all $i \in \{x, 4 + x, 2 \cdot 4 + x, \dots\}$. Let the frame of $p^*(\mathrm{ct})$ be the trivial partition of $\mathbf{heights}_{\mathrm{ct}}$ into singleton segments.
  - Now use some appropriate diagonalisation of $B_{\mathrm{tg}}$ to thin out $p^*(\mathrm{tg})$ in such a way that in the resulting $p(\mathrm{tg})$ fulfils modesty[19] (which only requires reducing creatures to singletons or to minimal subcreatures with the desired norm, which in turn means such that the splitting nodes plus leaves are order-isomorphic to $2^{\leq k}$) while still fulfilling the requirements on the $\limsup$ of the norms.
  - (It follows from the definitions of the forcing factors that the $p^*(\mathrm{tg})$ fulfil the $\limsup$ conditions for each $\mathrm{tg} \in \mathbf{typegroups}_{\limsup}$, and so do the $p(\mathrm{tg})$ after diagonalisation.)
- For nm, we let $d(p)(4k) = 0$ for all $k < \omega$ and pick some increasing sequence $s_{4k}$ (with $s_0 = 1$) such that $\lim_{k\to\infty} \frac{s_{4k}}{4k+1} = 0$. We will let $S_{4k}(p) := \{\beta_1, \dots, \beta_{s_{4k}}\}$, so $\lim_{k\to\infty} \frac{|S_{4k}(p)|}{4k+1} = 0$ is fulfilled. Note that without loss of generality $|S_{4k}(p)| = s_{4k}$ will be much smaller than $k + 1$.
  - We will define $p$ such that $p(\mathrm{nm}, 4k)$ has at least norm $k$. For each $\alpha \in S_{4k}(p)$, pick a set $X_\alpha \subseteq J_{4k}$ of size $3^{(4k+1)\cdot 2^{k\cdot n^P_{<(4k,0)}}}$ (which means $\mu_{4k}(X) = 2^{k\cdot n^P_{<(4k,0)}}$) disjoint from $X_{\alpha'}$ for each $\alpha' \in S_{4k}$ with $\alpha \neq \alpha'$. We let $p(\alpha, (4k, j))$ be equal to the full base set for each $j \in X_\alpha$ and some arbitrary singletons elsewhere. The full base sets have cell norms larger than $n^B_{(4k,j)} > 2^{k\cdot n^P_{<(4k,0)}}$, so the whole compound creature $p(\mathrm{nm}, 4k)$ has norm $k$ and the $\liminf$ condition is fulfilled.
  - The choice of these $X_\alpha$ is possible because we only require

$$s_{4k} \cdot 3^{(4k+1)\cdot 2^{k\cdot n^P_{<(4k,0)}}} < (k+1) \cdot 3^{(4k+1)\cdot 2^{k\cdot n^P_{<(4k,0)}}}$$

[19] Recall Definition 4.3 for the definition of modesty.

many different sublevels to choose from to do that, and by our definition, $J_{4k} = 3^{(4k+1)\cdot 2^{4k\cdot n^P_{<(4k,0)}}}$ is larger than that. $\square$

Before proceeding, we recall Definition 3.13 (x) for the definition of $\mu_{4k}$ and the following combinatorial result from [FGKS17, Lemma 2.2.2].

**Lemma 4.12.** *Given $\ell \leq 4k$ and a family $\langle X_i \mid 1 \leq i \leq \ell \rangle$ of subsets of $J_{4k}$, there is a family $\langle X_i^* \mid 1 \leq i \leq \ell \rangle$ of pairwise disjoint sets such that for each $1 \leq i \leq \ell$, $X_i^* \subseteq X_i$ and $\mu_{4k}(X_i^*) \geq \mu_{4k}(X_i) - 1$.*

**Lemma 4.13.** *Given two conditions $p, q \in \mathbb{Q}$ with disjoint supports, identical (or compatible) frames and identical sequences of halving parameters, there is a condition $r$ stronger than both.*

*Proof.* Since $\frac{|S_{4k}(p)|}{4k+1}$ and $\frac{|S_{4k}(q)|}{4k+1}$ must both converge to 0, there is some $k_0$ such that $\frac{|S_{4k}(p)|}{4k+1} \leq \frac{1}{2}$ and $\frac{|S_{4k}(q)|}{4k+1} \leq \frac{1}{2}$ for all $k \geq k_0$. Define $p' \leq p$ and $q' \leq q$ as the conditions resulting from extending the trunk to $4k_0$ (and choosing arbitrary singletons within all non-trivial creatures below).

We first define the pseudo-condition $r^*$ as simply the union of $p'$ and $q'$ together with the finest frame coarser than the frames of $p(\mathrm{ct})$ and $q(\mathrm{ct})$. Of course, $r^*$ might not fulfil modesty. For each $\mathrm{tg} \in \mathbf{typegroups}_{\limsup}$, we use diagonalisation to thin out $r^*(\mathrm{tg})$ and pick appropriately small subcreatures in $r^*(\mathrm{ct})$ in such a way that the resulting $r(\mathrm{tg})$ fulfils modesty.

As for each $\alpha \in \mathrm{supp}(r^*) \cap A_{\mathrm{ct}}$, the minimal elements of the segments – which are the reference points for the norms – can only have shrunk, it follows that $r^*(\mathrm{ct})$ is indeed a valid condition.

For $r^*(\mathrm{nm})$, we need to do a bit more. Assume without loss of generality that for all $k \geq k_0$, $\|p'(\mathrm{nm}, 4k)\|_{\mathrm{nm},4k} \geq 2$ and $\|q'(\mathrm{nm}, 4k)\|_{\mathrm{nm},4k} \geq 2$. We do the following procedure for each $k \geq k_0$:

- Let $S_{4k}(p') =: \{\alpha_1, \ldots, \alpha_c\}$ and $S_{4k}(q') =: \{\beta_1, \ldots, \beta_d\}$ and note that $c, d \leq 2k$ by our choice of $k_0$.
- Let $n_p := \|p'(\mathrm{nm}, 4k)\|_{\mathrm{nm},4k}$, $n_q := \|q'(\mathrm{nm}, 4k)\|_{\mathrm{nm},4k}$ and $n := \min(n_p, n_q)$. For each $1 \leq i \leq c$ and $1 \leq j \leq d$, there must be sets $A_i \subseteq J_{4k}$ respectively $B_j \subseteq J_{4k}$ such that they witness the stacked creature norm of $p'(\alpha_i, 4k)$ respectively $q'(\beta_j, 4k)$ being at least $n$. We remark that since $n$ is at least 2, we know that $\mu_{4k}(A_i)$ and $\mu_{4k}(B_j)$ are at least $2^{2\cdot n^P_{<(4k,0)}} + d(4k)$, and hence $|A_i|$ and $|B_j|$ are at least $3^{(4k+1)\cdot(2^{2\cdot n^P_{<(4k,0)}} + d(4k))}$.
- We can apply Lemma 4.12 to the family $\langle A_1, \ldots, A_c, B_1, \ldots, B_d \rangle$ (since $c + d \leq 4k$) to get a family $\langle A_1^*, \ldots, A_c^*, B_1^*, \ldots, B_d^* \rangle$ of pairwise disjoint subsets of $J_{4k}$ such that for each $1 \leq i \leq c$ and each $1 \leq j \leq d$,
  - $\mu_{4k}(A_i^*) \geq 2^{n\cdot n^P_{<(4k,0)}} + d(4k) - 1$ and $\mu_{4k}(B_j^*) \geq 2^{n\cdot n^P_{<(4k,0)}} + d(4k) - 1$,
  - for each $a \in A_i^*$, $\|p'(\alpha_i, (4k, a))\|^{\mathrm{cell}}_{(4k,a)} \geq 2^{n\cdot n^P_{<(4k,0)}} + d(4k)$ (as before), and
  - for each $b \in B_j^*$, $\|q'(\beta_j, (4k, b))\|^{\mathrm{cell}}_{(4k,b)} \geq 2^{n\cdot n^P_{<(4k,0)}} + d(4k)$ (as before).

  Define $r(\mathrm{nm}, 4k)$ by keeping the creatures in these sublevel index sets and replacing the others by arbitrary singletons.

- It follows that $\|r(\mathrm{nm}, 4k)\|_{\mathrm{nm},4k} \geq n - 1$.

Hence the resulting $r$ is indeed a condition, and $r$ is stronger than both $p$ and $q$ by construction. □

**Corollary 4.14.** *Given a condition $p \in \mathbb{Q}$ and any $\alpha \in A \setminus \mathrm{supp}(p)$, there is a $q \leq p$ with $\mathrm{supp}(q) = \mathrm{supp}(p) \cup \{\alpha\}$.*

*Proof.* By Lemma 4.11, there is a condition $p_\alpha$ with support $\{\alpha\}$ (and, for $\alpha \in A_{\mathrm{ct}}$, a frame compatible with $p$) which we can amalgamate with $p$.

We only have to be careful in ony specific case – if $\alpha \in A_{\mathrm{nm}}$ and $\mathrm{supp}(p) \cap A_{\mathrm{nm}} \neq \varnothing$ (if $\mathrm{supp}(p) \cap A_{\mathrm{nm}} = \varnothing$, the amalgamation is straightforward). In this case, we replace $p_\alpha$ by the condition with identical creatures, but the halving parameters of $p$ instead; since $p$ is a condition, the halving parameters must be small enough such that $\liminf_{k\to\infty} \|p(\mathrm{nm}, 4k)\|_{\mathrm{nm},4k} = \infty$, and hence the same must hold for $p_\alpha$ with the same halving parameters.

Applying Lemma 4.13 to $p$ and $p_\alpha$, the resulting $q$ is as required. □

We can now define the generic sequences added by the forcing.

**Definition 4.15.** Let $G$ be a $\mathbb{Q}$-generic filter. For each $\mathrm{tg} \in$ **typegroups**, each type $\mathrm{t} \in \mathrm{tg}$ and each $\alpha \in A_{\mathrm{t}} \subseteq A_{\mathrm{tg}}$, let $\dot{y}_\alpha$ be the name for

$$\{(L, z) \mid L \in \mathbf{heights}_{\mathrm{tg}}, \exists p \in G\colon\ \mathrm{trklgth}(p) > L \wedge p(\alpha, L) = \{z\}\}.$$

For $\mathrm{tg} \neq \mathrm{nm}$, there is the equivalent, clearer representation (with $x$ the appropriate element of $\{1, 2, 3\}$)

$$\{(k, z) \mid 4k + x \in \mathbf{heights}_{\mathrm{tg}}, \exists p \in G\colon\ \mathrm{trklgth}(p) > 4k + x \wedge p(\alpha, 4k + x) = \{z\}\}.$$

We write $\dot{y}$ for $\langle \dot{y}_\alpha \mid \alpha \in A\rangle$.

We note a few simple facts about generics and possibilities.

**Fact 4.16.** Let $p \in \mathbb{Q}$ and $L \in$ **heights**.

- For $\eta \in \mathrm{poss}(p, {<}L)$, $p \wedge \eta \leq p$.[20]
- $p \wedge \eta$ and $p \wedge \eta'$ are incompatible if $\eta, \eta' \in \mathrm{poss}(p, {<}L)$ are distinct.
- $p \wedge \eta$ forces that $\dot{y}$ extends $\eta$, i. e. that $\dot{y}_\alpha$ extends $\eta(\alpha)$ for all $\alpha \in \mathrm{supp}(p)$. In particular, $p$ forces that $\dot{y}$ extends $p^{<\mathrm{trklgth}(p)}$.
- $\eta \in \mathrm{poss}(p, {<}L)$ iff $p$ does not force that $\eta$ is incompatible with $\dot{y}$.
- $\mathbb{Q}$ forces that $\dot{y}$ is defined everywhere. (This follows from Corollary 4.14.)

**Lemma 4.17.** *Given $q \leq p$ and $\eta \in \mathrm{poss}(q, {<}L)$, there is a unique $\vartheta \in \mathrm{poss}(p, {<}L)$ such that $q \wedge \eta \leq p \wedge \vartheta$.*

*Proof.* Recall Definition 4.6: Since the possibilities are structurally somewhat more complicated for ct, we need to take that into account when trimming $\eta$ to get $\vartheta$.

Set $L^* := \min\{M \in \mathbf{heights}_{\mathrm{ct}} \mid M \geq L, M$ is the minimum of a segment of $p(\mathrm{ct})\}$. Let $\vartheta{\restriction}_{A \setminus A_{\mathrm{ct}}} := \eta{\restriction}_{\mathrm{supp}(p) \setminus A_{\mathrm{ct}}}$ and $\vartheta{\restriction}_{A_{\mathrm{ct}}} := \eta{\restriction}_{\mathrm{supp}(p) \cap A_{\mathrm{ct}}}({<}L^*)$ (this restriction is necessary for technical reasons, because $p(\mathrm{ct})$ in general could have a finer frame than

[20] Recall Definition 4.7 for the definition of $p \wedge \eta$.

$q(\mathrm{ct})$). Uniqueness follows from the incompatibility of $p \wedge \vartheta$ and $p \wedge \vartheta'$ for distinct $\vartheta, \vartheta' \in \mathrm{poss}(p, {<}L)$. □

The first important fact about $\mathbb{Q}$ we will prove is the following:

**Lemma 4.18.** *Assuming* CH, $\mathbb{Q}$ *is* $\aleph_2$*-cc.*

*Proof.* Assume that $Z := \langle p_i \mid i < \omega_2 \rangle$ is a family of conditions. Using the $\Delta$-system lemma for families of countable sets and CH, we can find $\Delta \subseteq A$ and thin out $Z$ to a subset of the same size such that for any distinct $p, q \in Z$,

- $\Delta = \mathrm{supp}(p) \cap \mathrm{supp}(q)$,
- for all $k < \omega$, $d(p)(4k) = d(q)(4k)$,
- the frames of $p(\mathrm{ct})$ and $q(\mathrm{ct})$ are identical, and
- $p$ and $q$ are identical on $\Delta$, i. e. for all $\alpha \in \Delta$ and all $L \in \mathbf{heights}_{\mathrm{tg}(\alpha)}$, $p(\alpha, L) = q(\alpha, L)$.

By Lemma 4.13 (applied to $p$ and $q{\restriction}_{A \smallsetminus \Delta}$), there is some $r \in \mathbb{Q}$ stronger than both $p$ and $q$, hence $Z$ is not an antichain. □

**Lemma 4.19.** *Assume that* $B \subseteq A$ *and either* $A_{\mathrm{nm}} \subseteq B$ *or* $A_{\mathrm{nm}} \cap B = \varnothing$ *(i. e.* $A_{\mathrm{nm}}$ *is not split by* $B$*). Let* $\mathbb{Q}_B \subseteq \mathbb{Q}$ *consist of all* $p \in \mathbb{Q}$ *with* $\mathrm{supp}(p) \subseteq B$*. Then* $\mathbb{Q}_B$ *is a complete subforcing poset of* $\mathbb{Q}$*. (*$\mathbb{Q}_B$ *has the same structural properties as* $\mathbb{Q}$*, such as continuous reading,* $\omega^\omega$*-bounding etc.; of course, it does not necessarily force the same statements, e. g. it may force a smaller continuum.)*

*Proof.* It is clear that the "stronger" relation and incompatibility work as required for a complete embedding. We have to show that given $q \in \mathbb{Q}$, there is some $\pi(q) =: p \in \mathbb{Q}_B$ such that any $p' \in \mathbb{Q}_B$ with $p' \leq p$ is compatible with $q$ in $\mathbb{Q}$.

Let $\pi\colon \mathbb{Q} \to \mathbb{Q}_B$ be the projection mapping each $q \in \mathbb{Q}$ to $\pi(q) := q{\restriction}_{(\mathrm{supp}(q) \cap B)}$. Let $p := \pi(q)$ and fix an arbitrary $p' \in \mathbb{Q}_B$ stronger than $p$. Let $p^* := q{\restriction}_{A \smallsetminus B}$ and apply Lemma 4.13 to $p'$ and $p^*$ (keeping in mind that their frames are necessarily compatible, in case that is relevant) to get an $r \in \mathbb{Q}$ stronger than $p'$ and $q$. □

The proof is even more straightforward if we additionally assume that $A_{\mathrm{ct}} \subseteq B$ or $A_{\mathrm{ct}} \cap B = \varnothing$ (i. e. $A_{\mathrm{ct}}$ is not split by $B$, either); in that case, it is already clear from the product structure of $\mathbb{Q} = \mathbb{Q}_B \times \mathbb{Q}_{A \smallsetminus B}$ that $\mathbb{Q}_B$ is a complete subforcing poset of $\mathbb{Q}$.

## 5. Bigness

One key concept for many of the following proofs is the fact that by our construction, creatures at height $L$ are much, much bigger than creatures at height $L^-$ and much, much smaller than creatures at height $L^+$.[21] The exact nature of this size difference is encoded in the sequence $n^B_L$. While this concept is referred to as *completeness* in the older [GS93], we will be using the modern and more standard terminology of *bigness* from [FGKS17], while unifying the different concepts and generalising them even further in Definition 5.4.

**Definition 5.1.** Fix positive integers $c$ and $d$.

[21] Recall that $L^-$ and $L^+$ denote the predecessor and successor of a height $L$, respectively.

(i) We say a non-empty set $C$ and a norm

$$\|\cdot\|\colon (\mathcal{P}(C) \smallsetminus \{\varnothing\}) \to \mathbb{R}_{\geq 0}$$

on the subsets of $C$ are *$c$-big* (synonymously, *have $c$-bigness*) if the following holds: For each non-empty $X \subseteq C$ and each colouring $\chi\colon X \to c$ of $X$, there is a non-empty $Y \subseteq X$ such that $\chi\restriction_Y$ is constant and $\|Y\| \geq \|X\| - 1$. Equivalently, $(C, \|\cdot\|)$ is $c$-big if for each non-empty $X \subseteq C$ and each partition $X = X_1 \cup X_2 \cup \ldots \cup X_c$, there is some $i \in \{1, 2, \ldots, c\}$ such that $\|X_i\| \geq \|X\| - 1$.

(ii) We say $(C, \|\cdot\|)$ is *$(c, d)$-big* (synonymously, *has $(c, d)$-bigness*) if the following holds: For each non-empty $X \subseteq C$ and each colouring $\chi\colon X \to c$ of $X$, there is a non-empty $Y \subseteq X$ such that $|\operatorname{ran}\chi\restriction_Y| \leq d$ and $\|Y\| \geq \|X\| - 1$. Equivalently, $(C, \|\cdot\|)$ is $(c, d)$-big if for each non-empty $X \subseteq C$ and each partition $X = X_1 \cup X_2 \cup \ldots \cup X_c$, there is some $d$-tuple $\{i_1, i_2, \ldots, i_d\} \subseteq \{1, 2, \ldots, c\}$ such that $\|X_{i_1} \cup X_{i_2} \cup \ldots \cup X_{i_d}\| \geq \|X\| - 1$.[22]

(iii) We say $(C, \|\cdot\|)$ is *strongly $c$-big* (synonymously, *has strong $c$-bigness*) if in the above, even $\|Y\| \geq \|X\| - 1/c$ (respectively $\|X_i\| \geq \|X\| - 1/c$) holds.

Since the colouring and the partition formulations of properties (i) and (ii) above, respectively, are evidently equivalent, we will use whichever is more suited for that particular proof.

**Fact 5.2.** A few simple facts about bigness:

- If $(C, \|\cdot\|)$ has (strong) $c$-bigness, it also has (strong) $c'$-bigness for any $c' \leq c$.
- A simple example of a norm with $c$-bigness is $\log_c |\cdot|$.
- Modifying the norm to be $\log_c|\cdot|/c$ gives us strong $c$-bigness.
- An example of a $(c, d)$-big norm is $\log_{c/d} |\cdot|$.

The first fact can be generalised as follows:

**Lemma 5.3.** *If $c/d \leq b$ and $\|\cdot\|$ is $b$-big, then $\|\cdot\|$ is also $(c, d)$-big.*

*Proof.* Let $X = X_1 \cup X_2 \cup \ldots \cup X_c$; we have to find a $d$-tuple $\{i_1, i_2, \ldots, i_d\} \subseteq \{1, 2, \ldots, c\}$ such that $\|X_{i_1} \cup X_{i_2} \cup \ldots \cup X_{i_d}\| \geq \|X\| - 1$. Since $c/d \leq b$, we have $c \leq d \cdot b$. We regroup the partition of $X$ as

$$X = (X_1 \cup X_2 \cup \ldots \cup X_b) \cup (X_{b+1} \cup \ldots \cup X_{2b}) \cup \ldots \cup (X_{(d-1)b+1} \cup \ldots \cup X_{db}),$$

where $X_k := \varnothing$ for $c < k \leq d \cdot b$. Define $Y_j := \bigcup_{0 \leq i \leq d-1} X_{ib+j}$ for $1 \leq j \leq b$. Then

$$X = Y_1 \cup Y_2 \cup \ldots \cup Y_b$$

and, since $\|\cdot\|$ is $b$-big, there is some $j_0$ such that $\|Y_{j_0}\| \geq \|X\| - 1$. Then $\{j_0, b+j_0, 2b+j_0, 3b+j_0, \ldots, (d-1)b+j_0\}$ is the $d$-tuple we had to provide. Possibly, some of these indices are not even necessary – namely if they point to empty $X_k$; in that case, pick arbitrary replacement indices pointing towards actually existing sets. $\square$

[22] This partition formulation of $(c, d)$-bigness is precisely the definition of $(c, d)$-completeness from [GS93, Definition 2.2].

**Definition 5.4.** We extend the definition of $(c, d)$-bigness and strong $c$-bigness in the following way:

(i) We say $(C, \|\cdot\|)$ is *$e$-strongly $c$-big* (synonymously, *has $e$-strong $c$-bigness*) if $X_i \subseteq X$ is even such that $\|X_i\| \geq \|X\| - {}^1\!/\!_e$.
(ii) We say $(C, \|\cdot\|)$ is *$e$-strongly $(c, d)$-big* (synonymously, *has $e$-strong $(c, d)$-bigness*) if the $d$-tuple is even such that $\|X_{i_1} \cup X_{i_2} \cup \ldots \cup X_{i_d}\| \geq \|X\| - {}^1\!/\!_e$.

(We omit the equivalent colouring formulations of the same definitions.)

**Observation 5.5.** Dividing a $(c, d)$-big or $c$-big norm by $e$ yields an $e$-strongly $(c, d)$-big or $e$-strongly $b$-big norm, respectively.

Note that if ${}^c\!/\!_d \leq b$ and $\|\cdot\|$ has $e$-strong $b$-bigness, using the same method as in the preceding lemma gives us a $d$-tuple such that $\|X_{i_1} \cup X_{i_2} \cup \ldots \cup X_{i_d}\| \geq \|X\| - {}^1\!/\!_e$ and hence even $e$-strong $(c, d)$-bigness.

We have defined the norms of the various forcing factors in such a way that they have $n^B_L$-bigness at height $L$:

**Theorem 5.6.** *Recall the definitions of the norms in Definition 3.4 (for $\mathbb{Q}_{\text{nn}}$), Definition 3.5 (for $\mathbb{Q}_{\text{cn}}$), Definition 3.3 (for $\mathbb{Q}_{\text{slalom}}$), Definition 3.11 (for $\mathbb{Q}_{\text{ct},\kappa_{\text{ct}}}$) and Definition 3.13 (for $\mathbb{Q}_{\text{nm},\kappa_{\text{nm}}}$).*

*(i) For each $\mathrm{t} \in \mathbf{types}_{\text{modular}}$ and each $L \in \mathbf{heights}_{\text{tg(t)}}$, $(\mathbf{POSS}_{\mathrm{t},L}, \|\cdot\|_{\mathrm{t},L})$ has $n^B_L$-bigness. For $\mathrm{t} \in \{\text{nn}, \text{cn}\}$, we even have strong $n^B_L$-bigness. Letting $L := 4k + 2$, for $\xi \in \mathbf{types}_{\text{slalom}}$, we even have $n^B_{4k+2}$-strong $g_\xi(k)$-bigness (hence also $n^P_{<4k+2}$-strong $g_\xi(k)$-bigness) at height $4k + 2$.*
*(ii) Given a condition $p(\text{ct}) \in \mathbb{Q}_{\text{ct},\kappa_{\text{ct}}}$, for each segment $\bar{K} := \langle K_1, \ldots, K_m\rangle$ of its frame, $(p(\text{ct}, \bar{K}), \|\cdot\|_{\text{ct},K_1})$ has $n^B_{K_1}$-bigness.*
*(iii) For each $L \in \mathbf{heights}_{\text{nm}}$, $(\mathbf{POSS}_{\text{nm},L}, \|\cdot\|^{\text{cell}}_L)$ has strong $n^B_L$-bigness.*

*Proof.* For $\mathrm{t} = \xi \in \mathbf{types}_{\text{slalom}}$, the norm is the exemplary norm with $g_\xi((L - 2)/4)$-bigness from Fact 5.2 divided by $n^B_{4k+2}$, so by Observation 5.5, we have $n^B_{4k+2}$-strong $g_\xi(k)$-bigness. (Clearly, since $g_\xi((L-2)/4) \geq n^B_L$, this also implies $n^B_L$-bigness.)

For $\mathrm{t} = \text{nn}$, let $X \subseteq \mathbf{POSS}_{\text{nn},L}$ and fix a partition $X = X_1 \cup X_2 \cup \cdots \cup X_{n^B_L}$. Consider $\|X_i\|^{\text{intersect}}_L$ and let $r$ be the maximal such intersect norm (letting $i^*$ be some such index with $\|X_{i^*}\|^{\text{intersect}}_L = r$); hence $\|X_i\|^{\text{intersect}}_L \leq r$ for all $i \in n^B_L$, witnessed by sets $Y_i$. Then $Y := \bigcup Y_i$ witnesses that $\|X\|^{\text{intersect}}_L \leq n^B_L \cdot r$; hence

$$\|X\|_{\text{nn},L} \leq \frac{\log(n^B_L \cdot r)}{n^B_L \log n^B_L} \leq \frac{\log \max_{i<n^B_L} \|X_i\|^{\text{intersect}}_L}{n^B_L \log n^B_L} + \frac{1}{n^B_L} = \|X_{i^*}\|_{\text{nn},L} + \frac{1}{n^B_L}$$

and hence $i^*$ is an index such that $\|X_{i^*}\|_{\text{nn},L} \geq \|X\|_{\text{nn},L} - {}^1\!/\!_{n^B_L}$.

For $\mathrm{t} = \text{cn}$, we first remark that

$$\|M\|_{\text{cn},L} = \frac{\log|M|}{r \cdot n^B_L \log 3n^B_L} - s$$

for some positive $r, s$ only depending on $L$. Given $X \subseteq \mathbf{POSS}_{\mathrm{cn},L}$ and a colouring $c\colon X \to n^B_L$, there is some $c$-homogeneous $Y \subseteq X$ with $|Y| \geq |X|/n^B_L$ and hence

$$\|Y\|_{\mathrm{cn},L} \geq \frac{\log|X| - \log n^B_L}{r \cdot n^B_L \log 3n^B_L} - s \geq \|X\|_{\mathrm{cn},L} - \frac{\log n^B_L}{r \cdot n^B_L \log 3n^B_L} \geq \|X\|_{\mathrm{cn},L} - \frac{1}{n^B_L}.$$

For $\mathrm{t} = \mathrm{ct}$, the claim follows from Lemma 3.10, since (letting $K_1 =: 4k+3$) modesty (see Definition 4.3) ensures that $p(\mathrm{ct}, \bar{K})$ will contain at most a $(k-1)$-tuple of creature segments[23] which are non-trivial, all of which have a norm of at least $k$.[24]

For nm, the cell norm is exactly the exemplary norm with strong $n^B_L$-bigness from Fact 5.2. ☐

We remark that the $n^B_L$ thus precisely describe the (strong) bigness properties at height $L$.

**Corollary 5.7.** *Let $p \in \mathbb{Q}$, $\alpha \in \mathrm{supp}(p) \setminus A_{\mathrm{ct}}$ and let $L \in \mathbf{heights}_{\mathrm{tg}(\alpha)}$ be a relevant height. Then for each colouring $c\colon p(\alpha, L) \to n^B_L$, there is a $c$-homogeneous $q(\alpha, L) \subseteq p(\alpha, L)$ such that $\|q(\alpha, L)\|_{\mathrm{tg}(\alpha),L} \geq \|p(\alpha, L)\|_{\mathrm{tg}(\alpha),L} - 1$ for $\mathrm{tg}(\alpha) \neq \mathrm{nm}$ and $\|q(\alpha, L)\|^{\mathrm{cell}}_L \geq \|p(\alpha, L)\|^{\mathrm{cell}}_L - 1$ for $\mathrm{tg}(\alpha) = \mathrm{nm}$.*

*The same holds for $\mathrm{supp}(p) \cap A_{\mathrm{ct}}$: Let $\bar{L} = \langle L_1, \dots, L_m\rangle$ be a segment of the frame of $p(\mathrm{ct})$ such that $p(\mathrm{ct}, \bar{L})$ is non-trivial. Then for each colouring $c\colon p(\mathrm{ct}, \bar{L}) \to n^B_{L_1}$, there is a $c$-homogeneous $q(\mathrm{ct}, \bar{L}) \subseteq p(\mathrm{ct}, \bar{L})$ such that*

$$\|q(\alpha, \bar{L})\|_{\mathrm{ct},L_1} \geq \|p(\alpha, \bar{L})\|_{\mathrm{ct},L_1} - 1$$

*for all $\alpha \in \mathrm{supp}(p, \mathrm{ct}, \bar{L})$.*

Note that using the fact that $n^B_L$ is big with respect to $n^P_{<L}$ and $n^R_{<L}$, this can be iterated downwards. (We will not use the following consideration directly, but a similar one will come up later on.) First note that a colouring $c\colon \mathrm{poss}(p, {\leq}L) \to n^R_{<L}$ can be reinterpreted as a colouring $d\colon p(\alpha, L) \to (n^R_{<L})^{\mathrm{poss}(p,<L)}$. Since $(n^R_{<L})^{\mathrm{poss}(p,<L)} \leq (n^R_{<L})^{n^P_L} \leq n^B_L$, we can use the preceding corollary to make the colouring independent of the possibilities at height $L$. (For ct, keep in mind we have to treat tuples of creature segments[25] as units.) Iterating this downwards allows, for instance, the following:

- Given a colouring $c\colon \mathrm{poss}(p, {<}L) \to n^R_{<L'}$ for some relevant heights $L' < L$, we can strengthen $p(\alpha_K, K)$ to $q(\alpha_K, K)$ for all $L' \leq K < L$, decreasing the corresponding norms by at most 1, such that the colouring $c$ restricted to $\mathrm{poss}(q, {<}L)$ only depends on $\mathrm{poss}(q, {<}L')$. (The number of colours here limits how far we can iterate this downwards.)
- In particular, if $c\colon \mathrm{poss}(p, {<}L) \to 2$ for some relevant height $L$, we can find $q \leq p$ such that $\mathrm{poss}(q, {<}L)$ is $c$-homogeneous.

Finally, we will require one similar specific consequence of strong bigness:

[23] The phrase "tuple of creature segments" refers to the finitely many creature segments which compose $p(\mathrm{ct}, \bar{K})$.

[24] The demand imposed on their norms in the definition of modesty is only necessary to be able to apply the cited technical lemma without modifications.

[25] The phrase "tuple of creature segments" refers to the finitely many creature segments which compose $p(\mathrm{ct}, \bar{L})$.

**Lemma 5.8.** *Let $H$ be a finite subset of* $\mathbf{heights} \setminus \mathbf{heights}_{\mathrm{ct}}$ *and for each $L \in H$, assume we are given some type $\mathrm{t}_L \in \mathrm{tg}(L)$ and some $C_L \subseteq \mathbf{POSS}_{\mathrm{t}_L,L}$. Let $K$ be the minimum of $H$ and $F\colon \prod_{L\in H} C_L \to n^B_K$. Then there are $D_L \subseteq C_L$, with the norm of $D_L$ decreasing by at most $1/n^B_L$ when compared to $C_L$, such that the value of $F$ is constant on $\prod_{L\in H} D_L$.*

*Proof.* The case $|H| = 1$ is trivial, so assume $|H| \geq 2$ and let $M$ be the maximum of $H$. We construct $D_L$ by downwards induction on $L \in H$. Then $F$ can be written as a function from $C_M$ to $(n^B_K)^P$, where $P := \prod_{L\in H, L\neq M} C_L$. Since $(n^B_K)^{|P|} \leq (n^B_K)^{n^P_{<M}} \leq n^B_M$, we can use strong $n^B_M$-bigness to find $D_L$.

Continue the downwards induction with $H' := H \setminus \{M\}$. □

## 6. Continuous and Rapid Reading

We now prove the main properties which allow us to show that $\mathbb{Q}$ is proper and $\omega^\omega$-bounding.

**Definition 6.1.** Let $p \in \mathbb{Q}$ and let $\dot\tau$ be a $\mathbb{Q}$-name for an ordinal. We say that $p$ *decides $\dot\tau$ below the height $L$* if $p \wedge \eta$ decides $\dot\tau$ for each $\eta \in \mathrm{poss}(p, {<}L)$; in other words, there is a function $T\colon \mathrm{poss}(p, {<}L) \to \mathrm{Ord}$ with $p \wedge \eta \Vdash \dot\tau = T(\eta)$ for each $\eta \in \mathrm{poss}(p, {<}L)$.

We say that $p$ *essentially decides* $\dot\tau$ if there is some height $L$ such that $p$ decides $\dot\tau$ below $L$. Let $\dot r$ be a $\mathbb{Q}$-name for a countable sequence of ordinals. We say that $p$ *continuously reads* $\dot r$ if $p$ essentially decides each $\dot r(n)$.

Let $\dot s$ be a $\mathbb{Q}$-name for an element of $2^\omega$. We say that $p$ *rapidly reads* $\dot s$ if for each $L \in \mathbf{heights}$, $\dot s{\restriction}_{n^R_{<L}}$ is decided below $L$.

For $B \subseteq A$, we say that $p$ continuously reads $\dot r$ *only using indices in $B$* if $p$ continuously reads $\dot r$ and the value of $T(\eta)$ depends only on $\eta{\restriction}_B$. Analogously, we say that $p$ continuously reads $\dot r$ *not using indices in $B$* if $p$ continuously reads $\dot r$ only using indices in $A \setminus B$. (The same terminology will be used for "rapidly" instead of "continuously".)

**Observation 6.2.** The name "continuous reading" comes from the following consideration: For a fixed condition $p$, the possibilities form an infinite tree $T_p$; the set of branches $[T_p]$ carries a natural topology. A condition $p$ continuously reads some $\dot r$ iff there is a function $f\colon T_p \to \mathrm{Ord}^{<\omega}$ in the ground model such that for the natural (continuous) extension $F\colon [T_p] \to \mathrm{Ord}^\omega$ of $f$, $p \Vdash \dot r = F(\dot y)$, where $\dot y$ is the generic branch in $[T_p]$. In our case, the tree is finitely splitting and hence $T_p$ is compact, so continuity and uniform continuity coincide.

Rapid reading then is equivalent to a kind of Lipschitz continuity. We remark that the $n^R_{<L}$ describe "how rapidly" $p$ reads $\dot s$, i. e. they can be interpreted as corresponding to the Lipschitz constants.

**Lemma 6.3.** *If $p$ continuously (or rapidly) reads $\dot r$ and $q \leq p$, then $q$ continuously (or rapidly) reads $\dot r$. (The same holds if we add "only using indices in $B$" or "not using indices in $B$".)*

*Proof.* This follows immediately from Lemma 4.17. □

**Lemma 6.4.** *If $q \leq^* p$ and $p$ essentially decides $\dot{\tau}$, then $q$ also essentially decides $\dot{\tau}$.*[26]

*Proof.* (The proof may be skipped by the experienced reader.)

Since $q \leq^* p$, the frame of $q$ must be eventually coarser (and, without loss of generality, actually coarser) than the frame of $p$ (because if not, then we could strengthen the frame of $q$ in a way incompatible with the frame of $p$ and get $r \leq q$ incompatible with $p$). $p$ forces that $\dot{\tau}$ is decided below some height $L$. Without loss of generality, $q$ is coarser than $p$ everywhere above $L$; let $L^* \geq L$ be the minimum of the first segment of the frame of $q$(ct) which is entirely above $L$. Clearly, $p$ also forces that $\dot{\tau}$ is decided below $L^*$; so for each $\eta \in \mathrm{poss}(p, {<}L^*)$, we have $p \wedge \eta \Vdash \dot{\tau} = t$ for some $t \in \mathrm{Ord}$.

Since $q \leq^* p$ and since $L^*$ is the minimum of segments in the frames of both $p$ and $q$ (which ensures that the possibilities of $p$(ct) and $q$(ct) below $L^*$ have the same length), it is clear that $\mathrm{poss}(q{\restriction}_{\mathrm{supp}(p)}, {<}L^*) \subseteq \mathrm{poss}(p, {<}L^*)$ (because if not, then there would be an $r \leq q$ incompatible with $p$). Let $\vartheta \in \mathrm{poss}(q, {<}L^*)$. There is a unique $\eta \in \mathrm{poss}(q{\restriction}_{\mathrm{supp}(p)}, {<}L^*) \subseteq \mathrm{poss}(p, {<}L^*)$ such that $\vartheta = \eta{\restriction}_{\mathrm{supp}(p)}$. By $q \wedge \vartheta \leq^* p \wedge \eta$ it follows that $q \wedge \vartheta \Vdash \dot{\tau} = t$ must also hold. □

**Lemma 6.5.** *In the ground model, let $\kappa := \max(\aleph_0, |B|)^{\aleph_0}$ for some $B \subseteq A$. Then in the extension, there are at most $\kappa$ many reals which are read continuously only using indices in $B$; more formally, letting $G$ be a $\mathbb{Q}$-generic filter, there are at most $\kappa$ many reals $r$ such that there is a $p \in G$ and a name $\dot{s}$ such that $p$ continuously reads $\dot{s}$ only using indices in $B$ and such that $\dot{s}[G] = r$.*

*Proof.* (The proof may be skipped by the experienced reader.)

The argument is a variation of the usual "nice names" consideration. Given $p$ continuously reading some $\dot{s}$, we can define the canonical name $\dot{s}'$ continuously read by $p' := p{\restriction}_B$ such that $p$ forces $\dot{s} = \dot{s}'$. (We can do this by the following procedure: Let $L_n$ be the height such that $\dot{s}(n)$ is decided below $L_n$. For each $\eta \in \mathrm{poss}(p, {<}L_n)$, we have $p \wedge \eta \Vdash \dot{s}(n) = x_n^\eta$ for some $x_n^\eta$. Define $\dot{s}'(n)$ as the name containing all pairs $\langle \check{x}_n^\eta, p \wedge \eta \rangle$.)

Hence it suffices to prove that there are at most $\kappa$ many names of reals continuously read in this manner. There are at most $\kappa$ many countable subsets of $B$ and hence at most $\kappa$ many conditions $p'$ with $\mathrm{supp}(p') \subseteq B$, because

- there are countably many heights,
- for each such height $L \in \mathbf{heights}_{\mathrm{tg}}$, we have at most countably many indices in $B \cap A_{\mathrm{tg}}$, and
- for each such index $\alpha$, we have to choose one of finitely many creatures (very often: singletons) to be $p'(\alpha, L)$.

[26] A brief note on notation: "$q \leq^* p$" states that while $q$ may not *actually* be a stronger condition than $p$, any strengthening of $q$ will still be compatible with $p$ (and hence there will always be joint strengthenings of the two), which means that, in a sense, $q$ is "morally" a stronger condition than $p$. (Consequently, "$p =^* q$" just means "$q \leq^* p$ and $q \geq^* p$".)

Given any such $p'$, there are only $2^{\aleph_0}$ many possible ways to continuously read a real $\dot{s}'$ with respect to $p'$ (by picking the decision heights $L_n$ and the values $x_n^\eta$ for each of finitely many $\eta \in \mathrm{poss}(p', {<}L_n)$). $\square$

We will now first prove that given a condition continuously reading some $\dot{r} \in 2^\omega$, we can find a stronger condition rapidly reading $\dot{r}$, and only afterwards prove that we can densely find conditions continuously reading any $\dot{\tau} \in \mathrm{Ord}^\omega$. (This sequence of proofs, the same as in [FGKS17], makes for an easier presentation.)

**Theorem 6.6.** *Given $p$ continuously reading $\dot{r} \in 2^\omega$, there is a $q \leq p$ rapidly reading $\dot{r}$. (The same is true if we add "only using indices in $B$".)*

*Proof.* For each height $L$, we define:

($*_1$) $K_{\mathrm{dec}}(L)$ is the maximal height such that $\dot{r}\restriction_{n^R_{<K_{\mathrm{dec}}(L)}}$ is decided below $L$ by $p$.

The function $K_{\mathrm{dec}}$ is non-decreasing, and continuous reading already implies that $K_{\mathrm{dec}}$ is unbounded. (If it were bounded by $K$, that would mean that for any $K' \geq K$, $\dot{r}(n^B_{K'})$ were not essentially decided by $p$.) $K_{\mathrm{dec}}$ can, however, grow quite slowly. ($p$ rapidly reading $\dot{r}$ translates to $K_{\mathrm{dec}}(L) \geq L$ for all $L$.)

For all heights $K \leq L$ we define

$$\dot{x}^L_K := \dot{r}\restriction_{n^R_{<\min(K, K_{\mathrm{dec}}(L))}}$$

(which is, by definition, decided below $L$).

There are at most $2^{n^R_{<K}}$ many possible values for $\dot{x}^L_K$, since $n^R_{<\min(K,K_{\mathrm{dec}}(L))} \leq n^R_{<K}$.

In the following, we will only consider relevant heights. Recall that relevant heights are those that are either in $\mathbf{heights}_{\mathrm{ct}}$ and the minimum of a segment of the frame of $p(\mathrm{ct})$, or are in $\mathbf{heights}_{\mathrm{tg}}$ for some $\mathrm{tg} \neq \mathrm{ct}$ and are such that there is an $\alpha_L \in \mathrm{supp}(p) \cap A_{\mathrm{tg}}$ with a non-trivial $p(\alpha_L, L)$. For a relevant height $L \notin \mathbf{heights}_{\mathrm{ct}}$, we will use $\alpha_L$ to refer to the corresponding index.

<u>Step 1:</u> Fix a relevant $L$. We will choose, by downwards induction on all relevant $L' \leq L$, objects $C^L_{L'}$ (which will be either creatures $C^L_{L'} \subseteq p(\alpha_{L'}, L')$ or tuples of creature segments $C^L_{L'} \subseteq p(\mathrm{ct}, \mathrm{segm}(L'))$) and functions $\psi^L_{L'}$.[27]

<u>Step 1a:</u> To start the induction, for $L' = L$ we set $C^L_L := p(\alpha_L, L)$ respectively $C^L_L := p(\mathrm{ct}, \mathrm{segm}(L))$. We let $\psi^L_L$ be the function with domain $\mathrm{poss}(p, {<}L)$ assigning to each $\eta \in \mathrm{poss}(p, {<}L)$ the corresponding value of $\dot{x}^L_L$. (This means that $p \wedge \eta \Vdash \dot{x}^L_L = \psi^L_L(\eta)$ for each $\eta \in \mathrm{poss}(p, {<}L)$.)

<u>Step 1b:</u> We continue the induction on $L'$ and write $C' := C^L_{L'}$, $\psi' := \psi^L_{L'}$, $\dot{x}' := \dot{x}^L_{L'}$ for short.

Our plan is as follows:

- We will pick a creature $C'$ stronger than $p(\alpha_{L'}, L')$ respectively a tuple of creature segments $C'$ stronger than $p(\mathrm{ct}, L')$ such that the corresponding norm decreases by at most 1.

[27] Recall Definition 4.2 for the notation $p(\mathrm{ct}, \mathrm{segm}(L'))$. The phrase "tuple of creature segments" refers to the finitely many creature segments which compose $p(\mathrm{ct}, \mathrm{segm}(L'))$ (and subsets thereof).

- $\psi'$ will be a function with domain $\mathrm{poss}(p, {<}L')$ such that, loosely speaking,

  modulo $\langle C_K^L \mid L' \leq K < L\rangle$, each $\eta \in \mathrm{poss}(p, {<}L')$ decides $\dot{x}'$ to be $\psi'(\eta)$,

  or, more precisely, that $p \wedge \eta$ forces $\dot{x}' = \psi'(\eta)$ if the generic $\dot{y}$ is compatible with $C_K^L$ for all non-trivial heights $K$ with $L' \leq K < L$.[28]

We will define $C', \psi'$ as follows: Let $L''$ be the smallest relevant height above $L'$. By induction, we already have that $\psi'' := \psi_{L''}^L$ is a function with domain $\mathrm{poss}(p, {<}L'')$ such that modulo $\langle C_K^L \mid L'' \leq K < L\rangle$, each $\eta \in \mathrm{poss}(p, {<}L'')$ decides $\dot{x}'' := \dot{x}_{L''}^L$ to be $\psi''(\eta)$.

Let $\psi_*''(\eta)$ be the restriction of $\psi''(\eta)$ to $n^R_{<\min(L', K^{\mathrm{dec}}(L))}$. This means that $\psi_*''$ maps each $\eta \in \mathrm{poss}(p, {<}L'')$ to a restriction of $\dot{x}''$ – a potential value for $\dot{x}'$.

We can refactor $\psi_*''$ as a function $\psi_*'' : X \times Y \to Z$, where $X := \mathrm{poss}(p, {<}L')$, $Y := p(\alpha_{L'}, L')$ respectively $Y := p(\mathrm{ct}, \mathrm{segm}(L'))$ and $Z$ is the set of possible values of $\dot{x}'$, which has at most size $2^{n^R_{<L'}}$. This implicitly defines a function from $Y$ to $Z^X$; with $|Z^X| \leq 2^{n^P_{<L'} \cdot n^R_{<L'}}$, we can by Corollary 5.7 use bigness at height $L'$ to find $C' \subseteq p(\alpha_{L'}, L')$ respectively $C' \subseteq p(\mathrm{ct}, \mathrm{segm}(L'))$ (with the norm decreasing by at most 1) such that $\psi_*''$ does not depend on the height $L'$. From this, we get a natural definition of $\psi'$.

<u>Step 2:</u> We perform a downwards induction as in step 1 (always in the original $p$) from each relevant height $L$, thus defining for each relevant $K < L$ the creatures/tuples of creature segments $C_K^L$ and a function $\psi_K^L$ fulfilling

($*_2$) modulo $\langle C_{K'}^L \mid K \leq K' < L\rangle$, each $\eta \in \mathrm{poss}(p, {<}K)$ decides $\dot{x}_K^L$ to be $\psi_K^L(\eta)$.

The corresponding norms of these creatures/tuples of creature segments decrease by at most 1.

<u>Step 3:</u> For a given $K$, there are only finitely many possibilities for both $C_K^L$ and $\psi_K^L$. So by König's Lemma there necessarily exists a sequence $\langle C_K^*, \psi_K^* \mid K \text{ relevant}\rangle$ such that

($*_3$) for each $L$, there is $L^* > L$ such that for all $K \leq L$, $\langle C_K^{L^*}, \psi_K^{L^*}\rangle = \langle C_K^*, \psi_K^*\rangle$.

(These $\langle C_K^*, \psi_K^* \mid K \text{ relevant}\rangle$ thus form an infinite branch in the tree of all $\langle C_K^L, \psi_K^L\rangle$.)

<u>Step 4:</u> To define $q$, we replace all creatures and tuples of creature segments of $p$ by $C_K^* \subseteq p(\alpha_K, K)$ respectively $C_K^* \subseteq p(\mathrm{ct}, \mathrm{segm}(K))$. Thus $q$ has the same support as $p$, the same trunk, the same frame and the same halving parameters, and all corresponding norms decrease by at most 1, hence $q$ actually is a condition. We now claim that $q$ rapidly reads $\dot{r}$, i. e. we claim that each $\eta \in \mathrm{poss}(q, {<}K)$ decides $\dot{r}{\restriction}_{n^R_{<K}}$.

<u>Step 5:</u> To show this, we fix $K$ and pick a $K' > K$ such that $K_{\mathrm{dec}}(K') \geq K$. According to its definition Eq. ($*_1$), this means that $\dot{r}{\restriction}_{n^R_{<K}}$ is decided below $K'$. Now pick $L^* > K'$ per Eq. ($*_3$) and note that per Eq. ($*_2$), $\dot{x}_K^{L^*}$ is decided below $K$ by each $\eta \in \mathrm{poss}(p, {<}K)$ to be $\psi_K^{L^*}(\eta)$, modulo $\langle C_{K''}^{L^*} \mid K \leq K'' < L^*\rangle$. Since $K_{\mathrm{dec}}(K') \geq K$ and $L^* \geq K'$ (from which $K_{\mathrm{dec}}(L^*) \geq K_{\mathrm{dec}}(K')$ follows), we have

[28] We could also introduce a term referring to "$p$, but replacing all $p(\alpha_K, K)$ by $C_K^L$" here, but for notational simplicity, we eschew this.

$\min(K^{\mathrm{dec}}(L^*), K) = K$ and hence $\dot{x}_K^{L^*} = \dot{r}\restriction_{n^R_{<K}}$. As we had $K_{\mathrm{dec}}(K') \geq K$, $\dot{x}_K^{L^*}$ is already decided below $K'$ by the original condition $p$. Hence, in "modulo $\langle C_{K''}^{L^*} \mid K \leq K'' < L^* \rangle$", we can actually disregard any $K'' > K'$.

However, by Eq. ($*_3$) we know that $q$ has as its creatures and tuples of creature segments $C_L^{L^*} = C_L^*$ for all relevant $L < K'$. Hence $q$ forces that the generic $\dot{y}$ be compatible with $C_L^{L^*}$ for all non-trivial $K \leq L < K'$. From that, we immediately have that $\psi_K^{L^*} = \psi_K^*$ correctly computes $\dot{x}_K^{L^*} = \dot{r}\restriction_{n^R_{<K}}$ modulo $q$, and hence $q$ decides $\dot{r}\restriction_{n^R_{<K}}$ below $K$. As Step 5 holds for any $K$, $q$ rapidly reads $\dot{r}$. $\square$

## 7. Unhalving and the Proof of Continuous Reading

This section will contain proofs constructing a fusion sequence of conditions in $\mathbb{Q}$. While the lemmata and theorems could be formulated more generally, this would not give any additional insight, as they are only of a technical character. Since the structure of the possibilities in the ct factor is a bit unpleasant to work with, we will anchor these fusion constructions at the easiest possible nm levels of a condition, which are those which lie exactly between the maximal height of one segment in the frame of the ct factor of the condition and the minimal height of the frame segments immediately succeeding it.

**Definition 7.1.** Given a condition $p \in Q$, we call a lim inf level $4k$ (respectively $4k+1 \in \mathbf{heights}_{*\mathrm{n}}$ respectively $4k+2 \in \mathbf{heights}_{\mathrm{slalom}}$) *$p$-agreeable* if the heights $4k-1$ and $4k+3$ in $\mathbf{heights}_{\mathrm{ct}}$ are such that $4k-1 = \max(\mathrm{segm}(4k-1))$ and $4k+3 = \min(\mathrm{segm}(4k+3))$.

Restricting our constructions to use these heights as the stepping stones makes the possibilities easier to think about.

This section will also be the only time we actually use the halving parameters, in the form of the following operation on conditions:

**Definition 7.2.** Given a condition $q \in \mathbb{Q}$ and $4h < \omega$, define $r := \mathrm{half}(q, \geq 4h)$ as the condition obtained by replacing the halving parameters $d(q)(4k)$ of $q$ by

$$\begin{aligned} d(r)(4k) &:= d^*(q)(4k) \\ &:= d(q)(4k) + \frac{\min\{\|q(\alpha, 4k)\|_{4k}^{\mathrm{stack}} \mid \alpha \in \mathrm{supp}(q, \mathrm{nm}, 4k)\} - d(q)(4k)}{2} \end{aligned}$$

for all $4k \geq 4h$.

It is clear that for $r := \mathrm{half}(q, \geq 4h)$, the compound creature $r(\mathrm{nm}, 4k)$ is identical to $q(\mathrm{nm}, 4k)$ for each $4k < 4h$ and that for $4k \geq 4h$, the norm of the compound creature $r(\mathrm{nm}, 4k)$ has decreased by exactly $1/n^P_{<(4k,0)}$ compared to the norm of $q(\mathrm{nm}, 4k)$ (respectively, has remained 0 in case $4h \leq 4k < \mathrm{trklgth}(q)$).

The point of this is the following: Given $q \in \mathbb{Q}$ with relatively large nm norms and $r \leq \mathrm{half}(q, \geq 4h)$ such that some nm norms of $r$ are rather small, we can find an "unhalved" version $s$ of $r$ such that $s \leq q$, $s$ has relatively large nm norms and $s =^* r$. We will use this unhalving operation in the first part of the proof of continuous reading.

**Lemma 7.3.** *Fix $M \in \mathbb{R}$ and $h < \omega$. Given $q \in \mathbb{Q}$ such that $\|q(\mathrm{nm}, 4k)\|_{\mathrm{nm},4k} \geq M$ for all $4k \geq 4h$ as well as $r \leq \mathrm{half}(q, {\geq}4h)$ such that $\mathrm{trklgth}(r) = 4h$ and $\|r(\mathrm{nm}, 4k)\|_{\mathrm{nm},4k} > 0$ for all $4k \geq 4h$, there are $s \in \mathbb{Q}$ and $h^* > h$ such that*

*(i)* $s \leq q$,
*(ii)* $\mathrm{trklgth}(s) = 4h$,
*(iii)* $\|s(\mathrm{nm}, 4k)\|_{\mathrm{nm},4k} \geq M$ *for all* $4k \geq 4h^*$,
*(iv)* $s$ *is identical to* $r$ *above* $(4h^*, 0)$, *which means:* $s(\alpha, L) = r(\alpha, L)$ *for each sensible choice of* $\alpha \in \mathrm{supp}(r) = \mathrm{supp}(s)$ *and* $L \in$ **heights** *(and their halving parameters and frames are identical above* $(4h^*, 0)$*),*
*(v)* $\|s(\mathrm{nm}, 4k)\|_{\mathrm{nm},4k} \geq M - {}^1/_{n^P_{<(4k,0)}} \geq M - {}^1/_{n^P_{<(4h,0)}}$ *for all* $4h \leq 4k < 4h^*$, *and*
*(vi)* $\mathrm{poss}(s, {<}(4h^*, 0)) = \mathrm{poss}(r, {<}(4h^*, 0))$.

*Taken together, (iv) and (vi) imply $s =^* r$ (see Figure 5 and also recall the footnote in Lemma 6.4) and hence by Lemma 6.4, if $r$ essentially decides some $\dot\tau$, then so does $s$.*

*Proof.* Let $h^\dagger \geq h$ such that $\|r(\mathrm{nm}, 4k)\|_{\mathrm{nm},4k} > M$ for all $4k \geq 4h^\dagger$. Set $h^* := h^\dagger + 1$. Define $s$ to be identical to $r$ except for the fact that for all $4h \leq 4k < 4h^*$, we replace the halving parameters $d(r)(4k)$ by $d(q)(4k)$. (This means that for $4h \leq 4k < 4h^*$ we have $d(s)(4k) = d(q)(4k)$.)

It is clear that (i)–(iv) and (vi) are true; it remains to show that (v) holds. Fix $k$ such that $4h \leq 4k < 4h^*$; we have to show the inequality

$$\begin{aligned} \|s(\mathrm{nm}, 4k)\|_{\mathrm{nm},4k} &= \frac{\log_2\left(\min\{\|s(\alpha, 4k)\|^{\mathrm{stack}}_{4k} \mid \alpha \in \mathrm{supp}(s, \mathrm{nm}, 4k)\} - d(s)(4k)\right)}{n^P_{<(4k,0)}} \\ &\geq M - \frac{1}{n^P_{<(4k,0)}}. \end{aligned}$$

Recall the definition of $d^*$ in the preceding definition; as $d^*(q)(4k)$ were the halving parameters of $\mathrm{half}(q, {\geq}4h)$ and $r \leq \mathrm{half}(q, {\geq}4h)$, we know that $d(r)(4k) \geq d^*(q)(4k)$.

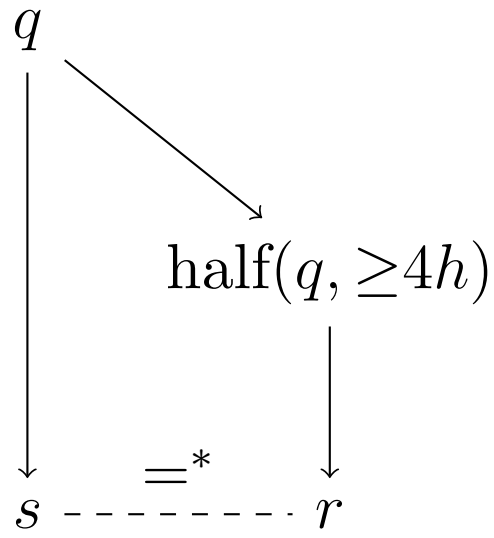

FIGURE 5. The construction in Lemma 7.3.

Since we assumed $\|r(\mathrm{nm}, 4k)\|_{\mathrm{nm},4k} > 0$, we know that

$$0 < \frac{\log_2\big(\min\{\|r(\alpha,4k)\|^{\mathrm{stack}}_{4k} \mid \alpha \in \mathrm{supp}(r,\mathrm{nm},4k)\} - d(r)(4k)\big)}{n^P_{<(4k,0)}}$$
$$= \frac{\log_2\big(\min\{\|s(\alpha,4k)\|^{\mathrm{stack}}_{4k} \mid \alpha \in \mathrm{supp}(s,\mathrm{nm},4k)\} - d(r)(4k)\big)}{n^P_{<(4k,0)}}$$

Fixing any $\beta \in \mathrm{supp}(s,\mathrm{nm},4k) = \mathrm{supp}(r,\mathrm{nm},4k)$, this shows

$$0 < \log_2\big(\|s(\beta,4k)\|^{\mathrm{stack}}_{4k} - d(r)(4k)\big)$$

and thus

$$\|s(\beta,4k)\|^{\mathrm{stack}}_{4k} > d(r)(4k) \geq d^*(q)(4k)$$
$$= d(q)(4k) + \frac{\min\{\|q(\alpha,4k)\|^{\mathrm{stack}}_{4k} \mid \alpha \in \mathrm{supp}(q,\mathrm{nm},4k)\} - d(q)(4k)}{2}$$

Hence (recalling $d(q)(4k) = d(s)(4k)$)

$$\|s(\beta,4k)\|^{\mathrm{stack}}_{4k} - d(s)(4k) \geq \frac{\min\{\|q(\alpha,4k)\|^{\mathrm{stack}}_{4k} \mid \alpha \in \mathrm{supp}(q,\mathrm{nm},4k)\} - d(q)(4k)}{2}$$

for any $\beta \in \mathrm{supp}(s,\mathrm{nm},4k) = \mathrm{supp}(r,\mathrm{nm},4k)$. Taking $\log_2$ and then dividing by $n^P_{<(4k,0)}$ yields

$$\frac{\log_2\big(\|s(\beta,4k)\|^{\mathrm{stack}}_{4k} - d(s)(4k)\big)}{n^P_{<(4k,0)}} \geq \|q(\mathrm{nm},4k)\|_{\mathrm{nm},4k} - \frac{1}{n^P_{<(4k,0)}}$$

and consequently (since this holds for any $\beta$)

$$\|s(\mathrm{nm},4k)\|_{\mathrm{nm},4k} \geq \|q(\mathrm{nm},4k)\|_{\mathrm{nm},4k} - \frac{1}{n^P_{<(4k,0)}} \geq M - \frac{1}{n^P_{<(4k,0)}},$$

proving (v). □

To prove that we can densely find conditions continuously reading a given name, we will first prove the following auxiliary lemma.

**Lemma 7.4.** *Let $\dot\tau$ be an arbitrary $\mathbb{Q}$-name and let $p^* \in \mathbb{Q}$ and $\ell^* < \omega$ and $M^* \geq 1$ be such that $4\ell^*$ is $p^*$-agreeable and $\|p^*(\mathrm{nm},4k)\|_{\mathrm{nm},4k} \geq M^* + 1$ holds for all $4k \geq 4\ell^*$. Then there is a condition $q$ such that:*

- *$q \leq p^*$,*
- *$q$ essentially decides $\dot\tau$,*
- *below $(4\ell^*, 0)$, $q$ and $p^*$ are identical on $\mathrm{supp}(p^*)$, and any $\alpha \in \mathrm{supp}(q) \setminus \mathrm{supp}(p^*)$ only enter the support (i. e. have their first non-trivial creature) above $(4\ell^*, 0)$ (as a consequence, $4\ell^*$ also is $q$-agreeable), and*
- *$\|q(\mathrm{nm},4k)\|_{\mathrm{nm},4k} \geq M^*$ for all $4k \geq 4\ell^*$.*

*Proof.* The proof consists of three parts.

**Part 1: finding intermediate deciding conditions by applying the unhalving lemma (Lemma 7.3)**

Suppose we are given $p \in \mathbb{Q}$, $\ell < \omega$ and $M \geq 1$ such that $4\ell$ is $p$-agreeable and $\|p(\mathrm{nm},4k)\|_{\mathrm{nm},4k} \geq M + 1$ for all $4k \geq 4\ell$. We construct an extension $r(p, 4\ell, M)$ of $p$ with certain properties:

First, enumerate $\mathrm{poss}(p, {<}(4\ell, 0))$ as $(\eta^1, \dots, \eta^m)$ and note that $m \leq n^P_{<(4\ell,0)}$. Setting $p^0 := q^0 := p$, we now inductively construct conditions $p^1 \geq \dots \geq p^m$ and auxiliary conditions $\tilde{q}^1, q^1, \dots, \tilde{q}^m, q^m$ such that for each $n < m$, the following properties hold:

(1) $\tilde{q}^{n+1}$ is derived from $p^n$ by replacing everything below $(4\ell, 0)$ (in $\mathrm{supp}(p)$) with $\eta^{n+1}$.
  - By (3) below, we will have $\mathrm{trklgth}(\tilde{q}^{n+1}) = 4\ell$.
  - For $n = 0$, $\tilde{q}^1$ is just $p^0 \wedge \eta^1$; but for $n \geq 1$, $\eta^{n+1}$ will not actually be in $\mathrm{poss}(p^n, {<}4\ell)$, so we cannot formally use that notation.
  - Note that in general, $\mathrm{supp}(p^n)$ will be larger than $\mathrm{supp}(p)$, so we do not replace *everything* below $4\ell$ with $\eta^{n+1}$, but only the part that is in $\mathrm{supp}(p)$.

(2) $q^{n+1} \leq \tilde{q}^{n+1}$. (Note that, obviously, $q^{n+1} \not\leq q^n$, since their trunks are different and the conditions are hence incompatible.)

(3) $\mathrm{trklgth}(q^{n+1}) = 4\ell$. (This means that by strengthening $\tilde{q}^{n+1}$ to $q^{n+1}$, we do not increase the trunk lengths.)

(4) $\|q^{n+1}(\mathrm{nm}, 4k)\|_{\mathrm{nm},4k} \geq M + 1 - {}^{n+1}/_{n^P_{<(4\ell,0)}}$ for all $4k \geq 4\ell$.

(5) One of the following two cases holds:
  - It is possible to choose $q^{n+1}$ such that it essentially decides $\dot{\tau}$.
  - Otherwise, $q^{n+1} := \mathrm{half}(\tilde{q}^{n+1}, {\geq}4\ell)$

  More explicitly: If the "decision" case is possible under the side conditions (2)–(4), then we use it (i. e. strengthen the condition to decide). If not, only then do we halve – and thereby certainly satisfy (2)–(4).

(6) We define $p^{n+1}$ as follows: Below $(4\ell, 0)$, $p^{n+1}$ is identical to $p$ on $\mathrm{supp}(p)$; above (including) $(4\ell, 0)$ as well as outside $\mathrm{supp}(p)$, $p^{n+1}$ is identical to $q^{n+1}$. In detail:
  - For all $\alpha \in \mathrm{supp}(p)$, $p^{n+1}(\alpha, L) := p(\alpha, L)$ for all sensible $L < (4\ell, 0)$.
  - For all $\alpha \in \mathrm{supp}(p)$, $p^{n+1}(\alpha, L) := q^{n+1}(\alpha, L)$ for all $L \geq (4\ell, 0)$.
  - For all $\beta \in \mathrm{supp}(q) \setminus \mathrm{supp}(p)$, $p^{n+1}(\beta) := q^{n+1}(\beta)$.

  (Note that as we required $\mathrm{trklgth}(q^{n+1})$ to remain $4\ell$, any newly added indices $\beta$ can only start having non-trivial creatures starting with height $(4\ell, 0)$ by modesty.)

(7) $p^{n+1} \leq p^n$, so the $\langle p^n \mid n \leq m \rangle$ are a descending sequence of conditions.

Ultimately, we define $r(p, 4\ell, M) := p^m$ (the last of the $p^n$ constructed above). $r := r(p, 4\ell, M)$ fulfils $r \leq p$ and $\|r(\mathrm{nm}, 4k)\|_{\mathrm{nm},4k} \geq M$ for all $4k \geq 4\ell$. As $r$ differs from $p$ only above $(4\ell, 0)$, it is also clear that $4\ell$ is $r$-agreeable.

Furthermore, $r$ has the following important decision property:

$$(*_4)\qquad \begin{aligned} &\text{If } \eta \in \mathrm{poss}(r, {<}(4\ell, 0)) \text{ and if there is an } s \leq r \wedge \eta \text{ such that } s \text{ essentially} \\ &\quad \text{decides } \dot{\tau},\ \mathrm{trklgth}(s) = 4\ell \text{ and } \|s(\mathrm{nm}, 4k)\|_{\mathrm{nm},4k} > 0 \text{ for all } 4k \geq 4\ell, \\ &\quad \text{then } r \wedge \eta \text{ already essentially decides } \dot{\tau}. \end{aligned}$$

To prove Eq. ($*_4$), note the following: $\eta$ canonically corresponds to $\eta \cap \mathrm{supp}(p) =: \eta^{n+1} \in \mathrm{poss}(p, {<}(4\ell, 0))$, therefore $s \leq r \wedge \eta \leq q^{n+1} \leq \tilde{q}^{n+1}$. We thus only have to show that $q^{n+1}$ was constructed using the "decision" case. Assume, towards an indirect proof, that this was not the case; so $q^{n+1}$ came about by halving $\tilde{q}^{n+1}$. Since $s$ is stronger than $\mathrm{half}(\tilde{q}^{n+1}, {\geq}4\ell)$, we can use Lemma 7.3 and unhalve $s$ to obtain some $s' \leq \tilde{q}^{n+1}$ with large norm such that $s' =^* s$. This means we could have used the "decision" case after all, which finishes this step of the proof.

**Part 2: iterating the intermediate conditions to define $q$**

Given $p^*$, $\ell^*$ and $M^*$ as in the lemma's statement, we inductively construct conditions $p_n$ and accompanying $\ell_n < \omega$ for each $n \geq 0$. Let $p_0 := p^*$ and $\ell_0 := \ell^*$. Given $p_n$ and $\ell_n$ such that $4\ell_n$ is $p_n$-agreeable, define $p_{n+1}$ and $\ell_{n+1}$ as follows:

- Choose $\ell_{n+1} > \ell_n$ such that:
  - $4\ell_{n+1}$ is $p_n$-agreeable,
  - $\|p_n(\mathrm{nm}, 4k)\|_{\mathrm{nm},4k} \geq M^* + n + 1$ for all $4k \geq 4\ell_{n+1}$, and
  - for each $\alpha \in \mathrm{supp}(p_n, (4\ell_n, 0)) \setminus A_{\mathrm{nm}}$ of type t, there is a height $L$ with $(4\ell_n, J_{4\ell_n} - 1) < L < (4\ell_{n+1}, 0)$ such that $\|p_n(\alpha, L)\|_{\mathrm{t},L} \geq M^* + n + 1$.
- Set $p_{n+1} := r(p_n, 4\ell_{n+1}, M^* + n + 1)$. (By the construction of $r$ in the previous part, it follows that $4\ell_{n+1}$ then also is $p_{n+1}$-agreeable.)

Thus $\langle p_n \mid n < \omega \rangle$ is a descending sequence of conditions, which converges to a condition $q \in \mathbb{Q}$. To verify that $q$ is indeed a condition, note the following: By construction, we have $\|q(\mathrm{nm}, 4k)\|_{\mathrm{nm},4k} \geq M^* + n$ for all $4k \geq 4\ell_{n+1}$. For all other types t and all indices $\alpha \in \mathrm{supp}(q) \cap A_{\mathrm{t}}$, we have assured the existence of a subsequence of creatures of strictly increasing norms of $q(\alpha)$, since below any $(4\ell_{n+1}, 0)$, $q$ is equal to $p_{n+1}$ (and also to $p_n$). Thus, $q$ is indeed a condition. Clearly, $q \leq p^*$ also holds.

In the next and final part, we will show that $q$ essentially decides $\dot{\tau}$ (proving the lemma). The following property will be central to the proof:

$$(*_5)\quad \begin{array}{l}\text{If } \eta \in \mathrm{poss}(q, {<}4(\ell_m, 0)) \text{ for some } m \text{ and if there is an } r \leq q \wedge \eta \text{ such that } r \\ \quad\text{essentially decides } \dot{\tau},\ \mathrm{trklgth}(r) = 4\ell_m \text{ and } \|r(\mathrm{nm}, 4k)\|_{\mathrm{nm},4k} > 0 \text{ for all} \\ \quad 4k \geq 4\ell_m, \text{ then } q \wedge \eta \text{ already essentially decides } \dot{\tau}.\end{array}$$

To prove Eq. ($*_5$), note that $\eta$ canonically corresponds to some $\eta^{n+1}$ which was already considered as a possible trunk when constructing the intermediate condition $p_m := r(p_{m-1}, 4\ell_m, M^* + m)$, so we can use Eq. ($*_4$) to conclude Eq. ($*_5$).

**Part 3: using bigness to thin out $q$ and prove its essential decision property**

The final part of the proof is essentially a rerun of the proof of Theorem 6.6. This is the main reason we proved rapid reading before continuous reading, as the idea of the proof is easier to digest in the rather simpler Theorem 6.6, in our opinion. The difference is that this time, we do not homogenise with respect to the potential values for some names, but instead with respect to whether $q \wedge \eta$ essentially decides $\dot{\tau}$ or not.

Step 1: Fix a relevant height $L > (4\ell_0, 0)$. We will choose, by downwards induction on all relevant $L'$ with $(4\ell_0, 0) \leq L' \leq L$, objects $C^L_{L'}$ (again, either creatures $C^L_{L'} \subseteq q(\alpha_{L'}, L')$ or tuples of creature segments $C^L_{L'} \subseteq q(\mathrm{ct}, \mathrm{segm}(L'))$) and subsets of possibilities $B^L_{L'}$.

Step 1a: To start the induction, for $L' = L$ we set $C^L_L := q(\alpha_L, L)$ respectively $C^L_L := q(\mathrm{ct}, \mathrm{segm}(L))$. We let $B^L_L$ be the set of all $\eta \in \mathrm{poss}(q, {<}L)$ such that $q \wedge \eta$ essentially decides $\dot{\tau}$.

Step 1b: We continue the induction downwards on the relevant heights $L'$ with $(4\ell_0, 0) \leq L' < L$. We construct $C^L_{L'}$ and $B^L_{L'}$ such that the following holds:

- $C^L_{L'}$ is a strengthening of $q(\alpha_{L'}, L')$ respectively $q(\mathrm{ct}, \mathrm{segm}(L'))$ such that the corresponding norm decreases by at most 1.
- $B^L_{L'}$ is a subset of $\mathrm{poss}(q, {<}L')$ such that for each $\eta \in B^L_{L'}$ and each $x \in C^L_{L'}$, we have $\eta^\frown x \in B^L_{L'^+}$, and analogously for each $\eta \in \mathrm{poss}(q, {<}L') \setminus B^L_{L'}$ and each $x \in C^L_{L'}$, we have $\eta^\frown x \notin B^L_{L'^+}$. (We will call this property "homogeneity".) Since we only concern ourselves with relevant heights, $B_{L'^+}$ might not be explicitly defined by this process – if not, just take the smallest relevant height $L''$ above $L'$ and cut off the elements of $B^L_{L''}$ at height $L'^+$ to get $B^L_{L'^+}$.

Just as in the case of the proof of rapid reading in Theorem 6.6, we can find such objects using bigness:

- Define $L''$ to be the smallest relevant height above $L'$.
- By induction, there is a function $F$ mapping each $\eta \in \mathrm{poss}(q, {<}L'')$ to $\{\in B^L_{L''}, \notin B^L_{L''}\}$.
- We thin out $q(\alpha_{L'}, L')$ to $C^L_{L'}$, decreasing the norm by at most 1, such that for each $\nu \in \mathrm{poss}(q, {<}L')$, each extension of $\nu$ compatible with $C^L_{L'}$ has the same $F$-value $F^*(\nu)$.
- This in turn defines $B^L_{L'}$.

Step 2: We perform a downwards induction as in step 1 (always in the original $q$) from each relevant height $L$ above $(4\ell_0, 0)$. Given a relevant height $K$ such that $(4\ell_0, 0) \leq K < L$ and $\eta \in \mathrm{poss}(q, {<}K)$, and given that $q \wedge \eta$ essentially decides $\dot\tau$ and that $\eta' \in \mathrm{poss}(q, {<}L)$ extends $\eta$, it is clear that $q \wedge \eta'$ also essentially decides $\dot\tau$. We thus have:

$(*_6)$ If $q \wedge \eta$ essentially decides $\dot\tau$ for $\eta \in \mathrm{poss}(q, {<}K)$, then $\forall\, L > K : \eta \in B^L_K$.

Step 3: We now show the converse, namely:

$(*_7)$ Whenever $\eta \in B^L_{L'}$ for some relevant height $L$ with $L' = (4\ell_m, 0) \leq L$ (for some $m$), then $q \wedge \eta$ essentially decides $\dot\tau$.

To prove Eq. $(*_7)$, derive a condition $r$ from $q$ by using $\eta$ as the trunk and replacing creatures respectively tuples of creature segments at relevant heights $K$ (with $L' \leq K \leq L$) with $C^L_K$. Now, since all $\eta' \in \mathrm{poss}(r, {<}L) \subseteq \mathrm{poss}(q, {<}L)$ are in $B^L_L$, all $q \wedge \eta' \geq^* r \wedge \eta'$ essentially decide $\dot\tau$, and consequently, so does $r$. Noting that $\|r(\mathrm{nm}, 4k)\|_{\mathrm{nm},4k} > 0$ for all $4k \geq 4\ell_m$, we can use Eq. $(*_5)$ to get that $q \wedge \eta$ essentially decides $\dot\tau$.

Hence, to show that $q$ essentially decides $\dot\tau$, by Eq. $(*_7)$ it suffices to show that for all $\eta \in \mathrm{poss}(q, {<}(4\ell_0, 0))$ there is a height $L$ such that $\eta \in B^L_{(4\ell_0,0)}$.

Step 4: As in Theorem 6.6, we choose an "infinite branch" $\langle C^*_K, B^*_K \mid K \text{ relevant}\rangle$. (Recall that this means that for each height $L_0$, there is some $L > L_0$ such that, for all $K \leq L_0$, $(C^L_K, B^L_K) = (C^*_K, B^*_K)$.) By replacing the creatures and tuples of creature segments of $q$ at relevant heights $K$ with $C^*_K$, we obtain a condition $q^*$.

Step 5: To show that $q$ essentially decides $\dot\tau$, we thus have to show (as noted in Step 3) that $\eta \in B^*_{(4\ell_0,0)}$ for all $\eta \in \mathrm{poss}(q, {<}(4\ell_0, 0)) = \mathrm{poss}(q^*, {<}(4\ell_0, 0))$.

Fix any such $\eta$. Find an $r \leq q^* \wedge \eta$ deciding $\dot{\tau}$. Without loss of generality, for some $m$, $\operatorname{trklgth}(r) = 4\ell_m$ and $\|r(\mathrm{nm}, 4k)\|_{\mathrm{nm},4k} > 0$ for all $4k \geq 4\ell_m$. Let $\eta' \supseteq \eta$ be the trunk of $r$ restricted to $\operatorname{supp}(q, (4\ell_m, 0))$, which ensures $\eta' \in \operatorname{poss}(q, <(4\ell_m, 0))$ and $r \leq q \wedge \eta'$. By Eq. ($*_5$), $q \wedge \eta'$ already essentially decides $\dot{\tau}$.

Now pick some relevant $L > (4\ell_m, 0)$ such that $(C^L_K, B^L_K) = (C^*_K, B^*_K)$ for all relevant $K \leq (4\ell_m, 0)$. According to Eq. ($*_6$), $\eta' \in B^*_K$ and by homogeneity $\eta \in B^*_{(4\ell_0,0)}$ (since $\eta'$ is an extension of $\eta$). Hence by Eq. ($*_7$), $q \wedge \eta$ also essentially decides $\dot{\tau}$, which completes the proof. □

**Observation 7.5.** Given a condition $p^*$, a $p^*$-agreeable level $4\ell^*$, a number $M^*$ and a name $\dot{\tau}$ as in Lemma 7.4, and additionally a finite subset $F$ of the set of $A_{\mathrm{nn}} \cup A_{\mathrm{cn}} \cup \bigcup_\xi A_\xi$, we can first find an agreeable level $\ell^{**} \geq \ell^*$ such that for each $\alpha \in F$, there is a creature with index $\alpha$ at a level below $4\ell^{**}$ with norm $\geq M^*$ and such that there is also a segment below $4\ell^{**}$ with norm $\geq M^*$.

Then applying Lemma 7.4 after replacing $\ell^*$ with $\ell^{**}$ will now yield a condition $q$ (as in the lemma) which additionally has the property that for all indices in $F$, there is a creature below $4\ell^{**}$ with large norm, and analogously for the ct-part.

We can now use the preceding lemma to prove continuous reading.

**Theorem 7.6.** *Let $\dot{r}$ be a $\mathbb{Q}$-name for an element of* $\mathrm{Ord}^\omega$ *in $V$ and $p \in \mathbb{Q}$. Then there is a $q \leq p$ continuously reading $\dot{r}$.*

*Proof.* We will iteratively construct conditions $p_n$ in a similar way as in Part 2 of Lemma 7.4. Given $p^*, \ell^*, M^*, \dot{\tau}$ as in Lemma 7.4, we will denote the condition resulting from the application of that lemma by $s(p^*, \ell^*, M^*, \dot{\tau})$.

Set $p_{-1} := p$. Let $p_0 := s(p_{-1}, \ell_{-1}, 1, \dot{r}(0))$, where $\ell_{-1}$ is the minimal $\ell$ such that

- $4\ell$ is $p_{-1}$-agreeable,
- $4\ell \geq \operatorname{trklgth}(p_{-1})$, and
- $\|p_{-1}(\mathrm{nm}, 4k)\|_{\mathrm{nm},4k} \geq 2$ for all $4k \geq 4\ell$.

Given $p_n$ and $\ell_{n-1}$ such that $4\ell_{n-1}$ is $p_n$-agreeable, $p_n$ essentially decides $\dot{r}\restriction_{\{0,\dots,n\}}$ and $\|p_n(\mathrm{nm}, 4k)\|_{\mathrm{nm},4k} \geq n+1$ for all $4k \geq 4\ell_{n-1}$ (which is evidently true for $n = 0$), we define $p_{n+1}$ and $\ell_n$ as follows:

- Let $\ell_n > \ell_{n-1}$ be the minimal $\ell$ such that:
  - $4\ell_n$ is $p_n$-agreeable,
  - $p_n$ decides $\dot{r}(n)$ (or, equivalently, $\dot{r}\restriction_{\{0,\dots,n\}}$) below $(4\ell, 0)$,
  - $\|p_n(\mathrm{nm}, 4k)\|_{\mathrm{nm},4k} \geq n + 2$ for all $4k \geq 4\ell$, and
  - for each $\alpha \in \operatorname{supp}(p_n, (4\ell_{n-1}, 0)) \setminus A_{\mathrm{nm}}$ of type t, there is a height $L$ with $(\ell_{n-1}, J_{4\ell_{n-1}} - 1) < L < (4\ell_n, 0)$ such that $\|p_n(\alpha, L)\|_{\mathrm{t},L} \geq n$.
- Let $p_{n+1} := s(p_n, \ell_n, n + 1, \dot{r}(n + 1))$.

Lemma 7.4 ensures that $p_{n+1} \leq p_n$, that $p_{n+1}$ essentially decides $\dot{r}(n+1)$ (and thus $\dot{r}\restriction_{\{0,\dots,n+1\}}$) and that it fulfils $\|p_{n+1}(\mathrm{nm}, 4k)\|_{\mathrm{nm},4k} \geq n + 1$ for all $4k \geq 4\ell_n$.

Similar to Part 2 of Lemma 7.4, $\langle p_n \mid n < \omega\rangle$ is a descending sequence of conditions converging to a condition $q \in \mathbb{Q}$. By construction, $q$ continuously reads $\dot{r}$. □

Theorem 7.6 and Theorem 6.6 taken together show that for any $p \in \mathbb{Q}$ and any $\mathbb{Q}$-name $\dot{r}$ for a real, there is a $q \leq p$ rapidly reading $\dot{r}$. Even more important are the following consequences of the previous two sections, which prove the first, easier parts of this paper's main theorem, Theorem 1.1:

**Lemma 7.7.** *$\mathbb{Q}$ satisfies the following variant of the finite version of Baumgartner's axiom A: There is a family $(\leq_{F,n})$ of binary relations, indexed by finite subsets $F \subseteq A$ and natural numbers $n$, such that the following holds:*

- *$q \leq_{F,n} p$ implies $q \leq p$.*
- *Whenever $\langle p_n \mid n \in \omega \rangle$ is a fusion sequence based on a sequence $\langle F_n \mid n \in \omega \rangle$ (i. e. the conditions $p_n$ are decreasing, the finite sets $F_n$ increasing, for each $n$ we have $p_{n+1} \leq_{F_n,n} p_n$, and the union $\bigcup_n F_n$ covers $\bigcup_n \operatorname{supp}(p_n)$), then the sequence $\langle p_n \mid n \in \omega \rangle$ has a canonical weakest lower bound (to which the $p_n$ converge in a natural sense).*
- *For all $p$, $n$, $F$ and all $\mathbb{Q}$-names $\dot{\alpha}$ of an ordinal, there is a condition $q \leq_{F,n} p$ and a finite set $E$ such that $q \Vdash \dot{\alpha} \in E$.*

*Hence $\mathbb{Q}$ is proper and $\omega^\omega$-bounding. Assuming* CH *in the ground model, $\mathbb{Q}$ moreover preserves all cardinals and cofinalities.*

*Proof.* Define the relations $\leq_{F,n}$ as follows: $q \leq_{F,n} p$ if there is some $\ell \geq n$ such that

- $q \leq p$;
- $4\ell$ is $p$-agreeable;
- $p$ and $q$ are identical below $(4\ell, 0)$ on $F$;
- $\|q(\mathrm{nm}, 4k)\|_{\mathrm{nm},4k} > n$ for all $4k \geq 4\ell$;
- there is a ct-segment of $p$ below $4\ell$ with norm $\geq n$; and
- for all $\alpha \in F$ which are "modular", i. e. neither in $A_{\mathrm{ct}}$ nor in $A_{\mathrm{nm}}$, there is a creature with index $\alpha$ below level $4\ell$ with norm $\geq n$.

Note that for any sequence $p_n \geq_{F_n,n} p_{n+1} \geq_{F_{n+1},n} \dots$ with witnesses $\ell_n$, $\ell_{n+1}$, for each index $\alpha \in F_n$ there will be a creature of norm $\geq n$ which will stay unchanged throughout the fusion sequence, and will hence also appear in the fusion limit of this sequence; hence any fusion sequence has a limit.

Given a condition $p$ and name $\dot{\tau}$ of a sequence of ordinals, one can inductively construct a fusion sequence $\langle p_n \mid n \in \omega \rangle$ (using Lemma 7.4 and Observation 7.5 in every step, as well as the usual bookkeeping to ensure that the union $\bigcup_n F_n$ covers $\bigcup_n \operatorname{supp}(p_n)$) that converges to a stronger condition $q$ such that there is a family $\langle E_n \mid n \in \omega \rangle$ of finite sets $q \Vdash \forall n \in \omega \colon \dot{\tau}(n) \in E_n$.

From Lemma 4.18 we conclude that $\mathbb{Q}$ preserves all cardinals and cofinalities $\geq \aleph_2$, and since it is proper, it also preserves $\aleph_1$. This proves the "moreover" part of Theorem 1.1. □

**Lemma 7.8.** *Assuming* CH *in the ground model, in the extension $\mathfrak{d} = \aleph_1$ and $\operatorname{cov}(\mathcal{N}) = \aleph_1$.*

*Proof.* Since $\mathbb{Q}$ is $\omega^\omega$-bounding, it forces $\mathfrak{d}$ to be $\aleph_1$. To prove the second part of the statement, we show that each new real is forced to be contained in a ground model null set, so the $\aleph_1$ many Borel null sets of the ground model cover the reals (in other words, $\mathbb{Q}$ adds no random reals) and hence $\mathrm{cov}(\mathcal{N})$ is forced to be $\aleph_1$.

Let $\dot r$ be a $\mathbb{Q}$-name for a real and $p \in \mathbb{Q}$. Let $q \leq p$ read $\dot r$ rapidly, which means that for each $L \in$ **heights**, $\dot r{\restriction}_{n^R_{<L}}$ is determined by $\eta \in \mathrm{poss}(q, {<}L)$; let $X^q_L$ be the set of possible values of $\dot r{\restriction}_{n^R_{<L}}$. For notational simplicity, consider only heights $\ell$ of the form $(4k, 0), 4k+1, 4k+2, 4k+3$ and identify $(4k, 0)$ with $4k$. Then it follows that $|X^q_\ell| \leq n^P_{<\ell} < n^R_{<\ell} < 2^{n^R_{<\ell}}/\ell$, where the last inequality holds by our general requirement on the $n^R_{<\ell}$. This means that the relative size of $X^q_\ell$ is bounded by $1/\ell$ and hence $\langle X^q_\ell \mid \ell < \omega \rangle$ can be used to define the ground model null set $N_q := \{s \in 2^\omega \mid \forall\, \ell < \omega \colon s{\restriction}_{n^R_\ell} \in X^q_\ell\}$. By definition, $q \Vdash \dot r \in N_q$. ☐

This proves (M1) of Theorem 1.1.

**Lemma 7.9.** *In the extension,* $2^{\aleph_0} = \kappa_{\mathrm{ct}}$.

*Proof.* If $\alpha, \beta \in A_{\mathrm{ct}}$ are distinct, then the reals $\dot y_\alpha$ and $\dot y_\beta$ are forced to be different, hence there are at least $\kappa_{\mathrm{ct}}$ many reals in the extension. But every real in the extension is read continuously by Theorem 7.6, hence by Lemma 6.5 there are at most $\kappa_{\mathrm{ct}}^{\aleph_0} = \kappa_{\mathrm{ct}}$ many reals in the extension. ☐

This proves (M6) of Theorem 1.1. It remains to prove points (M2)–(M5) of Theorem 1.1, which we will do in the following sections.

## 8. $\mathrm{cof}(\mathcal{N}) \leq \kappa_{\mathrm{cn}}$

To show $\mathrm{cof}(\mathcal{N}) \leq \kappa_{\mathrm{cn}}$, we prove that $\mathbb{Q}$ has the Laver property over the intermediate forcing poset

$$\mathbb{Q}_{\text{non-ct}} := \left( \prod_{\mathrm{t} \,\in\, \mathbf{types}_{\mathrm{modular}}} \mathbb{Q}_{\mathrm{t}}^{\kappa_{\mathrm{t}}} \right) \times \mathbb{Q}_{\mathrm{nm},\, \kappa_{\mathrm{nm}}}$$

(and hence also the Sacks property, since it is $\omega^\omega$-bounding). We will use the same equivalent formulation as in [FGKS17, Lemmas 6.3.1–2], namely, we will prove:

**Lemma 8.1.** *Given a condition* $p \in \mathbb{Q}$*, a name* $\dot r \in 2^\omega$ *and a function* $g \colon \omega \to \omega$ *in* $V$*. Then there is a* $q \leq p$ *and a name* $\dot T \subseteq 2^{<\omega}$ *for a leafless tree such that:*

- $q$ *reads* $\dot T$ *continuously not using any indices in* $A_{\mathrm{ct}}$*,*
- $q \Vdash \dot r \in [\dot T]$*, and*
- $\left|\dot T \cap 2^{g(n)}\right| < n + 2$ *for all* $n < \omega$*.*

*Proof.* We first note that we can increase $g$ without loss of generality, since if $g_1(n) \leq g_2(n)$ for all $n$ and $\dot T$ witnesses the lemma for $g_2$, then the same $\dot T$ also witnesses the lemma for $g_1$.

We can also assume without loss of generality that $p$ is modest and rapidly reads $\dot r$, i. e. $\mathrm{poss}(p, {<}L)$ determines $\dot r{\restriction}_{n^R_{<L}}$ for all heights $L$. Considering this, we can find a strictly increasing sequence of segment-initial heights $L_n$ (i. e. $\min(\mathrm{segm}(L_n)) = L_n$) such that $g(n) = n^R_{<L_n}$ for all $n < \omega$ (increasing $g$ when necessary).

Hence, each $\eta \in \mathrm{poss}(p, {<}L_n)$ defines a value $\dot{R}^n(\eta)$ for $\dot{r}\restriction_{g(n)}$. We split each $\eta$ into two components, $\eta_{\mathrm{ct}}$ and $\eta_{\mathrm{rmdr}}$ (i. e. the non-ct remainder). If we fix the $\eta_{\mathrm{ct}}$ component of $\eta$, then $\dot{R}^n(\cdot, \eta_{\mathrm{ct}})$ is a name not depending on the ct component, i. e. not using any indices in $A_{\mathrm{ct}}$. (More formally: Given an $\eta_{\mathrm{rmdr}}$ compatible with the generic filter such that $(\eta_{\mathrm{rmdr}}, \eta_{\mathrm{ct}}) = \eta \in \mathrm{poss}(p, {<}L_n)$, $\dot{R}^n(\eta_{\mathrm{rmdr}}, \eta_{\mathrm{ct}})$ evaluates to $\dot{R}^n(\eta)$.)

We will now construct a stronger condition $q$ and an increasing sequence $\langle i_n \mid n < \omega\rangle$ of natural numbers (such that each $L_{i_n}$ is a segment-initial height) with the following properties: Given some $i_{n+1}$, let $i_n < m \leq i_{n+1}$ and $\eta \in \mathrm{poss}(q, {<}L_{i_{n+1}})$. Such an $\eta$ extends a unique $\eta^m$ in the set of possibilities $\mathrm{poss}(q, {<}L_m)$ cut off at height $L_m$, which we call $\mathrm{poss}^\dagger(q, {<}L_m)$. Restricting this $\eta^m$ to the ct component yields $\eta^m_{\mathrm{ct}} := \eta^m\restriction_{A_{\mathrm{ct}}}$.[29] Then $q \wedge \eta$ forces the name $\dot{R}^m(\cdot, \eta^m_{\mathrm{ct}})$ to be evaluated to $\dot{r}\restriction_{g(m)}$, and hence $q$ forces $\dot{r}\restriction_{g(m)}$ to be an element of

$$\dot{T}^m := \{\dot{R}^m(\cdot, \eta^m_{\mathrm{ct}}) \mid \eta \in \mathrm{poss}(q, {<}L_{i_{n+1}})\},$$

which is a name not using any indices in $A_{\mathrm{ct}}$. It thus suffices to show that there are few such $\eta^m_{\mathrm{ct}}$, i. e. that letting $P_m := \{\eta^m_{\mathrm{ct}} \mid \eta \in \mathrm{poss}(q, {<}L_{i_{n+1}})\}$, for all $m < \omega$ we have $|P_m| < m + 2$.

The condition $q$ will have the same support as $p$. On $\mathrm{supp}(p) \setminus A_{\mathrm{ct}}$, we define $q$ to be equal to $p$. Hence we now only have to define $q$ on $\mathrm{supp}(p) \cap A_{\mathrm{ct}}$. We will inductively construct the sequence $\langle i_n\rangle$ and the new condition $q(\mathrm{ct})$ below $L_{i_n}$, and show that $|P_m| < m + 2$ holds for all $m \leq i_n$. To begin the induction, let $i_0 = 0$ and let $q(\mathrm{ct})$ below $L_0$ be identical to some arbitrary possibility in $\mathrm{poss}(p(\mathrm{ct}), {<}L_0)$, giving us $|P_{i_0}| = 1$.

By way of induction hypothesis, assume we already have $i_n$, $q$ is defined up to $L_{i_n}$ and $|P_m| < m+2$ holds for all $m \leq i_n$. (By our choice of $i_0 = 0$, all this is fulfilled for $n = 0$.) Keep in mind that each $L_i$ is the initial height in a segment of the frame of $p(\mathrm{ct})$.

Step 1: Let $\Sigma := \mathrm{supp}(p, \mathrm{ct}, L_{i_n}) \cap A_{\mathrm{ct}}$ and let $k$ be such that $L_{i_n} = 4k + 3$. (Note that hence $|\Sigma| < k$, though this is not important to this proof.) Let $c$ be minimal such that $\mathrm{nor}^{n^B_{L_{i_n}}, k}_{\mathrm{Sacks}}(c) = n$.[30] Let $i' := (i_n + 2) \cdot c^{|\Sigma|}$. For each $\alpha \in \Sigma$, find $L^\alpha > L_{i'}$ (with $L^\alpha \neq L^\beta$ for $\alpha \neq \beta$) such that $\mathrm{nor}^{n^B_{L_{i_n}}, k}_{\mathrm{Sacks}}(p(\alpha, L^\alpha)) \geq n$. Finally, let $i_{n+1} > i'$ be minimal such that $L^\alpha < L_{i_{n+1}}$ for all $\alpha \in \Sigma$.

Step 2: We define $q(\mathrm{ct})$ from $L_{i_n}$ up to (but excluding) $L_{i_{n+1}}$ as follows: For each $\alpha \in \Sigma$, we take $p(\alpha, L^\alpha)$ and shrink it such that $\mathrm{nor}^{n^B_{L_{i_n}}, k}_{\mathrm{Sacks}}(q(\alpha, L^\alpha)) = n$. For all other heights $L \in \mathbf{heights}_{\mathrm{ct}}$ with $L_{i_n} \leq L < L_{i_{n+1}}$, replace $p(\alpha, L)$ with an arbitrary singleton to get $q(\alpha, L)$. In particular, this means that for $L_{i_n} \leq L \leq L_{i'}$, $p(\alpha, L)$ is a singleton for each $\alpha \in \Sigma$.

[29] Note that these $\eta^m_{\mathrm{ct}}$ – i. e. elements of $\mathrm{poss}^\dagger(q, {<}L_m)$ restricted to $A_{\mathrm{ct}}$ – are still *not* possibilities in $\mathrm{poss}(q, \mathrm{ct}, {<}L_m)$, however, since such possibilities would actually go up to height $\max(\mathrm{segm}(L_{m-1}))$ instead of ending at height $L_m$. (They are, however, still possibilities in the "moral" sense, i. e. they are initial segments of the generic real $\dot{y}$.) Please excuse this minor abuse of notation. It makes sense when one considers that in the frame of $p$, $L_m$ *used* to be a segment-initial height, even if it no longer is in the frame of $q$.

[30] Careful readers may point out that the Sacks norm was defined only for sets in Definition 3.8; we ask these readers to read "$c$" as "$\{0, \ldots, c-1\}$" instead.

For the frame of $q(\mathrm{ct})$, take the segments in the frame of $p(\mathrm{ct})$ starting at (the segment starting with) $L_{i_n}$ and going up to, and including, the (segment ending with the) $\mathbf{heights}_{\mathrm{ct}}$-predecessor $L'$ of $L_{i_{n+1}}$; merge all of them to form a single segment in the frame of $q(\mathrm{ct})$.

Step 3: For those indices $\alpha$ in $\mathrm{supp}(p, \mathrm{ct}, L')$ which are outside of $\Sigma$ (i. e. those which enter the support of $p(\mathrm{ct})$ strictly above $L_{i_n}$ and up to $L'$), also choose arbitrary singletons to get a trivial $q(\alpha, \langle L_{i_n}, \dots, L' \rangle)$.[31] (Such indices will be in the support of $q(\mathrm{ct})$ from $L_{i_{n+1}}$ onwards.)

Step 4: We now just have to prove that $|P_m|$ is sufficiently small up to (and including) $i_{n+1}$. First, let $i_n < m < i'$; for such $m$, we did not add any possibilities to $q$ (as all new creature segments consist of singletons up to that height), so $|P_m| = |P_{i_n}| < i_n + 2 < m + 2$. Now consider $i' \leq m \leq i_{n+1}$. For each $\alpha \in \Sigma$, the number of possibilities in $q(\alpha, \langle L_{i_n}, \dots, L' \rangle)$ is exactly $c$. By the induction hypothesis we already know that $|P_{i_n}| < i_n + 2$, and due to the choice of $i'$, we altogether have

$$|P_m| \leq |P_{i_n}| \cdot c^{|\Sigma|} < (i_n + 2) \cdot c^{|\Sigma|} = i' < m + 2$$

and we are done with the induction. □

Having proved this, we now know that $\mathbb{Q}$ has the Sacks property over the intermediate forcing poset $\mathbb{Q}_{\text{non-ct}}$. By [BJ95, Theorem 2.3.12] (later restated as Theorem 10.3 in section 10, where we will use it a bit more extensively), this is equivalent to the fact that any null set in the model obtained by forcing with the entire $\mathbb{Q}$ is contained in a null set of the model obtained by forcing with $\mathbb{Q}_{\text{non-ct}}$, and hence we have shown that $\mathbb{Q} \Vdash \mathrm{cof}(\mathcal{N}) \leq \kappa_{\mathrm{cn}}$ by Lemma 6.5.

We will show $\mathbb{Q} \Vdash \mathrm{cof}(\mathcal{N}) \geq \kappa_{\mathrm{cn}}$ a bit later.

## 9. $\mathrm{non}(\mathcal{M}) = \kappa_{\mathrm{nm}}$

The following proof does not use any specifics of the creatures and possibilities; it only requires that $\mathbb{Q}_{\mathrm{nm}, \kappa_{\mathrm{nm}}}$ is the only part of the forcing poset involving a $\liminf$ construction.

**Lemma 9.1.** *The set of all reals that can be read continuously only using indices in $A_{\mathrm{nm}}$ is not meagre.*

*Proof.* Let $\dot{M}$ be a $\mathbb{Q}$-name for a meagre set. We can find $\mathbb{Q}$-names of nowhere dense trees $\dot{T}_n \subseteq 2^{<\omega}$ such that $\dot{M} \subseteq \bigcup_{n<\omega} [\dot{T}_n]$ is forced. We will show that there is a $\mathbb{Q}$-name for a real $\dot{r}$ which is continuously read only using indices in $A_{\mathrm{nm}}$ such that $\dot{r} \notin \dot{M}$; hence, the set of all such reals cannot be meagre.

First note that since $\mathbb{Q}$ is $\omega^\omega$-bounding and all $\dot{T}_n$ are nowhere dense, for each $n < \omega$, there is a ground model function $f_n \colon \omega \to \omega$ such that the following holds: For each $\rho \in 2^x$, there is a $\rho' \in 2^{f_n(x)}$ such that $\rho \subseteq \rho'$ and such that $\rho' \notin \dot{T}_n$ is forced. We find this family of functions as follows: Clearly, there is a family of names $\langle \dot{g}_n \rangle$ with this property; let $\dot{g}_\infty$ be a name for a function dominating each $\dot{g}_n$ above $n$ and let $f_\infty$ be a ground model function dominating $\dot{g}_\infty$ everywhere. As we may without loss of generality assume that all mentioned functions are strictly monotone, we can define the desired functions $f_n$ by $f_n(k) := f_\infty(\max(k, n))$.

[31] Recall that $q(\alpha, \bar{L})$ was defined in Definition 4.2.

We fix some $p \in \mathbb{Q}$ forcing the previously mentioned properties of $\dot{M}$ and $\langle \dot{T}_n \rangle$ and continuously reading all $\dot{T}_n$ (which is possible per Theorem 7.6). We will construct (in the ground model) $q \leq p$ and a name $\dot{r}$ for a real continuously read by $q$ only using indices in $A_{\mathrm{nm}}$ such that $q \Vdash \dot{r} \notin \dot{M}$.

We will define $q$ inductively as the limit of a fusion sequence $q_i$. Assume we have already defined $q$ in the form of a condition $q_i$ up to some $q_i$-agreeable $4k_i$, and that we have an $x_i < \omega$ and a $\mathbb{Q}$-name $\dot{z}_i$ for an element of $2^{x_i}$ such that $\dot{z}_i$ is decided by $\mathrm{poss}(q_i{\restriction}_{A_{\mathrm{nm}}}, {<}(4k_i, 0))$. (The real $\dot{r}$ will be defined as the increasing union of the $\dot{z}_i$.) Finally, assume that $q_i$ already forces $\dot{z}_i \notin \dot{T}_0 \cup \dot{T}_1 \cup \ldots \cup \dot{T}_{i-1}$. The idea is to now extend $\dot{z}_i$ to a longer name $\dot{z}_{i+1}$ which is forced by $q_{i+1}$ to avoid $\dot{T}_i$, as well.

To that end, enumerate $\mathrm{poss}(q_i, {<}(4k_i, 0))$ as $(\eta^0, \eta^1, \ldots, \eta^{m-1})$. Set $4k_i^0 := 4k_i$, $x_i^0 := x_i$, $\dot{z}_i^0 := \dot{z}_i$. By induction on $j$, $0 \leq j < m$, we deal with $\eta^j$: Assume we are given a name $\dot{z}_i^j$ for an element of $2^{x_i^j}$ that is decided by $\mathrm{poss}(q_i, {<}(4k_i^j, 0))$, and that we have already constructed a condition $q_i^j$ such that

- $q_i^j \leq q_i$,
- $q_i^j$ is identical to $q_i$ below $4k_i$,
- $4k_i^j$ is $q_i^j$-agreeable,
- between $4k_i^0 = 4k_i$ and $4k_i^j$, all creatures in $q_i^j{\restriction}_{A \setminus A_{\mathrm{nm}}}$ are singletons,
- $q_i^j \Vdash \dot{z}_i^{j-1} \subseteq \dot{z}_i^j$,
- $q_i^j \Vdash |\dot{z}_i^j| = x_i^j$, and
- $q_i^j \wedge \eta^\ell \Vdash \dot{z}_i^j \notin \dot{T}_i$ for all $0 \leq \ell < j$.

Let $x_i^{j+1} := f_n(x_i^j)$ and choose some $q_i^j$-agreeable $4k_i^{j+1}$ above $4k_i^j$ which is big enough for $p$ to determine $\dot{X} := \dot{T}_i{\restriction}_{x_i^{j+1}}$, i.e. there is some function $F$ from $\mathrm{poss}(p, {<}(4k_i^{j+1}, 0))$ to possible values of $\dot{X}$ (a consequence of continuous reading). We now define $q_i^{j+1}$ as follows: Below $4k_i^j$ and above $4k_i^{j+1}$, $q_i^{j+1}$ is identical to $q_i^j$. Between $4k_i^j$ and $4k_i^{j+1}$, we just leave $q_i^j(\mathrm{nm})$ as it is; and on $A \setminus A_{\mathrm{nm}}$, we just choose arbitrary singletons within $q_i^j$ to get $q_i^{j+1}$.

We now briefly consider the name $\dot{X}$: A possibility $\nu \in \mathrm{poss}(p, {<}(4k_i^{j+1}, 0))$ consists of

- $U$, the part below $4k_i^0$,
- $V$, the part in $A \setminus A_{\mathrm{nm}}$ between $4k_i^0$ and $4k_i^{j+1}$, and
- $W$, the part in $A_{\mathrm{nm}}$ between $4k_i^0$ and $4k_i^{j+1}$.

So we can write $\dot{X} = F(U, V, W)$. Under the assumption that the generic follows $\eta^j$ ($U = \eta^j$) and the singleton values of $q_i^{j+1}$ on $A \setminus A_{\mathrm{nm}}$ ($V = V^j$ for some fixed $V^j$), then there is a name $\dot{X}' = F(\eta^j, V^j, \cdot)$ depending only on indices in $A_{\mathrm{nm}}$ such that $q_i^{j+1} \Vdash \dot{X} = \dot{X}'$.

Recall that $\dot{z}_i^j$ is already determined by the nm-part of $\eta^j$ and that already $p$ forces that there is some extension $z' \in 2^{x_i^{j+1}}$ of that value of $\dot{z}_i^j$ such that $z' \notin \dot{X}'$. By picking (in the ground model) for each possible choice of $W$ some $z'(W) \in 2^{x_i^{j+1}} \setminus F(\eta_j, V^j, W)$ extending $\dot{z}_i^j$, we can define the name $\dot{z}_i^{j+1} := z'(\cdot)$ which depends only on indices in $A_{\mathrm{nm}}$ and is determined below $4k_i^{j+1}$. By construction, we have that $q_i^{j+1} \Vdash \dot{z}_i^j \subseteq \dot{z}_i^{j+1}$ and $q_i^{j+1} \wedge \eta^j \Vdash \dot{z}_i^{j+1} \notin \dot{T}_i$.

Repeating this construction for all $j$, $0 \leq j < m$, finally define $\dot{z}_{i+1} := \dot{z}_i^m$ and $x_{i+1} := x_i^m$ and let $4k_{i+1}$ be such that

- $4k_{i+1}$ is above $4k_i^m$,
- $4k_{i+1}$ is $q_i^m$-agreeable, and
- for all $\alpha \in \mathrm{supp}(q_i^m, (4k_i^m, 0)) \setminus A_{\mathrm{nm}}$ of type t, there is a height $(4k_i^m, 0) < L < (4k_{i+1}, 0)$ such that $\|p(\alpha, L)\|_{\mathrm{t},L} \geq i$.

Define $q_{i+1}$ to be equal to $q_i^m$ below $4k_i^m$ and equal to $p$ above $4k_i^m$. By our choice of $4k_{i+1}$, we have ensured that the lim sup part of the fusion condition $q := \bigcap_{i<\omega} q_i$ will actually be a condition (the lim inf part trivially is). By the construction of $q_{i+1}$, we have ensured that $q_{i+1}$ forces that $\dot{r} := \bigcup_{i<\omega} \dot{z}_i$ avoids $\dot{T}_0 \cup \dot{T}_1 \cup \ldots \cup \dot{T}_i$, and hence $q$ forces that $\dot{r}$ avoids $\dot{M}$. Finally, by the construction of the $\dot{z}_i$, they are continuously read by $q$ only using indices in $A_{\mathrm{nm}}$, and so is their union $\dot{r}$. □

**Corollary 9.2.** *$\mathbb{Q}$ forces* $\mathrm{non}(\mathcal{M}) \leq \kappa_{\mathrm{nm}}$.

*Proof.* By Lemma 6.5, the non-meagre set from Lemma 9.1 has size at most $\kappa_{\mathrm{nm}}$, and hence we have $\mathbb{Q} \Vdash \mathrm{non}(\mathcal{M}) \leq \kappa_{\mathrm{nm}}$. □

To prove $\mathrm{non}(\mathcal{M}) \geq \kappa_{\mathrm{nm}}$, we first define some meagre sets in the extension. Recall that for $\alpha \in A_{\mathrm{nm}}$, the generic object $\dot{y}_\alpha$ is a $\mathbf{heights}_{\mathrm{nm}}$-sequence of objects in $\mathbf{POSS}_{\mathrm{nm},L} = 2^{I_L}$, or equivalently an $\omega$-sequence of 0s and 1s. We define a name for a meagre set $\dot{M}_\alpha$ as follows: A real $r \in 2^\omega$ is in $\dot{M}_\alpha$ iff for all but finitely many $k < \omega$, there is an $i_k$ such that $r{\restriction}_{I_{(4k,i_k)}} \neq \dot{y}_\alpha(4k, i_k)$, or equivalently

$$\dot{M}_\alpha := \bigcup_{n<\omega} \bigcap_{k \geq n} \{r \in 2^\omega \mid r{\restriction}_{I_{4k}} \neq \dot{y}_\alpha{\restriction}_{I_{4k}}\}$$

(abusing the notation by letting $I_{4k} := \bigcup_{i \in J_{4k}} I_{(4k,i)}$), whence it is clear that $\dot{M}_\alpha$ is indeed a meagre set.

By the choice of $n^R_{<L}$ for the nm case in Definition 4.8, if $p$ rapidly reads $\dot{r}$, then for any $L \in \mathbf{heights}_{\mathrm{nm}}$, $\dot{r}{\restriction}_{I_L}$ is decided $\leq L$. Also note that if the cell norm $\|x\|^{\mathrm{cell}}_L$ of some creature $x$ is at least 1, then it follows that $|x| > n^P_{<L}$.

**Lemma 9.3.** *Let $\dot{r} \in 2^\omega$ be a name for a real and let $p$ rapidly read $\dot{r}$ not using the index $\alpha \in A_{\mathrm{nm}}$. Then $p \Vdash \dot{r} \in \dot{M}_\alpha$.*

*Proof.* We first remark that it suffices to prove that there is an $s \leq p$ such that $s \Vdash \dot{r} \in \dot{M}_\alpha$. Assume that we have shown this, and also assume that $p$ does *not* force $\dot{r} \in \dot{M}_\alpha$; then there is a $q \leq p$ forcing the contrary, and $q$ still rapidly reads $\dot{r}$ not using the index $\alpha$. Since we can thus find an $s \leq q$ which *does* force $\dot{r} \in \dot{M}_\alpha$, we have arrived at the desired contradiction.

We only have to find $s \leq p$ forcing $\dot{r} \in \dot{M}_\alpha$. As a matter of fact, we will only have to modify $p$ in very few places to arrive at the desired condition $s$. Without loss of generality, assume that $\alpha \in \operatorname{supp}(p)$. Recall that by the definition of $\mathbb{Q}_{\mathrm{nm},\kappa_{\mathrm{nm}}}$, there is some $k_1$ such that for any $k \geq k_1$, $\|p(4k)\|_{\mathrm{nm},4k} \geq 1$; as a consequence, for each stacked creature $p(\alpha, 4k)$, there is at least one $i \in J_{4k}$ with $\|p(\alpha,(4k,i))\|^{\mathrm{cell}}_{(4k,i)} \geq 1$; for each $k \geq k_1$, we pick some such $i_k$.

Consider one of these $(4k, i_k) =: L$. We know that $\dot r\restriction_{I_L}$ is decided $\leq L$ by $p$ – and actually even below $L$, since by modesty (ii), there can be at most one index $\beta$ such that $p(L,\beta)$ is non-trivial, $\alpha$ already is such an index and $p$ reading $\dot r$ does not depend on the index $\alpha$. Since there are at most $n^P_{<L}$ many possibilities below $L$ in $p$, there can be at most $n^P_{<L}$ many possible values for $\dot r\restriction_{I_L}$, and since $|p(\alpha, L)| > n^P_{<L}$, there must be some $x_k \in p(\alpha, L)$ different from all possible values of $\dot r\restriction_{I_L}$ under the reading by $p$.

We define the condition $s$ by replacing each $p(\alpha,(4k,i_k))$ with the singleton $\{x_k\}$. It is clear that $s$ is still a condition, as we have at most reduced each stacked creature's norm in $p(\alpha)$ by 1, which does not negatively affect the lim inf norm convergence. By definition, $s \Vdash \dot r \in \dot M_\alpha$ as required, since $\dot r\restriction_{I_{(4k,i_k)}}$ is different from $\dot y_\alpha$ for all $k \geq k_1$. ☐

**Corollary 9.4.** $\mathbb{Q}$ *forces* $\operatorname{non}(\mathcal{M}) \geq \kappa_{\mathrm{nm}}$.

*Proof.* Fix a condition $p$, some $\kappa < \kappa_{\mathrm{nm}}$ (without loss of generality $\kappa \geq \aleph_1$ – otherwise, there is nothing to prove) and a sequence of names of reals $\langle \dot r_i \mid i \in \kappa\rangle$. We find some $\alpha \in A_{\mathrm{nm}}$ such that $p \Vdash \{\dot r_i \mid i \in \kappa\} \subseteq \dot M_\alpha$.

For each $i \in \kappa$, fix a maximal antichain $A_i$ below $p$ such that each $a \in A_i$ rapidly reads $\dot r_i$. Recall that $\mathbb{Q}$ is $\aleph_2$-cc by Lemma 4.18. Since $\kappa_{\mathrm{nm}} > \kappa$, $S := \bigcup_{i\in\kappa}\bigcup_{a\in A_i}\operatorname{supp}(a)$ has size $\kappa < \kappa_{\mathrm{nm}}$ and we can find an index $\alpha \in A_{\mathrm{nm}} \setminus S$. Each $a \in A_i$ rapidly reads $\dot r_i$ not using the index $\alpha$; so by the preceding lemma, for each $i$, each $a \in A_i$ forces $\dot r_i \in \dot M_\alpha$ and so does $p$ (since $A_i$ is predense below $p$), finishing the proof. ☐

This proves (M2) of Theorem 1.1.

## 10. $\mathfrak{c}_{f_\xi,g_\xi} = \kappa_\xi$

**Definition 10.1.** Given $f, g \in \omega^\omega$ going to infinity such that $0 < g < f$, we call $S := \langle S_k \mid k < \omega\rangle \in ([\omega]^{<\omega})^\omega$ an $(f,g)$-*slalom* if $S_k \subseteq f(k)$ and $|S_k| \leq g(k)$ for all $k < \omega$, or in shorter notation, if $S \in \prod_{k<\omega}[f(k)]^{\leq g(k)}$.

We say a family of $(f,g)$-slaloms $\mathcal{S}$ is $(f,g)$-*covering* if for all $h \in \prod_{k<\omega} f(k)$ there is an $S \in \mathcal{S}$ such that $h \in^* S$ (i. e. $h(k) \in S_k$ for all but finitely many $k < \omega$).[32]

We then define the cardinal characteristic $\mathfrak{c}_{f,g}$, sometimes also denoted by $\mathfrak{c}^\forall_{f,g}$ and referred to as one of two kinds of *localisation cardinals*, as

$$\mathfrak{c}_{f,g} := \min\left\{ |\mathcal{S}| \;\middle|\; \mathcal{S} \subseteq \prod_{k<\omega}[f(k)]^{\leq g(k)}, \forall\, x \in \prod_{k<\omega} f(k)\ \exists S \in \mathcal{S} \colon x \in^* S \right\},$$

the minimal size of an $(f,g)$-covering family.

[32] Equivalently, $h \in S$ would lead to the same results.

A simple diagonalisation argument shows that under the assumptions above, $\mathfrak{c}_{f,g}$ is always uncountable. [GS93, section 1] contains a few simple properties following from the definition, but the only one we will be interested in here is monotonicity, in the following sense:

**Fact 10.2.** If $f \leq^* f'$ and $g \geq^* g'$, then $\mathfrak{c}_{f,g} \leq \mathfrak{c}_{f',g'}$.

For the following proof of $\mathfrak{c}_{f,g} \leq \operatorname{cof}(\mathcal{N})$, we recall a result from [Bar84] (as presented in [BJ95, Theorem 2.3.12]):

**Theorem 10.3** (Bartoszyński). *Let $M \subseteq N$ be transitive models of* ZFC*. *The following are equivalent:*

*(i) Every null set coded in $N$ is covered by a Borel null set coded in $M$.*
*(ii) Every convergent series of positive reals in $N$ is dominated by a convergent series in $M$.*
*(iii) For every function $h \in \omega^\omega \cap N$ there is a slalom $S \in \mathcal{C} \cap M$ such that $h(k) \in S(k)$ for almost all $k$.*

In this theorem, $\mathcal{C}$ is defined as the set of all slaloms $S$ such that

$$\sum_{k \geq 1} \frac{|S_k|}{k^2} < \infty,$$

which does not directly relate to our cardinal characteristics, but very nearly so:

**Fact 10.4.** Consider the following:

- We extend the definition of $(f,g)$-slalom and $(f,g)$-covering to allow $f \in (\omega+1)^\omega$. Write

$$\mathfrak{c}_{\omega,g} := \min \left\{ |\mathcal{S}| \;\middle|\; \mathcal{S} \subseteq \prod_{k<\omega} [\omega]^{\leq g(k)}, \forall\, x \in \omega^\omega \;\exists\, S \in \mathcal{S} \colon x \in^* S \right\},$$

  i. e. identify $\omega$ with the constant $\omega$-valued function. By Fact 10.2, we then have that $\mathfrak{c}_{f,g} \leq \mathfrak{c}_{\omega,g}$.
- Note that "$\mathfrak{c}_{\omega,g} \leq \kappa$" actually is a thinly veiled statement about the Sacks property, in the sense that it simply says that there is a model $M$ of size $\kappa$ over which the universe $V$ has the Sacks property.
- Finally, recall the well-known fact that the statement of the Sacks property is independent of the specific slalom size used (since for any two slalom size functions, the statements can be converted into each other by a simple coding argument). Hence it is clear that for any $g, g' \in \omega^\omega$ going to infinity with $0 < g, g'$, we have $\mathfrak{c}_{\omega,g} = \mathfrak{c}_{\omega,g'}$.

**Definition 10.5.** Let $\mathfrak{c}_{\omega,\mathcal{C}}$ be the minimal size of a family of slaloms in $\mathcal{C}$ covering all functions in $\omega^\omega$, i. e.

$$\mathfrak{c}_{\omega,\mathcal{C}} := \min\{|\mathcal{S}| \mid \mathcal{S} \subseteq \mathcal{C}, \forall\, x \in \omega^\omega \;\exists\, S \in \mathcal{S} \colon x \in^* S\}.$$

**Lemma 10.6.** *Let $g \in \omega^\omega$ be going to infinity with $0 < g$. Then $\mathfrak{c}_{\omega,g} = \mathfrak{c}_{\omega,\mathcal{C}}$.*

*Proof.* Letting $g^+(k) := k^2$ and $g^-(k) := \log k$, define

$$\mathcal{C}_+ := \left\{ S \;\middle|\; S \in^* \prod_{k<\omega} [\omega]^{\leq g^+(k)} \right\}$$

and

$$\mathcal{C}_- := \left\{ S \;\middle|\; S \in^* \prod_{k<\omega} [\omega]^{\leq g^-(k)} \right\}.$$

It is clear that $\mathcal{C}_+ \supseteq \mathcal{C} \supseteq \mathcal{C}_-$, which implies $\mathfrak{c}_{\omega,g^+} \leq \mathfrak{c}_{\omega,\mathcal{C}} \leq \mathfrak{c}_{\omega,g^-}$.

Now Fact 10.4 implies that $\mathfrak{c}_{\omega,g^+} = \mathfrak{c}_{\omega,g^-} = \mathfrak{c}_{\omega,g}$ for any $g \in \omega^\omega$ going to infinity with $0 < g$, and hence $\mathfrak{c}_{\omega,g^+} = \mathfrak{c}_{\omega,\mathcal{C}} = \mathfrak{c}_{\omega,g^-} = \mathfrak{c}_{\omega,g}$. □

**Lemma 10.7.** *For any $g \in \omega^\omega$ going to infinity with $0 < g$, $\mathfrak{c}_{\omega,g} = \operatorname{cof}(\mathcal{N})$.*

*Proof.* By Lemma 10.6, it is enough to show $\mathfrak{c}_{\omega,\mathcal{C}} = \operatorname{cof}(\mathcal{N})$. This is proven in [BJ95, Theorem 2.3.11 (2)]. □

**Theorem 10.8.** *Given $f, g$ as in Definition 10.1, $\mathfrak{c}_{f,g} \leq \operatorname{cof}(\mathcal{N})$.*

*Proof.* This follows immediately from Lemma 10.7 and Fact 10.4. □

**Observation 10.9.** Apart from this inequality, there are no limitations on the placement of the $\mathfrak{c}_{f_\xi,g_\xi}$ relative to the other cardinal characteristics in this paper.

Recall that $\kappa_\xi^{\aleph_0} = \kappa_\xi$ for any $\xi < \omega_1$. As a consequence of the theorem, note that if we were to omit the $\operatorname{cof}(\mathcal{N})$ forcing factors entirely, we would then get the following result for $\operatorname{cof}(\mathcal{N})$ in $V[G]$: Let $\lambda := \sup_{\xi<\omega_1} \mathfrak{c}_{f_\xi,g_\xi}$. Then it is clear that $\lambda \leq \operatorname{cof}(\mathcal{N})$ by Lemma 10.7 and $\operatorname{cof}(\mathcal{N}) \leq \lambda^{\aleph_0}$ by the fact that there are only $\lambda^{\aleph_0}$ many reals after forcing with $\mathbb{Q}_{\text{non-ct}}$ (recall section 8). If $\operatorname{cof} \lambda \geq \omega_1$, then GCH in the ground model implies that $\lambda^{\aleph_0} = \lambda$ and hence $\operatorname{cof}(\mathcal{N}) = \lambda$.

Before we prove the cardinal characteristics' inequalities, we need to show that there indeed is a congenial $\omega_1$-sequence of function pairs as defined in Definition 3.2. We can show even more:

**Lemma 10.10.** *There is a congenial sequence $\langle f_\xi, g_\xi \mid \xi < \mathfrak{c}\rangle$ of continuum many function pairs.*

*Proof.* (This proof is a modification and simplification of the construction in [GS93, Example 3.3].)

Recall that we need to show the following properties from Definition 3.2:

(i) For all $\xi$ and for all $k < \omega$, $n^B_{4k+2} \leq g_\xi(k) < f_\xi(k) \leq n^S_{4k+2}$.

(ii) For all $\xi$, $\lim_{k\to\infty} \frac{\log f_\xi(k)}{n^B_{4k+2} \cdot \log g_\xi(k)} = \infty$.

(iii) For all $\xi \neq \zeta$, either $\lim_{k\to\infty} \frac{f_\zeta(k)^2}{g_\xi(k)} = 0$ or $\lim_{k\to\infty} \frac{f_\xi(k)^2}{g_\zeta(k)} = 0$.

Also recall the definitions of $n^B_{4k+2}$ and $n^S_{4k+2}$ in Definition 4.8.

Let $\langle e_k \mid k < \omega\rangle$ be an increasing sequence such that $e_k > n^B_{4k+2} \geq 2$ for all $k < \omega$; we will canonically choose $e_k := n^B_{4k+2}+1$, but any other sequence with this property would work. Take the complete binary tree $T := 2^{<\omega}$ and enumerate $T \cap 2^k$ in lexicographic order as $\{s_k^1, \dots, s_k^{2^k}\}$. We now define a pair of functions $(f_\xi, g_\xi)$ for each branch $b_\xi \in [T]$ by the following rule: If $b_\xi{\restriction}_k = s_k^i$, then $f_\xi(k) := (n^B_{4k+2})^{e_k^{10i}}$ and $g_\xi(k) := (n^B_{4k+2})^{e_k^{10i-5}}$.

It is clear that by definition, $n^B_{4k+2} \leq g_\xi(k) < f_\xi(k)$, and recalling the fact that in Definition 4.8 we set

$$n^S_{4k+2} := (n^B_{4k+2})^{e_k^{10\cdot 2^k}},$$

we also have $f_\xi(k) \leq n^S_{4k+2}$. This proves property (i).

To show property (ii), we first note that given any $k \in \omega$, for some $1 \leq i \leq 2^k$ (depending on $b_\xi$ and $k$) we have $\frac{\log f_\xi(k)}{\log g_\xi(k)} = \frac{e_k^{10i}}{e_k^{10i-5}} = e_k^5$. As $n^B_{4k+2} < e_k$, we have $\frac{\log f_\xi(k)}{\log g_\xi(k)} \cdot \frac{1}{n^B_{4k+2}} > e_k^4$, which diverges to infinity.

Finally, consider $\xi \neq \zeta$, without loss of generality such that $b_\xi < b_\zeta$ in the natural lexicographic order on the branches of $[T]$. (If we have $b_\zeta < b_\xi$, we can just prove the other statement in property (iii) the same way.) Taking the first $k$ such that $b_\xi(k-1) \neq b_\zeta(k-1)$, there are $1 \leq i < j \leq 2^k$ such that $b_\xi{\restriction}_k = s_k^i$ and $b_\zeta{\restriction}_k = s_k^j$ (and analogously for any larger $k$). Then it follows that

$$F(\xi, \zeta, k) := \frac{f_\xi(k)^2}{g_\zeta(k)} = (n^B_{4k+2})^{2\cdot e_k^{10i} - e_k^{10j-5}}$$

and since $e_k > 2$ and $10i+1 < 10j-5$, we have $2\cdot e_k^{10i} - e_k^{10j-5} < e_k^{10i+1} - e_k^{10j-5} < -1$, hence $F(\xi, \zeta, k) < 1/n^B_{4k+2}$, which goes to 0 as $k$ goes to infinity – as required to show property (iii). $\square$

We point out once more that in order to make it easier to read, the construction above is actually slightly less general than the one in [GS93, Example 3.3]; in our case, the pairs of functions are not only pointwise "far apart", but instead even have the same ordering between them at each point. The reader can easily convince themself that the more general construction would also work in the same way.

**Lemma 10.11.** *$\mathbb{Q}$ forces that for all $\xi < \omega_1$, $\mathfrak{c}_{f_\xi, g_\xi} \geq \kappa_\xi$.*

*Proof.* Fix some $\xi < \omega_1$. Let $G$ be $\mathbb{Q}$-generic and let $\mathcal{S}$ be some family of $g_\xi$-slaloms in $V[G]$ of size less than $\kappa_\xi$. Each $S \in \mathcal{S}$ is read continuously only using indices in some countable subset $B_S$ of $A$ and there are fewer than $\kappa_\xi^{\aleph_0} = \kappa_\xi = |A_\xi|$ many $S$, so letting $B := \bigcup_{S \in \mathcal{S}} B_S$, all of $\mathcal{S}$ is read continuously only using indices in $B$ and there is some $\alpha \in A_\xi \smallsetminus B$.

Now assume towards a contradiction that there were some $g_\xi$-slalom $S^* \in V[G{\restriction}_B]$ covering the generic $\dot{y}_\alpha$. Working in $V$, this means that there is a $\mathbb{Q}{\restriction}_B$-name $\dot{S}^*$ and a condition $p \in \mathbb{Q}$ such that $\Vdash_{\mathbb{Q}{\restriction}_B}$ "$\dot{S}^*$ is a $g_\xi$-slalom" and $p \Vdash_{\mathbb{Q}}$ "$\dot{S}^*$ covers $\dot{y}_\alpha$".

But then we can find some $k < \omega$ such that $|p(\alpha, 4k+2)| > g_\xi(k)$. Find $q \leq p$ by first strengthening $p{\restriction}_B$ to decide $\dot{S}^*_k = T$ and then finding some $x \in p(\alpha, 4k+2) \smallsetminus T$ and replacing $p(\alpha, 4k+2)$ by $\{x\}$. The condition $q$ then forces the desired contradiction, proving that fewer than $\kappa_\xi$ many $g_\xi$-slaloms cannot suffice to cover all functions in $\prod_{k<\omega} f_\xi(k)$. $\square$

To prove the converse, we first have to prepare just a few more technical tools.

**Definition 10.12.** Let $p \in \mathbb{Q}$ and let $\dot{t}$ be a $\mathbb{Q}$-name for a function in $\prod_{k<\omega} n^S_{4k+2}$. We say that $p$ *punctually reads* $\dot{t}$ if for each $k < \omega$, $\dot{t}\restriction_{k+1}$ is decided below $4k + 2$ (including height $4k + 2$).

**Definition 10.13.** Let $\langle y_k \mid k \in \omega\rangle$ be a sequence of numbers such that $n^S_{4k+2} \leq 2^{y_k}$ and let $z(k) := \sum_{\ell\leq k} y_\ell$. After identifying $2^\omega$ with $\prod_{k\in\omega} 2^{y_k}$ and fixing maps from $2^{y_k}$ onto $n^S_{4k+2}$, we get a coherent family of surjective maps $C_k \colon 2^{z(k)} \to \prod_{\ell\leq k} n^S_{4\ell+2}$, which taken together describe a a natural "coding" map $C$ from $2^\omega$ onto $\prod_{k\in\omega} n^S_{4k+2}$ which satisfies $C_k(s\restriction_{z(k)}) = C(s)\restriction_{k+1}$ for all $k$.

**Corollary 10.14.** *Let $\dot{t}$ be a $\mathbb{Q}$-name for a function in $\prod_{k<\omega} n^S_{4k+2}$ and $p \in \mathbb{Q}$. Then there is a $q \leq p$ punctually reading $\dot{t}$.*

*Proof.* Using the map from the preceding definition, we can find a $\mathbb{Q}$-name $\dot{s}$ for an element of $2^\omega$ such that for any $k < \omega$, $C_k(\dot{s}\restriction_{z(k)}) = \dot{t}\restriction_{k+1}$. Find $q \leq p$ rapidly reading $\dot{s}$; by the preceding lemma, $q$ then punctually reads $\dot{t}$. □

**Definition 10.15.** For a modest $p \in \mathbb{Q}$, we call $4k + 2$ a *slalom-splitting level* if there is an $\alpha \in \operatorname{supp}(p)$ such that $|p(\alpha, 4k + 2)| > 1$. We refer to this unique index by $\alpha_k$, and the corresponding type by $\zeta_k < \omega_1$.

**Definition 10.16.** Fix some $\xi < \omega_1$. We call a condition $p \in \mathbb{Q}$ *$\xi$-prepared* if for all $k < \omega$, one of the following three statements holds:

- $4k + 2$ is not a slalom-splitting level of $p\restriction_{\operatorname{supp}(p)\setminus\{\xi\}}$,
- $f_{\zeta_k}(k)^2 < g_\xi(k)$, or
- $f_\xi(k)^2 < g_{\zeta_k}(k)$.

**Lemma 10.17.** *Fix some $\xi < \omega_1$ and let $p \in \mathbb{Q}$. Then there is a $\xi$-prepared $q \leq p$.*

*Proof.* We do the following steps for each $\zeta \neq \xi$.

Note that per property (iii) in Definition 3.2, we know that either

$$\lim_{k\to\infty} \frac{f_\zeta(k)^2}{g_\xi(k)} = 0 \quad \text{or} \quad \lim_{k\to\infty} \frac{f_\xi(k)^2}{g_\zeta(k)} = 0$$

and hence there must be some $k_\zeta$ such that $\frac{f_\zeta(k)^2}{g_\xi(k)} < 1$ for all $k \geq k_\zeta$ or $\frac{f_\xi(k)^2}{g_\zeta(k)} < 1$ for all $k \geq k_\zeta$. Now for $k < k_\zeta$, shrink each creature in $p(\zeta, 4k + 2)$ to an arbitrary singleton to get $q(\zeta, 4k + 2)$. The resulting $q$ is then $\xi$-prepared. □

**Lemma 10.18.** *$\mathbb{Q}$ forces that for all $\xi < \omega_1$, $\mathfrak{c}_{f_\xi,g_\xi} \leq \kappa_\xi$.*

*Proof.* Fix some $\xi < \omega_1$ and let $Z := \bigcup_{\kappa_{\mathrm{t}}\leq\kappa_\xi} A_{\mathrm{t}}$.[33] We will prove that the $g_\xi$-slaloms in $V^{\mathbb{Q}\restriction Z}$ cover $\prod_{k<\omega} f_\xi(k)$; this suffices since by Lemma 6.5 (and by the fact that $\mathbb{Q}\restriction_Z$ is a complete subforcing of $\mathbb{Q}$, see Lemma 4.19), $\Vdash_{\mathbb{Q}\restriction_Z} 2^{\aleph_0} \leq \kappa_\xi$ and hence $\Vdash_{\mathbb{Q}} (2^{\aleph_0})^{V^{\mathbb{Q}\restriction Z}} \leq \kappa_\xi$.

[33] Note, however, that the $\mathrm{t} \notin$ slalom are not especially relevant here. The case distinction below only cares about the $A_\zeta$ with $\zeta < \omega_1$ – and whether $\kappa_\zeta \leq \kappa_\xi$ or $\kappa_\zeta > \kappa_\xi$ –, but the definition is just cleaner in this more general formulation.

So let $\dot{t}$ be a $\mathbb{Q}$-name for a function in $\prod_{k<\omega} f_\xi(k)$ and let $p^* \in \mathbb{Q}$ be an arbitrary condition. Find $p \leq p^*$ such that $p$ punctually reads $\dot{t}$ and is $\xi$-prepared.

We will find a condition $q \leq p$ and a $\mathbb{Q}{\restriction}_Z$-name $\dot{S}$ for a $g_\xi$-slalom such that $q \Vdash$ "$\dot{S}$ covers $\dot{t}$".

To find $q$ and define $\dot{S}$, we go through the levels of the form $4k+2$ and make the following case distinction. (We know that one of the following cases must hold since $p$ is $\xi$-prepared.)

Case 0: $4k+2$ is not a slalom-splitting level of $p$.

In this case, we have that $|\operatorname{poss}(p, {<}4k+3)| = |\operatorname{poss}(p, {<}4k+2)|$ (since at level $4k+2$, there is only one possible extension for each possibility from below). Hence letting $q(\text{slalom}, 4k+2) := p(\text{slalom}, 4k+2)$ and defining

$$\dot{S}_k := \{x < \omega \mid \exists\, \eta \in \operatorname{poss}(q, {<}4k+3)\colon p \wedge \eta \Vdash \dot{t}(k) = x\},$$

we actually have a (ground model) set of size at most

$$|\operatorname{poss}(p, {<}4k+3)| = |\operatorname{poss}(p, {<}4k+2)| < n^P_{<4k+2} < n^B_{4k+2} \leq g_\xi(k),$$

and clearly $q \Vdash \dot{t}(k) \in \dot{S}_k$.

Case 1: $4k+2$ is a slalom-splitting level of $p$, but $\zeta_k$ is such that $\kappa_{\zeta_k} \leq \kappa_\xi$.

In this case, $\alpha_k \in Z$. We once more let $q(\text{slalom}, 4k+2) := p(\text{slalom}, 4k+2)$ and define $\dot{S}_k$ to be a $\mathbb{Q}{\restriction}_Z$-name satisfying

$$\Vdash_{\mathbb{Q}\restriction_Z} \dot{S}_k = \{x < \omega \mid \exists\, \eta \in \operatorname{poss}(q, {<}4k+3)\colon p \wedge \eta \Vdash \text{ ``}\dot{t}(k) = x \text{ and } \eta(\alpha_k) \subseteq \dot{y}_{\alpha_k}\text{''}\},$$

which means we only allow those possibilities $\eta$ which are compatible with the generic real $\dot{y}_{\alpha_k}$ added by the forcing factor $\mathbb{Q}_{\alpha_k}$. Similar to the previous case, this means that $\Vdash_{\mathbb{Q}\restriction_Z} |\dot{S}_k| \leq |\operatorname{poss}(p, {<}4k+2)|$, for the following reason: Let $\varepsilon := \dot{y}_{\alpha_k}$; then in the definition of $\dot{S}_k$ above, the only admissible possibilities $\eta \in \operatorname{poss}(p, {<}4k+3)$ are those of the form $\eta = \nu^\frown\varepsilon$ for some $\nu \in \operatorname{poss}(p, {<}4k+2)$. Hence $\Vdash_{\mathbb{Q}\restriction_Z} |\dot{S}_k| < g_\xi(k)$, and by definition $q \Vdash \dot{t}(k) \in \dot{S}_k$.

Case 2: $4k+2$ is a slalom-splitting level of $p$, $\kappa_{\zeta_k} > \kappa_\xi$ and $f_{\zeta_k}(k)^2 < g_\xi(k)$.

From $f_{\zeta_k}(k)^2 < g_\xi(k)$, we get the following:

$$\begin{aligned} |\operatorname{poss}(p, {<}4k+3)| &\leq |\operatorname{poss}(p, {<}4k+2)| \cdot f_{\zeta_k}(k) \leq n^P_{<4k+2} \cdot f_{\zeta_k}(k) \\ &< n^B_{4k+2} \cdot f_{\zeta_k}(k) < f_{\zeta_k}(k)^2 < g_\xi(k) \end{aligned}$$

By once more letting $q(\text{slalom}, 4k+2) := p(\text{slalom}, 4k+2)$ and defining

$$\dot{S}_k := \{x < \omega \mid \exists\, \eta \in \operatorname{poss}(q, {<}4k+3)\colon p \wedge \eta \Vdash \dot{t}(k) = x\},$$

we hence again have a (ground model) set of size at most $g_\xi(k)$, and by definition $q \Vdash \dot{t}(k) \in \dot{S}_k$.

Case 3: $4k+2$ is a slalom-splitting level of $p$, $\kappa_{\zeta_k} > \kappa_\xi$ and $f_\xi(k)^2 < g_\zeta(k)$.

This is the only case where we have to do any actual work to get $q$, as we cannot simply collect all potential values of $\dot{t}(k)$. Instead, we will first have to use the bigness properties of the norm to reduce the number of potential values. To begin, we remark that letting $m_k := |\operatorname{poss}(p, {<}4k+2)|$, $c := f_\xi(k)$ and $d := \lfloor {}^{g_\xi(k)}/_{m_k} \rfloor$, we have

$$\frac{c}{d} = \frac{f_\xi(k)}{\lfloor \frac{g_\xi(k)}{m_k} \rfloor} \leq f_\xi(k) \cdot \frac{2 \cdot |\operatorname{poss}(p, {<}4k+2)|}{g_\xi(k)} \leq f_\xi(k) \cdot \frac{2 \cdot n^P_{<4k+2}}{g_\xi(k)} < f_\xi(k)^2 < g_\zeta(k)$$

and since $\|\cdot\|_{\zeta,4k+2}$ has $n^B_{4k+2}$-strong $g_\zeta(k)$-bigness by Theorem 5.6, according to Observation 5.5 it also has $n^B_{4k+2}$-strong $(c,d)$-bigness.

Enumerate $\operatorname{poss}(p, < 4k+2) =: \{\eta_1, \dots, \eta_{m_k}\}$. We claim we can find a sequence of subsets $p(\alpha_k, 4k+2) = F^0_k \supseteq F^1_k \supseteq F^2_k \supseteq \dots \supseteq F^{m_k}_k$ and a sequence of sets $C_j$ with the following properties (for each $1 \leq j \leq m_k$):

(i) $\|F^j_k\|_{\zeta,4k+2} \geq \|F^{j-1}_k\|_{\zeta,4k+2} - {}^1/_{n^B_{4k+2}}$
(ii) $|C_j| \leq d$
(iii) $p \wedge (\eta^j{}^\frown x) \Vdash \dot{t}(k) \in C_j$ holds for all $x \in F^j_k$.

We know that, given $F^j_k$, for each $x \in F^j_k$, we have that $p \wedge (\eta^j{}^\frown x)$ decides $\dot{t}(k)$ by punctual reading of $\dot{t}$ (noting that $\eta^j{}^\frown x \in \operatorname{poss}(p, {<}4k+3)$). Since there are most $c$ many possible values for $\dot{t}(k)$, we can use $n^B_{4k+2}$-strong $(c,d)$-bigness of the norm $\|\cdot\|_{\zeta,4k+2}$ to find $F^{j+1}_k \subseteq F^j_k$ and $C_{j+1}$ with the desired properties, proving our claim.

Now, we define $F_k := F^{m_k}_k$. Since $m_k \leq n^P_{<4k+2} < n^B_{4k}$, by (i) we have that

$$\|F_k\|_{\zeta,4k+2} \geq \|p(\alpha_k, 4k+2)\|_{\zeta,4k+2} - \frac{m_k}{n^B_{4k+2}} \geq \|p(\alpha_k, 4k+2)\|_{\zeta,4k+2} - 1.$$

Hence, defining $q(\alpha_k, 4k+2) := F_k$ and $q(\text{slalom}, 4k+2) := p(\text{slalom}, 4k+2)$ elsewhere (i. e. on $A_{\text{slalom}} \setminus \{\alpha_k\}$) does not negatively affect the lim sup properties of the norm.

Finally, let $\dot{S}_k := \bigcup_{1 \leq j \leq m_k} C_j$. Then by (ii), we have that $|\dot{S}_k| \leq d \cdot m_k = g_\xi(k)$, and (iii) implies $q \Vdash \dot{t}(k) \in \dot{S}_k$, finishing the proof. □

This proves (M5) of Theorem 1.1.

## 11. non($\mathcal{N}$) $\geq \kappa_{\text{nn}}$ AND cof($\mathcal{N}$) $\geq \kappa_{\text{cn}}$

The proofs in this section are more or less identical to those in [FGKS17], although we sincerely hope we have improved the presentation.

To prove $\operatorname{non}(\mathcal{N}) \geq \kappa_{\text{nn}}$, we define some null sets in the extension, similar to the definition of $\dot{M}_\alpha$ in section 9. Recall that for $\alpha \in A_{\text{nn}}$, the generic object $\dot{y}_\alpha$ is a $\mathbf{heights}_{*\text{n}}$-sequence of subsets $\dot{R}_{\alpha,L}$ of $2^{I_L}$ of relative size $1 - 2^{-n^B_L}$.[34] Since the sequence of $n^B_L$, $L \in \mathbf{heights}_{*\text{n}}$, is strictly monotone, we have

$$\prod_{L \in \mathbf{heights}_{*\text{n}}} \left(1 - \frac{1}{2^{n^B_L}}\right) > 0$$

[34] Recall that $\mathbf{heights}_{*\text{n}} = \{4k+1 \mid k < \omega\}$.

and hence the set $\{r \in 2^\omega \mid \forall\, k < \omega \colon r{\restriction}_{I_{4k+1}} \in \dot{R}_{\alpha,4k+1}\}$ is positive. It follows that the set $\{r \in 2^\omega \mid \forall^\infty\, k < \omega \colon r{\restriction}_{I_{4k+1}} \in \dot{R}_{\alpha,4k+1}\}$ has measure one, and therefore its complement

$$\dot{N}_\alpha := \{r \in 2^\omega \mid \exists^\infty\, k < \omega \colon r{\restriction}_{I_{4k+1}} \notin \dot{R}_{\alpha,4k+1}\}$$

is a name for a null set.

Recall that (by Theorem 5.6) for each $L \in \mathbf{heights}_{*\mathrm{n}}$, $(\mathbf{POSS}_{\mathrm{nn},L}, \|\cdot\|_{\mathrm{nn},L})$ has strong $n^B_L$-bigness. We show a similar, more specific property:

**Lemma 11.1.** *Let $L \in \mathbf{heights}_{*\mathrm{n}}$, $X \subseteq \mathbf{POSS}_{\mathrm{nn},L}$ and $E \subseteq 2^{I_L}$. Let $X' := \{H \in X \mid H \cap E = \varnothing\}$. Then $\|X'\|^{\mathrm{intersect}}_L \geq \|X\|^{\mathrm{intersect}}_L - |E|$.*

*If additionally $|E| \leq n^P_{<L}$, then it follows that $\|X'\|_{\mathrm{nn},L} \geq \|X\|_{\mathrm{nn},L} - 1/\log n^B_L$.*

*Proof.* For the first part, assume that some $Y$ witnesses $\|X'\|^{\mathrm{intersect}}_L$; then $Y \cup E$ certainly witnesses $\|X\|^{\mathrm{intersect}}_L$.

For the second part, note that $n^P_{<L} \leq n^B_L/2$ and hence $|E| \leq (n^B_L)^{\|X\|_{\mathrm{nn},L}}/2$. Since $n^B_L \geq 2$ and assuming $X$ is non-trivial (without loss of generality, assume $\|X\|_{\mathrm{nn},L} \geq 2$), a trivial inequality gives

$$\frac{(n^B_L)^{\|X\|_{\mathrm{nn},L}}}{2} = \frac{(\|X\|^{\mathrm{intersect}}_L)^{1/n^B_L}}{2} \leq \left(1 - \frac{1}{2^{n^B_L}}\right) \cdot \|X\|^{\mathrm{intersect}}_L$$

and hence the first part implies

$$\begin{aligned} \|X'\|_{\mathrm{nn},L} &= \frac{\log \|X'\|^{\mathrm{intersect}}_L}{n^B_L \log n^B_L} \geq \frac{\log(\|X\|^{\mathrm{intersect}}_L - |E|)}{n^B_L \log n^B_L} \geq \frac{\log(\|X\|^{\mathrm{intersect}}_L / 2^{n^B_L})}{n^B_L \log n^B_L} \\ &= \frac{\log \|X\|^{\mathrm{intersect}}_L}{n^B_L \log n^B_L} - \frac{\log 2^{n^B_L}}{n^B_L \log n^B_L} = \|X\|_{\mathrm{nn},L} - \frac{1}{\log n^B_L}. \end{aligned}$$

□

**Lemma 11.2.** *Let $\dot{r} \in 2^\omega$ be a name for a real and let $p$ rapidly read $\dot{r}$ not using the index $\alpha \in A_{\mathrm{nn}}$. Then $p \Vdash \dot{r} \in \dot{N}_\alpha$.*

*Proof.* We first remark that as in Lemma 9.3, it suffices to prove that there is an $s \leq p$ such that $s \Vdash \dot{r} \in \dot{N}_\alpha$. Similar to that proof, we will only have to modify $p$ in very few places to get the desired condition $s$. Without loss of generality, assume that $\alpha \in \mathrm{supp}(p)$.

We will only modify $p$ at index $\alpha$ for infinitely many heights in $\mathbf{heights}_{*\mathrm{n}}$. Assume we have already modified $n$ many heights $L_0, \ldots, L_{n-1}$; pick some $L_n := 4k_n + 1 \in \mathbf{heights}_{*\mathrm{n}}$ such that $p(\alpha, L_n)$ is non-trivial and has a norm of at least $n$. By rapid reading, we know that $\dot{r}{\restriction}_{I_{L_n}}$ is decided $\leq L_n$ by $p$, and as in Lemma 9.3, by modesty and since $\dot{r}$ does not depend on the index $\alpha$, it is even decided below $L_n$. Hence the set $E_n$ of possible values for $\dot{r}{\restriction}_{I_{L_n}}$ has size at most $n^P_{<L}$; by the preceding lemma, replacing $p(\alpha, L_n)$ by $C_n := \{H \in p(\alpha, L_n) \mid H \cap E_n = \varnothing\}$ only decreases the norm by at most 1.

The condition $s$ resulting from replacing each $p(\alpha, L_n)$ by $C_n$ then fulfils $s \Vdash \dot{r} \in \dot{N}_\alpha$ by definition, as $s \Vdash \dot{r}{\restriction}_{I_{4k_n+1}} \notin \dot{R}_{\alpha,4k_n+1}$ holds for all $n < \omega$. □

**Corollary 11.3.** *$\mathbb{Q}$ forces* $\mathrm{non}(\mathcal{N}) \geq \kappa_{\mathrm{nn}}$.

*Proof.* The proof is identical to the proof of Corollary 9.4. □

To prove $\mathrm{cof}(\mathcal{N}) \geq \kappa_{\mathrm{cn}}$, we define null sets in the extension in the same way we did at the start of this section, namely

$$\dot{N}_\alpha := \{r \in 2^\omega \mid \exists^\infty k < \omega \colon r{\restriction}_{I_{4k+1}} \notin \dot{R}_{\alpha,4k+1}\}$$

for $\alpha \in A_{\mathrm{cn}}$. However, the purpose of these null sets will be quite different; rather than covering all reals in the extension which do not depend on the index $\alpha$, they will avoid being covered by any null set not depending on $\alpha$.

We wish to spare the reader the details of the combinatorial arguments from [FGKS17, section 9], and hence will only sketch the modifications necessary to see why the proofs in [FGKS17] still hold. The relevant result we will be using is [FGKS17, Lemma 10.2.1], in the following form:

**Lemma 11.4.** *Fix a height $L \in$ **heights**$_{*\mathrm{n}}$, an index $\alpha \in A_{\mathrm{cn}}$ and a creature $C \subseteq \mathbf{POSS}_{\alpha,L}$ such that $\|C\|_{\mathrm{cn},L} \geq 2$.*

*(i) Given $T \subseteq 2^{I_L}$ of relative size at least $1/2$, we can strengthen $C$ to a creature $D$ such that $T \not\subseteq X$ for all $X \in D$ and such that*

$$\|D\|_{\mathrm{cn},L} \geq \|C\|_{\mathrm{cn},L} - \frac{1}{2^{\min I_L} \cdot n^B_L}.$$

*(ii) Given a probability space $\Omega$ and a function $F \colon C \to \mathcal{P}(\Omega)$ mapping each $X \in C$ to some $F(X) \subseteq \Omega$ of measure at least $1/n^B_L$, we can strengthen $C$ to a creature $D$ such that $\bigcap_{X \in D} F(X)$ has measure at least $1/n^B_{L^+}$ and such that*

$$\|D\|_{\mathrm{cn},L} \geq \|C\|_{\mathrm{cn},L} - \frac{1}{2^{\min I_L} \cdot n^B_L}.$$

*Proof sketch.* In [FGKS17], Lemma 10.2.1 is an immediate consequence of two properties – (i) follows from Lemma 9.2.2 and (ii) follows from Eq. (9.1.4).

The second part is straightforward: By the considerations in Observation 3.6, it is clear that (ii) still follows for our new norm as long as $n^B_{L^+} > n^B_L \cdot 2^{n^S_L+1}$, which is the case by the definition of $n^B_{L^+}$.

The first part (Lemma 9.2.2) requires a bit more thought (since our modification to take the logarithm of the $\mathrm{nor}^{\div}_{I,b}$ complicates the direct argument). Define

$$\Delta_L := \binom{2^{|I_L|-1}}{2^{n^B_L-1}}.$$

The relevant statement in Lemma 9.2.2 then is: Given $C$ and $T$ as in (i), we can find $D$ such that $T \not\subseteq X$ for all $X \in D$ and such that

$$(*_8) \qquad |C \setminus D| \leq \Delta_L.$$

We first explain why Eq. ($*_8$) implies (i): Let $\delta := |C \setminus D|$. Then $\|C\|_{\mathrm{cn},L} \geq 2$ implies $|C| \geq 2\Delta_L$, and $\Delta_L \geq \delta$ implies $|D| \geq |C|/2 \geq \delta$, so $\delta/|D| \leq 1$. By the well-known fact that

$$\frac{\ln(z+\varepsilon) - \ln z}{\varepsilon} \leq \frac{1}{z}$$

and the fact that $|C| = |D| + \delta$, we get

$$\log|C| - \log|D| \leq \frac{\delta}{|D|} \cdot \log e \leq 2.$$

By all this, we know that having found such a $D$, the numerators of the fractions in $\|C\|_{\mathrm{cn},L}$ and $\|D\|_{\mathrm{cn},L}$ differ by at most 2, and hence the whole norms differ by at most $1/(2^{\min I_L} \cdot n^B_L)$ (actually, even less).

Finally, the reason why Eq. ($*_8$) holds is the same combinatorial consideration explained in [FGKS17, Lemma 9.2.2 (1)]. □

**Fact 11.5.** We require a few facts about the correspondence between trees of measure $1/2$ and null sets.

(i) Let $T \subseteq 2^{<\omega}$ be a leafless tree of measure $1/2$ (and recall that such trees bijectively correspond to closed sets of measure $1/2$). For $X \subseteq 2^\omega$, let $X + 2^{<\omega} := \bigcup\{X + r \mid r \in 2^{<\omega}\}$, the set of all rational translates of $X$ (where $X + r := X + r^\frown\langle 000\dots\rangle$). Then the set $N_T := 2^\omega \setminus ([T] + 2^{<\omega})$ is a null set closed under rational translations.

(ii) Conversely, given an arbitrary null set $N$, there is a leafless tree $T$ of measure $1/2$ such that $N \subseteq N_T$, since the complement of $N + 2^{<\omega}$ must contain a closed set of size $1/2$.

(iii) Let $k < \omega$ and $s \in T \cap 2^k$. We define the relative measure of $s$ in $T$ as $2^k \cdot \lambda([T] \cap [s])$. Analogously, for finite trees $T \subseteq 2^{\leq m}$ (such that there are no leaves below tree level $m$) we define the relative measure of $s \in T \cap 2^k$ for $k \leq m$ in the same way. (For $s \notin T$, the relative measure of $s$ in $T$ is 0, naturally.)

(iv) Given a leafless tree $T \subseteq 2^{<\omega}$, some $s \in T$ of positive relative measure and some $0 < \varepsilon < 1$, there is some extension $t$ of $s$ such that $t$ has relative measure $> \varepsilon$. Moreover, it follows that for all tree levels above the tree level of $t$ there is some extension $u$ of $t$ such that $u$ has relative measure $>\varepsilon$. (These statements are a simple consequence of Lebesgue's density theorem.)

Since the measure of a tree $T$ does not change if we remove any $s \in T$ of relative measure 0, we will be working with such trees instead:

**Definition 11.6.** We call $T \subseteq 2^{<\omega}$ a *sturdy tree* if it has measure $1/2$ and no $s \in T$ has relative measure 0 (in particular, this means $T$ is leafless).

The considerations from Fact 11.5 also hold for sturdy trees, so we will be working with those instead.

Finally, we remark that $2^{2^k}$ is an upper bound for the cardinality of the set $2^{\leq k}$. We can thus code any name for a sturdy tree $\dot{T}$ by a real $\dot{t} \in 2^\omega$ such that $\dot{T} \cap 2^k$ is determined by $\dot{t}{\restriction}_{2^{2^{k+1}}}$, and by the definition of $n^R_{<L}$, if a condition $p$ rapidly reads $\dot{t}$, then for each $\eta \in \mathrm{poss}(p, {\leq}L)$, $p \wedge \eta$ decides $\dot{T} \cap 2^{\max I_L}$; we abbreviate this fact by "$p$ rapidly reads $\dot{T}$".

**Lemma 11.7.** *Let $\dot{T}$ be a name for a sturdy tree and let $p$ rapidly read $\dot{T}$ not using the index $\alpha \in A_{\mathrm{cn}}$. Then $p \Vdash \dot{N}_\alpha \not\subseteq N_{\dot{T}}$, that is, $p$ forces that there is an $s \in \dot{N}_\alpha \cap [\dot{T}]$.*

*Proof.* Once again, it suffices to find a $q \leq p$ and a $\mathbb{Q}$-name $\dot{s}$ for a real such that $q \Vdash \dot{s} \in \dot{N}_\alpha \cap [\dot{T}]$. Without loss of generality, assume that $\alpha \in \mathrm{supp}(p)$. To achieve this, we will modify $p(\alpha)$ at infinitely many $*$n heights to get $q$ and thereafter define the required real $s$ inductively in the extension.

Let $L \in \mathbf{heights}_{*\mathrm{n}}$ be a height, above all the previously modified heights, such that $\|p(\alpha, L)\| \geq 3$. (The condition on the norm is necessary for us to be able to apply Lemma 11.4 (i) sufficiently often.) Let $\dot{T}^* := \dot{T} \cap 2^{\max I_L}$. By rapid reading, $p$ decides $\dot{T}^*$ below $L$ (since $\dot{T}$ does not depend on the index $\alpha$ and by modesty, there is no other index $\beta$ such that $p(\beta, L)$ is non-trivial). In particular, this means that the set $W$ of possible values of $\dot{T}^*$ has size at most $n^P_{<L}$.

We now enumerate all $U \in W$ and all $u \in U \cap 2^{\min I_L}$ with relative measure at least $1/2$ (measured in $U$). Clearly, there are at most $M := n^P_{<L} \cdot 2^{\min I_L}$ many such pairs $(U, u)$. Starting with $C_0 := p(\alpha, L)$, we will iteratively apply Lemma 11.4 (i) to the creature $C^n$ and the tree $u^\frown U\restriction_{2^{I_L}}$ to get a creature $C_{n+1} \subseteq C_n$ which then fulfils the following statement: For each $X \in C_{n+1}$, there is some $u' \in 2^{I_L} \setminus X$ such that $u^\frown u' \in U$, and

$$\|C_{n+1}\|_{\mathrm{cn},L} \geq \|C_n\|_{\mathrm{cn},L} - \frac{1}{2^{\min I_L} \cdot n^B_L}.$$

After going through all $M$ many possible choices of $(U, u)$, we arrive at $D := C_M$, which fulfils the following statement: For each $X \in D$ and each $(U, u)$ as above, there is some $u' \in 2^{I_L} \setminus X$ such that $u^\frown u' \in U$, and $\|D\|_{\mathrm{cn},L} \geq \|p(\alpha, L)\|_{\mathrm{cn},L} - 1$, since $n^P_{<L} < n^B_L$ and hence

$$\frac{M}{2^{\min I_L} \cdot n^B_L} = \frac{n^P_{<L} \cdot 2^{\min I_L}}{n^B_L \cdot 2^{\min I_L}} \leq 1.$$

Denote the condition which emerges after repeating the process above for infinitely many heights by $q$ (and note that $q \leq p$ and $q$ only differs from $p$ at index $\alpha$). We will now work in the forcing extension $V[G]$ (for some generic filter $G$ containing $q$) and construct some $s \in \dot{N}_\alpha \cap [\dot{T}]$. Recall that the requirements on $s$ are that it is a branch of $[\dot{T}]$ and that for infinitely many $L \in \mathbf{heights}_{*\mathrm{n}}$ we have $s\restriction_{I_L} \notin \dot{R}_{\alpha,L}$.

Start with $s_0 := \varnothing$ and $k_0 := 0$. Assume we have already defined $k_n$ and $s_n$ (which will be equal to $s\restriction_{k_n}$) such that $s_n \in \dot{T}$. Since $\dot{T}$ is a sturdy tree and hence has no nodes of relative measure 0, by Fact 11.5 (iv) there is some $k' > k_n$ and a $t \in \dot{T} \cap 2^{k'}$ such that $t$ extends $s_n$ and has relative measure at least $1/2$. Pick a height $L \in \mathbf{heights}_{*\mathrm{n}}$ such that $L$ was considered in the construction of $q$ and such that $\min I_L =: k'' > k'$. Also by Fact 11.5 (iv), there is (still) a $u \in \dot{T} \cap 2^{k''}$ such that $u$ extends $s_n$ and has relative measure at least $1/2$. Let $U := \dot{T} \cap 2^{\max I_L}$ and note that in the construction of $q$, we dealt with the pair $(U, u)$. Hence for all $X \in q(\alpha, L)$ (in particular, the $\dot{R}_{\alpha,L}$ chosen by the generic filter $G$), there is some $u' \in 2^{I_L} \setminus X$ such that $u^\frown u' \in U$. So we can set $s_{n+1} := u^\frown u'$ and $k_{n+1} := \max I_L$ and continue the induction; the resulting $s := \bigcup_{n<\omega} s_n$ is as required. □

**Corollary 11.8.** $\mathbb{Q}$ *forces* $\mathrm{cof}(\mathcal{N}) \geq \kappa_{\mathrm{cn}}$.

*Proof.* Fix a condition $p$, some $\kappa < \kappa_{\mathrm{cn}}$ and a sequence of names of null sets $\langle \dot{N}_i \mid i \in \kappa \rangle$ which $p$ forces to be a basis of null sets. As described above, for each $i \in \kappa$, we can assume that $\dot{N}_i = N_{\dot{T}_i}$ for some name for a sturdy tree $\dot{T}_i$. The rest of the proof is identical to the proof of Corollary 9.4. □

This proves (M4) of Theorem 1.1.

## 12. $\mathrm{non}(\mathcal{N}) \leq \kappa_{\mathrm{nn}}$

For the final proofs, we will require two more lemmata. First, we show that the slalom part of the forcing construction has a property similar to Lemma 11.4 (ii).

**Lemma 12.1.** *Fix a height $L \in \mathbf{heights}_{\mathrm{slalom}}$, a slalom type $\xi \in \mathbf{types}_{\mathrm{slalom}}$, an index $\alpha \in A_\xi$, and a creature $C \subseteq \mathbf{POSS}_{\xi,L}$ such that $\|C\|_{\xi,L} \geq 2$.*

*Given a probability space $\Omega$ and a function $F\colon C \to \mathcal{P}(\Omega)$ mapping each $X \in C$ to some $F(X) \subseteq \Omega$ of measure at least $1/n^B_L$, we can strengthen $C$ to a creature $D$ such that $\bigcap_{X \in D} F(X)$ has measure at least $1/n^B_{L^+}$ and such that*

$$\|D\|_{\xi,L} \geq \|C\|_{\xi,L} - \frac{1}{n^B_L}.$$

*Proof.* As in the proof of Lemma 11.4 (ii), we only require that a statement analogous to [FGKS17, Eq. (9.1.4)] holds (as $n^B_{L^+} > n^B_L \cdot 2^{n^S_L+1}$ is true for any $L$).

We already know that Eq. (9.1.4) holds for a norm with the basic structure[35] of $\frac{\log x}{\log 3b}$; the slalom norms have the basic structure $\frac{\log x}{\log g_\xi(k)}$ (with $L = 4k+2$), and in our construction in Lemma 10.10, each $g_\xi(k)$ is defined as some $e$-th power of $n^B_L$; each such exponent $e$ is assured to be at least 8 and even the smallest $n^B_{(0,0)} \geq 8$, hence $(n^B_L)^e \geq 3n^B_L$ and the same basic property holds for this norm structure, as well.[36] □

The other lemma is one more combinatorial statement about trees.

**Lemma 12.2.** *Given a tree $T \subseteq 2^{<\omega}$ of positive measure and an $\varepsilon > 0$, we call $s \in T \cap 2^k$* fat *if*

$$\lambda([T] \cap [s]) \geq \frac{1-\varepsilon}{2^k}.$$

*Then there is a $k^* < \omega$ such that for all $k \geq k^*$, there are at least $|T \cap 2^k| \cdot (1-\varepsilon)$ many fat nodes $s \in T \cap 2^k$.*

*Proof.* (This is the same proof as the one of [FGKS17, Lemma 10.5.3].)

Let $\mu := \lambda([T])$. Since $|T \cap 2^k| \cdot 2^{-k}$ decreasingly converges to $\mu$, there is some $k^*$ such that for all $k \geq k^*$, we have

$$(*_9) \qquad \frac{|T \cap 2^k|}{2^k} - \mu \cdot \varepsilon^2 \leq \mu.$$

Fix some $k \geq k^*$ and let $f$ be the number of fat $s \in T \cap 2^k$ (and $\ell := |T \cap 2^k| - f$ the number of non-fat $s$).

[35] "Basic structure" in the sense of "ignoring additional factors in the denominator or additive terms in the numerator, the norm is fundamentally of logarithmic character".

[36] Alternatively, we could simply amend Definition 3.2 (i) to demand $3n^B_{4k+2} \leq g_\xi(k)$, instead.

Note that

$$(*_{10}) \qquad \mu \leq f \cdot \frac{1}{2^k} + \ell \cdot \frac{1-\varepsilon}{2^k} = \frac{|T \cap 2^k|}{2^k} - \frac{\ell \cdot \varepsilon}{2^k}$$

and hence Eq. ($*_9$) and Eq. ($*_{10}$) together imply $\ell \leq \mu \cdot 2^k \cdot \varepsilon \leq |T \cap 2^k| \cdot \varepsilon$. Since $f + \ell = |T \cap 2^k|$, it follows that $f = |T \cap 2^k| - \ell \geq |T \cap 2^k|(1 - \varepsilon)$. □

Now recall that in section 8, we proved that $\mathbb{Q}$ had the Sacks property over the complete subforcing poset $\mathbb{Q}_{\text{non-ct}}$ (which consists of all conditions $p$ with $\operatorname{supp}(p) \cap A_{\text{ct}} = \varnothing$). In particular, this implied that any null set in the $\mathbb{Q}$-extension is already contained in some null set in the $\mathbb{Q}_{\text{non-ct}}$-extension.

We will now show that the set $R$ of all reals read rapidly only using indices in $A_{\text{nm}} \cup A_{\text{nn}}$ is not null;[37] by the consideration above, we can work entirely with $\mathbb{Q}_{\text{non-ct}}$ and show that it is not null there. As in the preceding section, we will work with sturdy trees instead of null sets.

**Lemma 12.3.** *Let $\dot{T}$ be a name for a sturdy tree and let $p \in \mathbb{Q}_{\text{non-ct}}$ continuously read $\dot{T}$. Then there is a $q \leq p$ in $\mathbb{Q}_{\text{non-ct}}$ and a name $\dot{r}$ for a real such that $q$ continuously reads $\dot{r}$ only using indices in $A_{\text{nm}} \cup A_{\text{nn}}$ (i. e.* not *using any indices in $A_{\text{cn}} \cup A_{\text{slalom}}$) and such that $q \Vdash \dot{r} \in [\dot{T}]$.*

*Proof.* We will construct $q$ and $\dot{r}$ by induction on $n < \omega$. For each $n$, we will define or show the following:

(i) We will define some $L_n := (4k_n, 0) \in \mathbf{heights}_{\text{nm}}$.[38]

(ii) We will define conditions $q_n \leq p$ such that
- $\|q_n(4k)\|_{\text{nm},4k} \geq n + 3$ for all $k \geq k_n$,
- $q_{n+1} \leq q_n$,
- $q_{n+1}$ and $q_n$ are identical on $\operatorname{supp}(q_n)$ below $L_n$ and any new $\alpha \in \operatorname{supp}(q_{n+1}) \setminus \operatorname{supp}(q_n)$ only enters the support of $q_{n+1}$ above $L_n$,
- $\|q_{n+1}(4k)\|_{\text{nm},4k} \geq n$ for all $k_n \leq k < k_{n+1}$, and
- for each $\alpha \in \operatorname{supp}(q_{n+1}, L_n) \setminus A_{\text{nm}}$ of type t, there is a height $L$ with $L_n < L < L_{n+1}$ such that $\|q_{n+1}(\alpha, L)\|_{\text{t},L} \geq n$.

Thus $\langle q_n \mid n < \omega \rangle$ will be a descending sequence of conditions converging to a condition $q$.

(iii) We will define some $i_n < \omega$ and a name $\dot{r}_n$ for an element of $\dot{T} \cap 2^{i_n}$ such that $q_n$ decides $\dot{r}_n$ below $L_n$ only using indices in $A_{\text{nm}} \cup A_{\text{nn}}$.

(iv) We will require that $i_n$ is not "too large" with respect to $L_n$ in the sense that $2^{i_n+2} < n^B_{L_n}$. (Since $n^B_L$ grows quickly and monotonously, it will suffice to show $2^{i_n+2} < 4k_n$.)[39]

(v) The $i_n$ will be such that $i_{n+1} > i_n$.

(vi) The $\dot{r}_n$ will be such that $\dot{r}_{n+1}$ is forced (by $q_{n+1}$) to extend $\dot{r}_n$.

Thus $q$ will force that $\dot{r} := \bigcup_{n<\omega} \dot{r}_n$ will the the desired branch in $[\dot{T}]$.

[37] Recall that we had demanded $\kappa_{\text{nm}} \leq \kappa_{\text{nn}}$.

[38] Since we are working in $\mathbb{Q}_{\text{non-ct}}$, in contrast to section 7 we do not have to worry about the heights being $p$-agreeable.

[39] The reason for the term "+2" will become apparent in Step 3 of the construction below.

(vii) Finally, we will also construct a name $\dot{T}_n$ which $q_n$ will force to be
- a subtree of $\dot{T}$ with stem $\dot{r}_n$ and relative measure greater than $1/2$ (i. e. $\lambda([\dot{T}_n]) > 1/2 \cdot 2^{-i_n}$)
- which is read continuously by $q_n$ in such a way that below $L_n$, the reading only uses indices in $A_{\mathrm{nm}} \cup A_{\mathrm{nn}}$.

Step 0: To start the induction, assume that $p$ already decides the stem of $\dot{T}$ to be $r_0$ and define $i_0 := 0$, $\dot{r}_0 := r_0$ and $\dot{T}_0 := \dot{T}$. Choose $L_0 = (4k_0, 0)$ such that $\|p(4k')\|_{\mathrm{nm},4k'} \geq 3$ for all $k' \geq k_0 > 1$ ("$> 1$" to ensure property (iv)) and let $q_0$ be the condition resulting from extending the trunk of $p$ to $L_0$. It is clear that properties (i)–(vi) are fulfilled by definition, and property (vii) holds since below $L_0$, there is only a single possibility, and hence the reading of $\dot{T}_0$ cannot depend on any indices in $A_{\mathrm{cn}} \cup A_{\mathrm{slalom}}$ below $L_0$.

In the following steps, assume we have constructed the required objects ($L_n = (4k_n, 0)$, $q_n$, $i_n$, $\dot{r}_n$ and $\dot{T}_n$) for some $n < \omega$; we will now proceed to construct them for $n+1$.

Step 1: Choose a height $L^* = (4k^*, 0)$ large enough such that for each $\alpha \in \mathrm{supp}(q_n, L_n) \setminus A_{\mathrm{nm}}$ of type t, there is a height $L$ with $L_n < L < L^*$ such that $\|q_n(\alpha, L)\|_{\mathrm{t},L} \geq n+1$.

It is forced (by $q_n$) that Lemma 12.2 holds for $\dot{T}_n$ and $\varepsilon := 1/(n^P_{<L_n} \cdot n^P_{<L^*})$. Hence there is a name for a tree level $\dot{m}$ such that from $\dot{m}$ upwards, there are many fat nodes in $\dot{T}_n$. We can use Lemma 7.4 to strengthen $q_n$ to $q'$ such that

- $q_n$ and $q'$ are identical below $L^*$,
- the nm norms of $q'$ remain at least $n+2$ starting from $4k^*$, and
- there is an $m^* > i_n$ such that $q' \Vdash m^* \geq \dot{m}$.

Hence Lemma 12.2 is forced to hold for this $m^*$ as well, and there is a name $\dot{F} \subseteq \dot{T}_n \cap 2^{m^*}$ for a "large" set of fat nodes. This $m^*$ will be our $i_{n+1}$.

Step 2: We apply Lemma 7.4 a second time to strengthen $q'$ to $q''$ such that

- ($q_n$ and) $q'$ and $q''$ are identical below $L^*$,
- the nm norms of $q''$ remain at least $n+1$ starting from $4k^*$, and
- $q''$ essentially decides $\dot{F}$, i. e. $q''$ decides $\dot{F}$ below some height $L^{**} = (4k^{**}, 0)$.

Since we already know $\dot{T}$ is read continuously by $p$ (and thus also by any stronger condition), we pick $L^{**}$ large enough such that $q''$ decides $\dot{T}_n \cap 2^{i_{n+1}}$ below $L^{**}$, and also such that the nm norms of $q''$ are at least $n+4$ starting from $4k^{**}$ and $4k^{**} > 2^{i_{n+1}+2}$. This $L^{**} = (4k^{**}, 0)$ will be our $L_{n+1} = (4k_{n+1}, 0)$. So far, we have defined $L_{n+1}$ and $i_{n+1}$ and fulfilled properties (i), (iv) and (v).

Step 3: The set $\dot{F}$ is forced to be a subset of $\dot{T}_n \cap 2^{i_{n+1}}$ of relative size at least $1-\varepsilon$, and both $\dot{F}$ and $\dot{T}_n \cap 2^{i_{n+1}}$ are decided by $q''$ below $L_{n+1}$. We also already know that $\dot{T}_n \cap 2^{i_{n+1}}$ does not depend on any indices in $A_{\mathrm{cn}} \cup A_{\mathrm{slalom}}$ below $L_n$. Hence we can construct a name $\dot{F}' \subseteq \dot{F}$, also not depending on such indices, such that $\dot{F}'$ has relative size at least $1 - \varepsilon \cdot n^P_{<L_n} = 1 - 1/n^P_{<L^*} \geq 1/2$, as follows:

Each $\eta \in \mathrm{poss}(q'', {<}L_{n+1})$ determines objects $F_\eta \subseteq S_\eta$ in the sense that

$$q'' \wedge \eta \Vdash \dot{F} = F_\eta \text{ and } \dot{T}_n \cap 2^{i_{n+1}} = S_\eta.$$

We call two possibilities $\eta, \eta' \in \operatorname{poss}(q'', {<}L_{n+1})$ equivalent if they differ only on indices in $A_{\mathrm{cn}} \cup A_{\mathrm{slalom}}$ below $L_n$. (Note that this implies $S_\eta = S_{\eta'}$.) Obviously, each equivalence class $[\eta]$ has size at most $n^P_{<L_n}$; for each such equivalence class, let $F'_{[\eta]} := \bigcap_{\vartheta \in [\eta]} F_\vartheta$; the relative size of any such $F'_{[\eta]}$ then is at least $1 - \varepsilon \cdot n^P_{<L_n}$. Hence the function mapping each $\eta$ to $F'_{[\eta]}$ defines a name $\dot{F}'$ (not depending on any indices in $A_{\mathrm{cn}} \cup A_{\mathrm{slalom}}$ below $L_n$) for a subset of $\dot{T}_n \cap 2^{i_{n+1}}$ of relative size at least $1/2$.

Since $\dot{T}_n$ is forced to have $\dot{r}_n \in 2^{i_n}$ as its stem and measure greater than $1/2 \cdot 2^{-i_n}$, the size of $\dot{T}_n \cap 2^{i_{n+1}}$ is forced to be greater than $2^{i_{n+1}-(i_n+1)}$, and the size of $\dot{F}'$ is then forced to be greater than $2^{i_{n+1}-(i_n+1)} \cdot 1/2 = 2^{i_{n+1}-i_n-2}$, which is greater than $2^{i_{n+1}}/n^B_{L_n}$ by property (iv).

So far, we have achieved the following: $\dot{T}_n \cap 2^{i_{n+1}}$ and its subset $\dot{F}'$ are decided by $q''$ below $L_{n+1}$ not using any indices in $A_{\mathrm{cn}} \cup A_{\mathrm{slalom}}$ below $L_n$; $q''$ forces each $s \in \dot{F}'$ to fulfil $\lambda([\dot{T}_n] \cap [s]) \geq (1-\varepsilon) \cdot 2^{-i_{n+1}}$; and as a subset of $2^{i_{n+1}}$, $\dot{F}'$ is forced to have measure greater than $1/n^B_{L_n}$.

<u>Step 4:</u> We define the condition $q^* \leq q''$ by replacing all lim sup creatures in $q''$ starting from $L^*$ and below $L_{n+1}$ by arbitrary singletons. So $q^*$ is identical to $q_n$ below $L_n$, and identical to $q''$ starting from $L_{n+1}$. Note that so far, the nm norms of $q^*$ remain at least $n+1$ starting from $4k_n$. In the next few (lengthy) steps, we will define $q_{n+1}$ from $q^*$ by modifying the creatures in $q^*$ starting from $L_n$ and below $L_{n+1}$ such that afterwards, the nm norms of $q_{n+1}$ will remain at least $n$ starting from $4k_n$, and there will be witnesses for lim sup norms at least $n$ between $L_n$ and $L_{n+1}$, as required to fulfil property (ii).

Since $q^*$ decides both $\dot{T}_n \cap 2^{i_{n+1}}$ and $\dot{F}'$ below $L_{n+1}$ not using any indices in $A_{\mathrm{cn}} \cup A_{\mathrm{slalom}}$ below $L_n$, we decompose the set of possibilities $\operatorname{poss}(q^*, {<}L_{n+1})$ into $U \times V \times W$ as follows:

- $U := \operatorname{poss}(q^*, {<}L_n) = \operatorname{poss}(q_n, {<}L_n)$,
- $V$ are the possibilities of $q^*$ starting from $L_n$ and below $L^*$, and
- $W$ are the possibilities of $q^*$ starting from $L^*$ and below $L_{n+1}$, for which we only have to consider the nm part, as the lim sup part has just been defined to be arbitrary singletons.

We will now proceed as follows: For each $\nu \in W$, we will perform an induction on the heights starting from $L_n$ up to $(L^*)^-$ to arrive at a candidate $D(\nu)$ for the creatures of $q_{n+1}$ between $L_n$ and $L^*$; we will then use bigness to see that for many $\nu \in W$, the candidates $D(\nu)$ will be equal, and we will use that fact to finally define $q_{n+1}$.

<u>Step 5:</u> Fix some $\nu \in W$. Recall that relevant heights (in the context of this proof) are those $L \in \mathbf{heights}_{\mathrm{tg}}$ for some tg $\neq$ ct such that there is some $\alpha_L \in \operatorname{supp}(q^*) \cap A_{\mathrm{tg}}$ with a non-trivial $q^*(\alpha_L, L)$. We will inductively go through all heights $L$ with $L_n \leq L < L^*$ (although we will only have to do something for relevant heights) and successively define conditions $q^L \leq q^*$ such that for any $L_n \leq K < L < L^*$

- $q^L \leq q^K$ and $q^K$ and $q^L$ are identical up to (including) $K$,
- the norm of $q^K(\alpha_K, K)$ decreased by at most 1 when compared with the norm of $q^*(\alpha_K, K)$, and

- the norm of $q^K(\alpha_L, L)$ decreased by at most $i/n^B_L$ when compared with the norm of $q^*(\alpha_L, L)$, where $i$ is the number of steps already performed in the induction (i. e. the number of heights between $L_n$ and $K$).

This means that the induction successively strengthens the non-trivial creature at height $L$ until the induction height is $L$ itself; after that step, the non-trivial creature at height $L$ is final and will no longer be modified.

We will also define functions $F^L$ mapping each $\eta \in U \times V$ to a subset $F^L(\eta)$ of $2^{i_{n+1}}$ such that

- $q^{L^-} \wedge (\eta, \nu) \Vdash F^L(\eta) \subseteq \dot{F}'$,
- $F^L(\eta)$ is of relative size at least $1/n^B_L$, and
- $F^L(\eta)$ does not depend on any indices in $A_{\text{cn}} \cup A_{\text{slalom}}$ below $L$.

The preparation for the induction (so that we can start with $L = L_n$) is simply to set $q^{L_n^-} := q^*$ and $F^{L_n} := \dot{F}'$.[40]

Now assume we are at some step $L_n \leq L < L^*$ of the iteration and have already defined $q^{K^-}$ and $F^K$ for all $L_n \leq K < L$. If $L$ is not a relevant height or if the associated index is in $A_{\text{nm}} \cup A_{\text{nn}}$, we do not have to do anything and can set $q^L := q^{L^-}$ and $F^{L^+} := F^L$. So assume the creature $C := q^{L^-}(\alpha_L, L)$ associated with the relevant height $L$ is of type cn or slalom.

We now further decompose $V$ (restricted to just those possibilities which are compatible with $C$) into

- $V^-$, the part below $L$,
- $C$, the part at height $L$, and
- $V^+$, the part strictly above $L$ (and below $L^*$).

Hence we can write every $\eta \in U \times V$ (which is compatible with $C$) as $(\eta^-, \eta^L, \eta^+)$, where $\eta^- \in U \times V^-$, $\eta^L \in C$ and $\eta^+ \in V^+$.

If we now fix $\eta^-$ and $\eta^+$, the function $F^L$ is reduced to an $F^{(\eta^-,\eta^+)}$ mapping each $X \in C$ to a subset of $2^{i_{n+1}}$ of relative size at least $1/n^B_L$. Hence we can use (depending on the type of $C$) either Lemma 11.4 or Lemma 12.1 to strengthen the creature $C$ to $D(\eta^-, \eta^+)$, decreasing the norm by at most $1/n^B_L$, such that

$$F^*(\eta^-, \eta^+) := \bigcap_{X \in D(\eta^-,\eta^+)} F^{(\eta^-,\eta^+)}(X)$$

is a set of relative size at least $1/n^B_{L^+}$.

If we now fix only $\eta^+$ and successively iterate this strengthening for all $\eta^- \in U \times V^-$, we ultimately arrive at some $D(\eta^+) \subseteq C$ with the norm decreasing by at most $n^P_{<L}/n^B_L < 1$ in total. Note that since $2^{n^S_L} < n^B_{L^+}$, there are fewer than $n^B_{L^+}$ many possible values for $D(\eta^+)$ and we can apply strong bigness in the form of Lemma 5.8 on the $V^+$ part to strengthen all $q^{L^-}(\alpha_K, K)$ for $L^+ \leq K < L^*$ to $q^L(K)$, decreasing the norm by at most $1/n^B_K$ at each height $K$, such that for each $\eta^+$ in the resulting smaller $\overline{V}^+$, we get the same $D := D(\eta^+)$.[41] This $D$ then will be the (final) value of $q^L(\alpha_L, L)$. If we now define

$$F^{L^+}(\eta) := \bigcap_{X \in D} F^L(\eta^-, X, \eta^+),$$

by the considerations above, this is a set of relative size at least $1/n^B_{L^+}$, does not depend on any indices in $A_{\mathrm{cn}} \cup A_{\mathrm{slalom}}$ below $L^+$, and is forced to be a subset of $\dot{F}'$ by $q^L \wedge (\eta, \nu)$.

[40] We ask the reader to excuse the abuse of notation here; a name and a function are, of course, not the same thing, but for all practical purposes, they might as well be in the context of this step of the proof.

Having now defined $q^L$ and $F^{L^+}$, we can proceed with the next step of the inductive construction.

Step 6: We perform the construction in Step 5 independently for each $\nu \in W$ (i. e. starting with the original $q^*$ each time). We thus get a (potentially) different $q_\nu^{(L^*)^-}$ for each $\nu$. Since the number of possible values for $q_\nu^{(L^*)^-}$ is less than $n^B_{L^*}$, we can now apply Lemma 5.8 again to thin out the creatures $q^*(\alpha_K, K)$ for $L^* \leq K < L_{n+1}$ to $q_*(\alpha_K, K)$, decreasing the norm by at most $1/n^B_K$ at each height $K$, such that for each $\nu$ in the resulting smaller $\overline{W}$, we get the same $q_{**} := q_\nu^{(L^*)^-}$. We can then finally define $q_{n+1} := q_{**}^{<L^*} {}^\frown q_*^{\geq L^*}$, which fulfils property (ii) by construction.

Step 7: Now, this $q_{n+1}$ forces the family of "terminal" $F_\nu^{(L^*)^-}$ (for $\nu \in \overline{W}$) to constitute a name $\dot{F}''$ for a subset of $\dot{F}' \subseteq 2^{i_{n+1}}$ of relative size greater than 0, and $q_{n+1}$ decides $\dot{F}''$ below $L_{n+1}$ not using any indices in $A_{\mathrm{cn}} \cup A_{\mathrm{slalom}}$ – due to the fact that below $L_n$, even the name $\dot{F}'$ did not depend on such indices; from $L_n$ up to $L^*$, we removed the dependence on such creatures height by height in Step 5; and from $L^*$ up to $L_{n+1}$, by Step 4 only singletons remain for such lim sup creatures, anyway.

Hence we can pick some name $\dot{r}_{n+1}$ for an arbitrary fixed element of $\dot{F}''$ (e. g. the first element in the natural lexicographic order), and this name fulfils properties (iii) (by construction) and (vi) (since $\dot{r}_{n+1}$ is a node in $\dot{T}_n$, whose stem is forced to be $\dot{r}_n$ by $q_n$).

Step 8: Since $q_{n+1}$ forces $\dot{r}_{n+1} \in \dot{F}''$, $\dot{r}_{n+1}$ is a fat node, which means $\dot{T}' := \dot{T}_n \cap [\dot{r}_{n+1}]$ is forced to have measure greater than $\frac{1-\varepsilon}{2^{i_{n+1}}}$. The tree $\dot{T}'$ is read continuously by $q_n$ and hence also by $q_{n+1}$; in particular, for each $j > i_{n+1}$, the finite initial tree $\dot{T}' \cap 2^j$ is decided below some $L_j$. For each $\eta \in \mathrm{poss}(q_{n+1}, {<}L_j)$, let $T^j_\eta$ be the corresponding value of $\dot{T}' \cap 2^j$ (which is a subset of $2^j$ with at least $2^j \cdot \frac{1-\varepsilon}{2^{i_{n+1}}}$ many elements). It is clear that for $j < j'$ and $\eta \in \mathrm{poss}(q_{n+1}, {<}L_j)$, $\eta' \in \mathrm{poss}(q_{n+1}, {<}L_{j'})$ such that $\eta \subseteq \eta'$, it is forced that the corresponding finite trees are also nested, i. e. $T^{j'}_{\eta'} \subseteq T^j_\eta$.

We now implement a reduction similar to Step 3 to eliminate the dependency on indices in $A_{\mathrm{cn}} \cup A_{\mathrm{slalom}}$: We call two possibilities $\eta, \eta' \in \mathrm{poss}(q_{n+1}, {<}L_j)$ equivalent if they differ only on indices in $A_{\mathrm{cn}} \cup A_{\mathrm{slalom}}$ below $L_{n+1}$. Since from $L^*$ up to $L_{n+1}$, there are only singletons for such lim sup creatures, each equivalence class $[\eta]$ has size at most $n^P_{<L^*}$. For each such equivalence class, let $T^j_{[\eta]} := \bigcap_{\vartheta \in [\eta]} T^j_\vartheta$. Note that by the nesting of the $T^j_\vartheta$, the $T^j_{[\eta]}$ are also nested (for $j, j', \eta, \eta'$ as above), and the size of $T^j_{[\eta]}$ is at least

$$\frac{2^j \cdot (1 - \varepsilon \cdot n^P_{<L^*})}{2^{i_{n+1}}} = \frac{2^j \cdot (1 - {}^1\!/n^P_{<L_n})}{2^{i_{n+1}}}.$$

[41] Keep in mind that since we are working in $\mathbb{Q}_{\text{non-ct}}$, there will be no $K \in \mathbf{heights}_{\text{ct}}$, and hence we can apply Lemma 5.8.

So $q_{n+1}$ forces the family of such $T^j_{[\eta]}$ (for $j > i_{n+1}$ and $\eta \in \operatorname{poss}(q_{n+1}, {<}L_j)$) to constitute a name $\dot{T}_{n+1}$ as required to fulfil property (vii). ☐

**Corollary 12.4.** $\mathbb{Q}$ *forces* $\operatorname{non}(\mathcal{N}) \leq \kappa_{\text{nn}}$.

*Proof.* Fix a condition $p$ and a sequence of names of null sets $\langle \dot{N}_i \mid i \in I \rangle$ which $p$ forces to be a basis of null sets. As described above, for each $i \in I$, we can assume that $\dot{N}_i = N_{\dot{T}_i}$ for some name for a sturdy tree $\dot{T}_i$. Let $\dot{R}$ consist of all reals read continuously only using indices in $A_{\text{nm}} \cup A_{\text{nn}}$.

By the preceding lemma, for each $\dot{T}_i$, there is a $q \leq p$ in $\mathbb{Q}_{\text{non-ct}}$ and an $\dot{r} \in \dot{R}$ such that $q \Vdash \dot{r} \in [\dot{T}_i]$ and hence $\mathbb{Q}_{\text{non-ct}} \Vdash \dot{r} \notin N_{\dot{T}_i}$; it follows that $\mathbb{Q}_{\text{non-ct}} \Vdash$ "$\dot{R}$ is not null" and hence also $\mathbb{Q} \Vdash$ "$\dot{R}$ is not null". ☐

This proves (M3) of Theorem 1.1, and hence completes the proof of that theorem entirely.

## 13. Failed Attempts, Limitations and Open Questions

To counteract the common habit of only talking about successes and withholding the failed attempts that went before, we want to give a very brief account of two results we attempted, but failed to achieve in the course of writing this paper.

For one, we wanted to add $\kappa_{\text{rp}}$ many factors which would carefully increase the cardinals $\mathfrak{r}$ and $\mathfrak{u}$ to $\kappa_{\text{rp}}$, a cardinal between $\kappa_{\text{cn}}$ and $\kappa_{\text{ct}}$. The plan was to use a forcing poset $\mathbb{Q}_{\text{rp}}$ (a variant of the forcing poset from [GS90]) in each factor. While it seemed quite simple to align the structure of $\mathbb{Q}_{\text{rp}}^{\kappa_{\text{rp}}}$ with the structure of $\mathbb{Q}_{\text{ct},\kappa_{\text{ct}}}$ to allow the proof of section 8 to function for both $\mathbb{Q}_{\text{rp}}^{\kappa_{\text{rp}}}$ and $\mathbb{Q}_{\text{ct},\kappa_{\text{ct}}}$, it was not clear why the $\kappa_{\text{ct}}$ many Sacks-like reals would preserve $\mathfrak{u} \leq \kappa_{\text{rp}}$, or indeed why the old reals would be unreapable even after multiplying the forcing poset with the product of merely two copies of $\mathbb{Q}_{\text{rp}}$.

The second idea we had was to add Cohen forcing to the construction to control the value of $\operatorname{cov}(\mathcal{M})$ in the resulting model. This would have complicated a lot of the proofs, since many things would then have turned into names dependent on the Cohen-generic filter; however, a more fundamental problem is that this approach destroys the Sacks property of the "upper" part of the construction:

**Lemma 13.1.** *Let $\mathbb{C}$ be the Cohen forcing poset and let $\mathbb{S}$ be the Sacks forcing poset. Then $V^{\mathbb{C}\times\mathbb{S}}$ does not have the Sacks property over $V^{\mathbb{C}}$.*

*More generally, consider two forcing posets $\mathbb{X}$ and $\mathbb{Y}$, where $\mathbb{X}$ adds an unbounded real $\dot{x}$ and $\mathbb{Y}$ adds another new real $\dot{y}$. Then $V^{\mathbb{X}\times\mathbb{Y}}$ does not have the Sacks property over $V^{\mathbb{X}}$.*

*Proof.* We prove the stronger claim. Let $\langle\dot\tau(n)\rangle := \langle\dot y|_{\dot x(n)}\rangle$, i. e. each $\dot\tau(n)$ is an initial segment of the new real $\dot y$.

Assume that we have some sequence of $\mathbb{X}$-names $\dot B_k$ for a $(k+1)$-slalom catching $\dot\tau(k)$:

$$\Vdash_{\mathbb{X}} |\dot B_k| = k+1$$
$$(p,q) \Vdash \forall\, k<\omega\colon \dot B_k \subseteq 2^{\dot x(k)} \wedge \dot\tau(k) \in \dot B_k$$

Let $n$ be the index of the first value of $\dot x$ not bounded by $p$. Let $T$ be the tree of potential initial segments of $\dot y$ below the condition $q$, i. e. $T := \{s \mid q \nVdash s \not\subseteq \dot y\}$; since $\dot y$ is a new real, $T$ must have unbounded width. Hence there is some $m$ such that the $\dot\tau(m)$ has at least $n+2$ many possible values. Fix such an $m$. Then find $p' \leq p$ forcing $\dot x(n) = m^* \geq m$, and let $p^* \leq p'$ be such that $p^*$ decides $\dot B_n$, i. e. $p^* \Vdash \dot B_n = B$ for some $B$.

Since $\dot\tau(m^*)$ has at least $n+2$ many possible values, there is some possible value $v$ that is not in $B$. But then there is a $q^* \leq q$ forcing $\dot\tau(m^*) = v$, and hence

$$(p^*,q^*) \Vdash \dot\tau(m^*) \notin B,$$

which is a contradiction. □

The point of this lemma is that if we do not have the Sacks property, $\operatorname{cof}(\mathcal{N})$ will increase.

We turn our attention towards related work and open questions. Several recent results [FFMM18, BCM21, KTT18, GKS19, KST19, GKMS19, GKMS2x, GKMS20] have constructed models in which eight or even all ten conceivably different cardinal characteristics in Cichoń's diagram take different values. The constructions involved are all finite support iterations, however, which necessarily means the left side of Cichoń's diagram must be less than or equal to the right side, in particular $\operatorname{non}(\mathcal{M}) \leq \operatorname{cov}(\mathcal{M})$ (since the cofinality of the iteration length lies between these two cardinal characteristics). In contrast, [FGKS17] and our improvement thereof have $\operatorname{non}(\mathcal{M}) > \operatorname{cov}(\mathcal{M})$.

However, as far as Cichoń's diagram is concerned, our creature forcing construction still has rather strict limitations as explained in the preceding section. Necessarily, $\mathfrak{d} = \aleph_1$ by the $\omega^\omega$-boundedness of the forcing posets involved; the only open question regarding Cichoń's diagram and our construction is whether it is possible to separate $\operatorname{cov}(\mathcal{N})$ from $\aleph_1$.

**Question A.** *Is it possible to modify the construction to achieve $\aleph_1 < \operatorname{cov}(\mathcal{N})$?*

Finally, our failed attempt to introduce $\mathfrak{r}$ and $\mathfrak{u}$ into the construction motivates the following general question:

**Question B.** *Are there any well-known cardinal characteristics which can be set via a* $\limsup$*-type creature forcing poset compatible with the structure of* $\mathbb{Q}$*?*

Institute of Discrete Mathematics and Geometry, TU Wien, Wiedner Hauptstrasse 8–10/104, 1040 Wien, Austria

*Email address*: martin.goldstern@tuwien.ac.at

*URL*: http://www.tuwien.ac.at/goldstern/

Institute of Discrete Mathematics and Geometry, TU Wien, Wiedner Hauptstrasse 8–10/104, 1040 Wien, Austria

*Email address*: mail@l17r.eu

*URL*: https://l17r.eu